\documentclass[a4paper]{LMCS}

\newif\iflmcs\lmcstrue % <--------------- TOGGLE THIS FOR FINAL VERSION OF LMCS
\usepackage{stmaryrd}
\usepackage[lutzsyntax]{virginialake}\vlsmallbrackets\aftrianglefalse
\usepackage{hyperref}

\def\doi{4 (1:9) 2008}
\lmcsheading%
{\doi}
{1--36}
{}
{}
{Aug.~\phantom{0}2, 2007}
{Mar.~31, 2008}
{}   

\begin{document}

\title[Normalisation Control in Deep Inference   via Atomic Flows]
      {Normalisation Control in Deep Inference\\ via Atomic
      Flows\rsuper *}

\author{Alessio Guglielmi}
\address{University of Bath, Bath BA2 7AY, UK}
\iflmcs\email{\{A.Guglielmi,T.E.Gundersen\}@Bath.Ac.UK}\fi

\author{Tom Gundersen}
%\iflmcs\address{University of Bath, Bath BA2 7AY, UK}\fi
%\iflmcs\email{T.E.Gundersen@Bath.Ac.UK}\fi

\thanks{This work was in part funded by an Overseas Research Scholarship and a Research Studentship, both from the University of Bath, and by the British Council Alliance Programme.}


% Mathematical Logic---Proof theory}

\begin{abstract}
We introduce `atomic flows': they are graphs obtained from derivations by tracing atom occurrences and forgetting the logical structure. We study simple manipulations of atomic flows that correspond to complex reductions on derivations. This allows us to prove, for propositional logic, a new and very general normalisation theorem, which contains cut elimination as a special case. We operate in deep inference, which is more general than other syntactic paradigms, and where normalisation is more difficult to control. We argue that atomic flows are a significant technical advance for normalisation theory, because 1) the technique they support is largely independent of syntax; 2) indeed, it is largely independent of logical inference rules; 3) they constitute a powerful geometric formalism, which is more intuitive than syntax.
\end{abstract}

\maketitle

%===============================================================================
\section{Introduction}

We are interested in normalising derivations in proof systems of propositional logic. As for natural deduction and the sequent calculus \cite{Gent:69:Investig:xi}, we intend normalisation as eliminating cuts, or, more in general, `detours', from derivations. Normalisation is performed by algorithms that, given a non-normal derivation, produce a normal derivation, free of detours. A typical detour can be depicted, in an abstract representation of a portion of a derivation, as on the left in
\[
\atomicflow{
( 0  , 4  )*{\afaidnw{}{}};
( 6  , 3  )*{\afvjd6{}a};
( 0  , 2.3)*{\aflabelright{\bar a}};
( 2  , 2  )*{\afvj4};
(-2  , 1  )*{\afvju6a{}};
( 4  ,-2  )*{\afaiunw{}{}};
(-3.5, 0  )*{\invisiblemark};
( 7.5, 0  )*{\invisiblemark}}
\quad\to\quad
\atomicflow{
(0  ,4  )*{\afvj8};
(0  ,4.3)*{\aflabelright a};
(1.5,0  )*{\invisiblemark}}
\quad,
\]
where the atom $\bar a$ is created and destroyed, respectively, by an axiom and a cut, represented as two horizontal bars. In many cases, the diagram on the left can be streamlined as on the right, and this is what happens in a typical normalisation step. However, inside concrete proof systems, this abstract, geometric simplicity is almost always severely obscured by syntactic bureaucracy.

Originally stimulated by Girard's \cite{Gira:87:Linear-L:wm}, we share with several colleagues the research objective of getting rid of bureaucracy in proof systems. We believe that, ultimately, derivations are geometric objects of some sort, for which current syntactic formalisms only offer imperfect, bureaucratic representations. Our hopes rely, in part, on the success of proof nets for multiplicative linear logic. Those proof nets are a geometric, largely bureaucracy-free proof system, where normalisation can be directly performed as depicted above. 

In classical propositional logic, the situation is more complicated. There are proof nets for classical propositional logic, but they are not a proof system (as defined by Cook and Reckhow in \cite{CookReck:79:The-Rela:mf}), because checking their correctness requires exponential time in their size. There is not much hope for improvement, because obtaining a proof system for classical propositional logic based on proof nets would imply that $\mathsf{coNP}$ is equal to $\mathsf{NP}$. So, proof nets are a source of inspiration, and especially those devised by Lamarche and Stra\ss burg-\break er in \cite{LamaStra:05:Naming-P:ov}, but they are not a solution to our quest of bureaucracy-free formalisms for logics as expressive as classical logic.

This paper is about a geometric, bureaucracy-free formalism, called \emph{atomic flows}. An atomic flow is a directed graph obtained from a derivation by only retaining information about the creation and destruction of atom occurrences. Atomic flows are, essentially, specialised Buss flow graphs \cite{Buss:91:The-Unde:uq}. Notably, the atomic flow of a derivation completely disregards all the logical relations and associated inference steps; so, an atomic flow is not a derivation, but only a very abstract representation of it. Since atomic flows and the techniques they induce are largely independent of syntax, we think that they will help us defining a bureaucracy-free formalism. In fact, they yield a geometric and bureaucracy-free understanding of normalisation, which is the most important proof-theoretic aspect of proof systems.

We show that the information contained in an atomic flow is sufficient to control several normalisation algorithms for its associated derivation. This means that a normalisation algorithm extracts from a derivation its atomic flow, and then decides the normalisation steps only based on the structure of the atomic flow. In particular, atomic flows provide convenient induction measures for termination. The advantage of atomic flows is their simplicity, compared to derivations, as Figure~\ref{FigExStrRed} (page~\pageref{FigExStrRed}) eloquently shows. They allow us to prove, in this paper, a new and more general normalisation result than cut elimination. Proving the same without atomic flows is conceivable, but perverse; finding the normalisation algorithm without atomic flows is, in our opinion, inconceivable. 

\newcommand{\SKS}{\mathsf{SKS}}
We consider derivations in the proof system $\SKS$ of the calculus of structures \cite{BrunTiu:01:A-Local-:mz}, which is a very general formalism based on deep inference \cite{Gugl:06:A-System:kl}. Because of their comparatively larger expressive power, normalisation in pure, unconstrained proof systems in the calculus of structures is much more challenging than in more disciplined formalisms, like the sequent calculus or natural deduction. In fact, deep-inference rules disregard the notion of root connective of a formula, which is a crucial asset for normalisation in non-deep inference formalisms.

We show that every $\SKS$ derivation can be \emph{streamlined}, \emph{i.e.}, we can remove all causal dependencies between axioms and cuts of the kind depicted above on the left. This result generalises cut elimination for two reasons: 1) $\SKS$ can faithfully embed derivations in most other proof systems of different formalisms, like the sequent calculus, hypersequents, natural deduction, resolution and others; 2) a cut-free proof is a special case of a streamlined derivation. Intuitively, this result is sort of a Craig's interpolation for derivations instead of formulae. We note that, contrary to traditional cut elimination, our result is \emph{symmetric} along the axis premiss-conclusion of a derivation (as is traditional in deep inference).

The core idea of our algorithms, as controlled by atomic flows, consists in slightly altering a derivation around a couple of matching axiom and cut, and composing it with itself. In principle, this is very similar to normalisation in natural deduction. This idea is due to Tiu and then used by Br\"unnler, in \cite{Brun:03:Atomic-C:oz}, to simplify his previous proof of cut elimination in the calculus of structures. In our case, the idea allows us to dispense entirely with the usual case analysis and permutation of rules, which is the standard routine in cut elimination proofs in the sequent calculus and elsewhere.

The observed independence from syntax leads us to believe that successful normalisation stems less than is usually believed from the mutual `harmony' between logical rules. In fact, in this paper, there is no concern whatsoever for this issue. Rather, the choice in designing logical rules seems to be essentially free, provided that they are linear and that they support atomic structural rules. Our technique relies on substituting formulae for occurrences of atoms, which is only manageable if the structural part of the proof system is atomic. As a matter of fact, complete and analytic proof systems with atomic structural rules and linear logical ones can only be designed in deep inference, for most logics.

Finally, we note that several normalisation algorithms can be designed and several normal forms can be obtained by employing atomic flows. We are seeing a \emph{robust} normalisation phenomenon: variations in the algorithms are possible and interesting, we are not dealing here with a delicate property of derivations that requires extra care and attention to the tiniest syntactic detail. Once the basics of the technique are mastered, derivations can be manipulated with ease. We think that atomic flows get closer than ever to the essence of normalisation in classical propositional logic.

After a brief introduction to deep inference, we introduce atomic flows and their reductions, and then the streamlining algorithms. The two crucial, simple ideas to understand in this paper are shown in Remark~\ref{RemCycle} (page~\pageref{RemCycle}) and in Definition~\ref{DefRedS} (page~\pageref{DefRedS}).

%===============================================================================
\section{Background on Deep Inference}

Deep inference is a relatively recent development in proof theory. It is a methodology according to which several formalisms can be defined with excellent structural properties. The calculus of structures \cite{Gugl:06:A-System:kl} is one of them and is now well developed for classical \cite{Brun:03:Atomic-C:oz,Brun:06:Cut-Elim:cq,Brun:06:Locality:zh,BrunTiu:01:A-Local-:mz,Brun:06:Deep-Inf:qy}, intuitionistic \cite{Tiu:06:A-Local-:gf}, linear \cite{Stra:02:A-Local-:ul,Stra:03:MELL-in-:oy}, modal \cite{Brun::Deep-Seq:ay,GoreTiu:06:Classica:uq,Stou:06:A-Deep-I:rt} and commutative/non-commutative logics \cite{Gugl:06:A-System:kl,Tiu:06:A-System:ai,Stra:03:Linear-L:lp,Brus:02:A-Purely:wd,Di-G:04:Structur:wy,GuglStra:01:Non-comm:rp,GuglStra:02:A-Non-co:lq,GuglStra:02:A-Non-co:dq,Kahr:06:Reducing:hc,Kahr:07:System-B:fk}. The basic proof complexity properties of the calculus of structures are known \cite{BrusGugl:07:On-the-P:fk}. The calculus of structures promoted the discovery of a new class of proof nets for classical and linear logic \cite{LamaStra:05:Construc:qq,LamaStra:05:Naming-P:ov,LamaStra:06:From-Pro:et,StraLama:04:On-Proof:ec} (see also \cite{Guir:06:The-Thre:qt}). There exist implementations in Maude of deep-inference proof systems \cite{Kahr:07:Maude-as:lr}. For a better introduction than this, we refer the reader to \cite{Brun:03:Atomic-C:oz}.

\newcommand{\fff}{\mathsf f}
\newcommand{\ttt}{\mathsf t}
\newcommand{\ot}{\mathbin\shortleftarrow}

%---------------------------------------
\begin{defi}
\emph{Formulae}, $\alpha$, $\beta$, $\gamma$, $\delta$ are freely built from: \emph{units}, $\fff$ (false), $\ttt$ (true); \emph{atoms}, $a$, $b$, $c$, $d$, $e$; \emph{disjunction} and \emph{conjunction}, ${\vlsbr[\alpha.\beta]}$ and $\vlsbr(\alpha.\beta)$. On the set of atoms a (non-identical) involution $\bar\cdot$ is defined, and dual atom occurrences, as $a$ and $\bar a$, can appear in formulae. We denote \emph{contexts}, \emph{i.e.}, formulae with a hole, by $\xi\vlhole$ and $\zeta\vlhole$; we also use \emph{multiple} contexts, $\xi\vlhole\cdots\vlhole$, \emph{i.e.}, formulae with many holes; for example, if $\xi\{a\}$ is $\vls(b.[a.c])$, then $\xi\vlhole$ is $\vls(b.[\vlhole.c])$, $\xi\{b\}$ is $\vls(b.[b.c])$ and $\xi\vlscn(a.d)$ is $\vls(b.[(a.d).c])$; if $\xi\{a\}\{b\}\{c\}$ is $\vls(b.[(a.d).c])$ then $\xi\{b\}\{c\}\{a\}$ is $\vls(c.[(b.d).a])$.
\end{defi}

%---------------------------------------
\begin{rem}
Negation is only defined for atoms, which is not a limitation thanks to De Morgan laws.
\end{rem}

Note that when we write $\xi\{a\}$, we mean that an occurrence of $a$ exists in the formula, we singled it out and we refer specifically to that occurrence. It is important to distinguish between an atom $a$ and a set of occurrences of atom $a$ inside a formula or a derivation. In the following, we mark in various ways occurrences of atoms, and we perform several substitutions of formulae in the place of atom occurrences.

\newcommand{\one}{{\mathchoice{\scriptstyle\mathbf1}
                              {\scriptstyle\mathbf1}
                              {\scriptstyle\mathbf1}
                              {\scriptscriptstyle\mathbf1}}}
\newcommand{\two}{{\mathchoice{\scriptstyle\mathbf2}
                              {\scriptstyle\mathbf2}
                              {\scriptstyle\mathbf2}
                              {\scriptscriptstyle\mathbf2}}}
\newcommand{\mk}[1]{{#1}^{\scriptscriptstyle\bullet}}
%---------------------------------------
\begin{defi}
\emph{Inference rules}, $\rho$, have one \emph{premiss} and one \emph{conclusion}, and their \emph{instances} are used in \emph{inference steps} to rewrite inside formulae. A \emph{derivation}, $\Phi$, from $\alpha$ (\emph{premiss}) to $\beta$ (\emph{conclusion}) is a chain of inference steps with $\alpha$ at the top and $\beta$ at the bottom, and is usually indicated by $\vlder{\Phi}{\mathcal S}{\beta}{\alpha}$, where $\mathcal S$ is the name of the deductive system or a set of inference rules; a \emph{proof} is a derivation from $\ttt$; besides $\Phi$, we denote derivations with $\Psi$. We denote with $\xi\{\Phi\}$ the result of including every formula of $\Phi$ into the context $\xi\vlhole$: since we adopt deep inference, $\xi\{\Phi\}$ is a valid derivation. We denote with $\Phi\{a\ot\alpha\}$, $\Phi\{\fff\ot\alpha\}$ and $\Phi\{\ttt\ot \alpha\}$ the operation of substituting $\alpha$ into a set of \emph{occurrences} of an atom $a$ or unit in $\Phi$; the result is not necessarily a valid derivation, because some instances of rules might break; which occurrences to replace is always made clear by suitable decorations of $a$, $\fff$ and $\ttt$, like $a^\one$ and $\mk\fff$.
\end{defi}

\newcommand{\KS}{\mathsf{KS}}
Now we define the two standard deductive systems for classical propositional logic in deep inference that are used throughout the paper. $\KS$ is analytic, in the sense that premisses only contain subformulae of conclusions, and $\SKS$ is not \cite{Brun:03:Atomic-C:oz,Brun:06:Cut-Elim:cq,Brun:06:Locality:zh,BrunTiu:01:A-Local-:mz}.

\newcommand{\ai}{\mathsf{ai}}
\newcommand{\aw}{\mathsf{aw}}
\newcommand{\ac}{\mathsf{ac}}
\newcommand{\aid}{{\ai{\downarrow}}}
\newcommand{\awd}{{\aw{\downarrow}}}
\newcommand{\acd}{{\ac{\downarrow}}}
\newcommand{\aiu}{{\ai{\uparrow}}}
\newcommand{\awu}{{\aw{\uparrow}}}
\newcommand{\acu}{{\ac{\uparrow}}}
\newcommand{\swi}{\mathsf{s}}
\newcommand{\med}{\mathsf{m}}
%---------------------------------------
\begin{defi}
System $\SKS$ in the calculus of structures is defined by the following \emph{structural} rules:
\[
\begin{array}{@{}c@{}c@{}c@{}}
      \vlinf{\aid}{}{\vls[a.{\bar a}]}{\ttt}&
\qquad\vlinf{\awd}{}a\fff&
\qquad\vlinf{\acd}{}a{\vls[a.a]}\\
\noalign{\smallskip}
      \emph{interaction}&
\qquad\emph{weakening}&
\qquad\emph{contraction}\\
\noalign{\bigskip}
      \vlinf{\aiu}{}\fff{\vls(a.{\bar a})}&
\qquad\vlinf{\awu}{}\ttt a&
\qquad\vlinf{\acu}{}{\vls (a.a)}a\\
\noalign{\smallskip}
      \emph{cointeraction}&
\qquad\emph{coweakening}&
\qquad\emph{cocontraction}\\
\end{array}\quad,
\]
and by the two \emph{logical} rules:
\[
\begin{array}{@{}c@{}c@{}}
\vlinf{\swi}{}{\vls[(\alpha.\beta).\gamma]}{\vls(\alpha.[\beta.\gamma])}&\qquad
\vlinf{\med}{}{\vls([\alpha.\gamma].[\beta.\delta])}
              {\vls[(\alpha.\beta).(\gamma.\delta)]}\\
\noalign{\smallskip}
\emph{switch}&\qquad\emph{medial}\\
\end{array}\quad.
\]
The rule cointeraction is also called an (\emph{atomic}) \emph{cut}. In addition to the rules shown, there is a rule $\vldownsmash{\vlinf={}\delta\gamma}$, such that $\gamma$ and $\delta$ are opposite sides in one of the following equations:
\vlstore{
\vls[\alpha.\beta]         &=\vls[\beta.\alpha]         \quad,&
\vls[\alpha.\fff]          &=\vls[\alpha]               \quad,\\
\vls(\alpha.\beta)         &=\vls(\beta.\alpha)         \quad,&
\vls(\alpha.\ttt)          &=\vls(\alpha)               \quad,\\
\vls[[\alpha.\beta].\gamma]&=\vls[\alpha.[\beta.\gamma]]\quad,&
\vls[\ttt.\ttt]            &=\vls[\ttt]                 \quad,\\
\vls((\alpha.\beta).\gamma)&=\vls(\alpha.(\beta.\gamma))\quad,&
\vls(\fff.\fff)            &=\vls(\fff)                 \quad\vldot}
\begin{align*}
\vlread
\end{align*}
We do not always show the instances of rule $=$, and when we do show them, we gather several contiguous instances into one. System $\KS$ is the same as $\SKS$, but without the rules $\aiu$, $\awu$ and $\acu$. A \emph{cut-free} derivation is a derivation where $\aiu$ is not used. All derivations in this paper are in $\SKS$, unless indicated otherwise.
\end{defi}

Note that all the structural rules only apply to atoms. As shown later, equivalent structural rules applying to formulae instead of atoms can be derived from the atomic ones together with the logical rules. The fact that we can work only with atomic structural rules is essential later on.

Instead of the term `axiom' we use `interaction'; the reason is that, in deep inference, axioms do not close derivation branches. However, it is not misleading to think of interaction instances as axiom instances in the sequent calculus. In several papers, including \cite{Brun:03:Atomic-C:oz}, the reader can find explanations of how reducing a proof in $\SKS$ to a proof in $\KS$ is a cut-elimination process in the traditional sense. In other words, the rules $\aiu$, $\awu$ and $\acu$ are, together, morally equivalent to a cut in the sequent calculus.

There are many $\SKS$ derivations in this paper, providing examples for the above definitions. Mastering the following, standard constructions of the calculus of structures is crucial: 1) moving a formula outside of a context, as in $\vldownsmash{\vlder{}{\{\swi\}}{\vls[\xi\{\fff\}.\alpha]}{\xi\{\alpha\}}}$; 2) bringing a formula inside a context, as in $\vlder{}{\{\swi\}}{\xi\{\alpha\}}{\vls(\xi\{\ttt\}.\alpha)}$. In the next remark, we appeal to both. The constructions in the rest of this subsection are needed later in the paper.

\newcommand{\sus}{\mathsf{ss}}
\newcommand{\contr}{\mathsf{c}}
\newcommand{\cod}{{\contr{\downarrow}}}
\newcommand{\cou}{{\contr{\uparrow}}}
%---------------------------------------
\begin{rem}\label{RemSupSwitch}
For any $\xi\vlhole$, $\zeta\vlhole$ and $\alpha$, by working inductively on the contexts $\xi\vlhole$ and $\zeta\vlhole$, we can build
\[
\vlder{}{\{\swi\}}{\vls[\xi\{\alpha\}.\zeta\{\fff  \}]}
                  {\vls(\xi\{\ttt  \}.\zeta\{\alpha\})}
\quad.
\]
We can do this according to the following two schemes:
\[
\vlderivation                                            {                      
\vlde{}{\{\swi\}}{\vls[\xi\{\alpha\}.\zeta\{\fff  \}]}  {                      
\vlde{}{\{\swi\}}{\vls\zeta\{ \xi\{\alpha\}        \}} {                      
\vlde{}{\{\swi\}}{\vls\zeta\{(\xi\{\ttt  \}.\alpha)\}}{
\vlhy            {\vls(\xi\{\ttt  \}.\zeta\{\alpha\})}}}}}
\qquad\hbox{and}\qquad
\vlderivation                                            {                      
\vlde{}{\{\swi\}}{\vls[\xi\{\alpha\}.\zeta\{\fff  \}]}  {                      
\vlde{}{\{\swi\}}{\vls\xi\{[\alpha.\zeta\{\fff  \}]\}} {                      
\vlde{}{\{\swi\}}{\vls\xi\{\zeta\{\alpha\}         \}}{
\vlhy            {\vls(\xi\{\ttt  \}.\zeta\{\alpha\})}}}}}
\quad.
\]
For example, for $\xi\vlhole=\vls([\vlhole.b].c)$ and $\zeta\vlhole=\vls[(d.\vlhole).e]$, consider
\[
\vlderivation                                          {
\vlin{=   }{}{\vls[([\alpha.b].c).[(d.\fff).e]]}      {
\vlin{\swi}{}{\vls[[(d.\fff).([\alpha.b].c)].e]}     {
\vlin{=   }{}{\vls[(d.[\fff.([\alpha.b].c)]).e]}    {
\vlin{\swi}{}{\vls[(([(\alpha.\ttt).b].c).d).e]}   {
\vlin{=   }{}{\vls[(((\alpha.[\ttt.b]).c).d).e]}  {
\vlin{\swi}{}{\vls[(([\ttt.b].c).(\alpha.d)).e]} {
\vlin{=}   {}{\vls(([\ttt.b].c).[(\alpha.d).e])}{
\vlhy        {\vls(([\ttt.b].c).[(d.\alpha).e])}}}}}}}}}
\quad.
\]
We define the following `macro' rule $\sus$, called \emph{super switch}, to be a shorthand for any derivation of the above form:
\[
\vlinf{\sus}{}{\vls[\xi\{\alpha\}.\zeta\{\fff\}]}
              {\vls(\xi\{\ttt\}.\zeta\{\alpha\})}
\quad.
\]
\end{rem}

%---------------------------------------
\begin{rem}\label{RemGenContr}
$\vlupsmash{\vlinf{\cod}{}\alpha{\vls[\alpha.\alpha]}}$ and $\vlupsmash{\vlinf{\cou}{}{\vls(\alpha.\alpha)}\alpha}$ are two other `macro' rules: they are called, respectively, \emph{contraction} and \emph{cocontraction}, and they apply to generic formulae instead of atoms. They can be derived, respectively, from $\{\acd,\med\}$ and $\{\acu,\med\}$. For an example of the latter, see Figure~\ref{FigExAF} (page~\pageref{FigExAF}).
\end{rem}

%===============================================================================
\section{Atomic Flows}

In this section, we define atomic flows and $\ai$-cycles inside them. These special cycles are circular dependencies between interactions and cointeractions, and they are particularly complex to deal with in the syntax. They turn out to be natural concepts in atomic flows, and understanding them is the key to our normalisation algorithm. The section closes with the definition of streamlined atomic flows and associated derivations.

%-------------------------------------------------------------------------------
\subsection{Atomic Flows and Derivations}

Atomic flows are somewhat similar to proof nets. However, we prove that, no matter how we freely build an atomic flow (as opposed to a proof net structure), the flow is associated with some derivation. So, atomic flows are always `sequentialisable', in proof-net parlance. In fact, atomic flows carry much less information than derivations do, because they do not keep track of the logical relations between the atoms they trace, only their structural information is retained (in the sense of structural rules, as opposed to logical ones).

We can think of atomic flows as composite diagrams that are freely generated from a set of six elementary diagrams. Technically, atomic flows are special kinds of labelled directed acyclic graphs, and the properties of their vertices are dictated by their labels, which we define as follows.

%---------------------------------------
\begin{defi}
We call the following six diagrams (\emph{atomic-flow}) \emph{labels}:
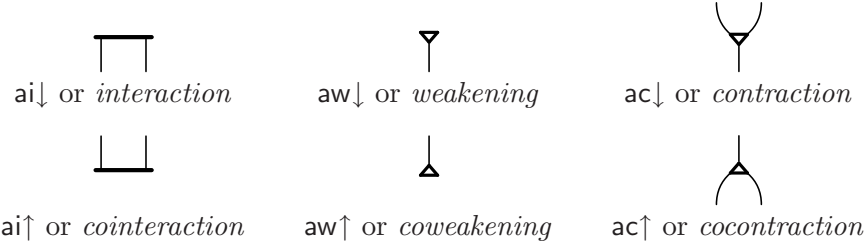
\[
\begin{array}{@{}c@{}c@{}c@{}}
      \vcenter{\afaid{}{}{}{}{}{}}&
\qquad\vcenter{\afawd{}{}{}{}}&
\qquad\vcenter{\afacd{}{}{}{}{}{}}\\
\noalign{\smallskip}
      \mbox{$\aid$ or \emph{interaction}}&
\qquad\mbox{$\awd$ or \emph{weakening}}&
\qquad\mbox{$\acd$ or \emph{contraction}}\\
\noalign{\bigskip}
      \vcenter{\afaiu{}{}{}{}{}{}}&
\qquad\vcenter{\afawu{}{}{}{}}&
\qquad\vcenter{\afacu{}{}{}{}{}{}}\\
\noalign{\smallskip}
      \mbox{$\aiu$ or \emph{cointeraction}}&
\qquad\mbox{$\awu$ or \emph{coweakening}}&
\qquad\mbox{$\acu$ or \emph{cocontraction}}\\
\end{array}\quad.
\]
Cointeraction is also called \emph{cut}.
\end{defi}

\newcommand{\ppl  }{{\mathchoice{\scriptstyle+}
                                {\scriptstyle+}
                                {\scriptstyle+}
                                {\scriptscriptstyle+}}}
\newcommand{\pmi  }{{\mathchoice{\scriptstyle-}
                                {\scriptstyle-}
                                {\scriptstyle-}
                                {\scriptscriptstyle-}}}
\newcommand{\three}{{\mathchoice{\scriptstyle\mathbf3}
                                {\scriptstyle\mathbf3}
                                {\scriptstyle\mathbf3}
                                {\scriptscriptstyle\mathbf3}}}
\newcommand{\four }{{\mathchoice{\scriptstyle\mathbf4}
                                {\scriptstyle\mathbf4}
                                {\scriptstyle\mathbf4}
                                {\scriptscriptstyle\mathbf4}}}
\newcommand{\five }{{\mathchoice{\scriptstyle\mathbf5}
                                {\scriptstyle\mathbf5}
                                {\scriptstyle\mathbf5}
                                {\scriptscriptstyle\mathbf5}}}
\newcommand{\six  }{{\mathchoice{\scriptstyle\mathbf6}
                                {\scriptstyle\mathbf6}
                                {\scriptstyle\mathbf6}
                                {\scriptscriptstyle\mathbf6}}}
\newcommand{\seven}{{\mathchoice{\scriptstyle\mathbf7}
                                {\scriptstyle\mathbf7}
                                {\scriptstyle\mathbf7}
                                {\scriptscriptstyle\mathbf7}}}
\newcommand{\eight}{{\mathchoice{\scriptstyle\mathbf8}
                                {\scriptstyle\mathbf8}
                                {\scriptstyle\mathbf8}
                                {\scriptscriptstyle\mathbf8}}}
\newcommand{\nine }{{\mathchoice{\scriptstyle\mathbf9}
                                {\scriptstyle\mathbf9}
                                {\scriptstyle\mathbf9}
                                {\scriptscriptstyle\mathbf9}}}
\newcommand{\card}[1]{\mathord\vert #1\mathord\vert}
\newcommand{\up}{{\mathit up}}
\newcommand{\lo}{{\mathit lo}}
%---------------------------------------
\begin{defi}
An (\emph{atomic}) \emph{flow} is a tuple $(V,E,\eta,\up,\lo)$ such that:
\begin{enumerate}
%-------------------
\item $V$ is a finite set of \emph{vertices}, denoted by $\nu$;
%-------------------
\item $E$ is a finite set of \emph{edges}, denoted by $\epsilon$;
%-------------------
\item $\eta\colon V\to\{\aid,\aiu,\awd,\awu,\acd,\acu\}$ maps vertices to their \emph{labels};
%-------------------
\item $\up\colon E\to V\cup\{\top\}$ and $\lo\colon E\to V\cup\{\bot\}$ are, respectively, the \emph{upper} and \emph{lower} maps, and $\top$ and $\bot$ are special vertices not belonging to $V$; we define, for every $\nu\in V\cup\{\top,\bot\}$, the set $L_\nu=\{\,\epsilon\mid\up(\epsilon)=\nu\,\}$ of \emph{lower edges of $\nu$}, the set $U_\nu=\{\,\epsilon\mid\lo(\epsilon)=\nu\,\}$ of \emph{upper edges of $\nu$}, and the set $E_\nu=L_\nu\cup U_\nu$ of \emph{edges of $\nu$};
%-------------------
\item if $\card S$ denotes the cardinality of set $S$, we have that
\begin{align*}
\mbox{if $\eta(\nu)=\aid$ then $\card{L_\nu}=2$ and $\card{U_\nu}=0$,}&\\
\mbox{if $\eta(\nu)=\aiu$ then $\card{L_\nu}=0$ and $\card{U_\nu}=2$,}&\\
\mbox{if $\eta(\nu)=\awd$ then $\card{L_\nu}=1$ and $\card{U_\nu}=0$,}&\\
\mbox{if $\eta(\nu)=\awu$ then $\card{L_\nu}=0$ and $\card{U_\nu}=1$,}&\\
\mbox{if $\eta(\nu)=\acd$ then $\card{L_\nu}=1$ and $\card{U_\nu}=2$,}&\\
\mbox{if $\eta(\nu)=\acu$ then $\card{L_\nu}=2$ and $\card{U_\nu}=1$;}&
\end{align*}
%-------------------
\item\label{ItAcycl} there is no sequence $\epsilon_1,\dots,\epsilon_h$ of edges of $V$ such that $\up(\epsilon_i)=\lo(\epsilon_{i+1\pmod h})$, for $1\le i\le h$;
%-------------------
\item\label{ItPol} there is a \emph{polarity assignment} $\pi\colon E\to\{\pmi,\ppl\}$ such that, for every $\nu\in V$,
\begin{enumerate}
%---------
\item if $\eta(\nu)\in\{\acd,\acu\}$ then $\pi(E_\nu)=\{\pmi\}$ or $\pi(E_\nu)=\{\ppl\}$;
%---------
\item if $\eta(\nu)\in\{\aid,\aiu\}$ then $\pi(E_\nu)=\{\pmi,\ppl\}$.
\end{enumerate}
\end{enumerate}
Besides $\epsilon$, we use small numerals $\one$, $\two$, \dots\ and colours to denote edges. Atomic flows are denoted with $A$, $B$, $C$ and $D$. Given an atomic flow $A$, we say that the sets $L_\top=\{\epsilon_1,\dots,\epsilon_h\}$ and $U_\bot=\{\epsilon'_1,\dots,\epsilon'_k\}$ contain, respectively, the \emph{upper} and \emph{lower edges of $A$}; in such a case, we can represent $A$ as
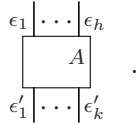
\[
\atomicflow{
( 0,10)*{\afvju4{\epsilon_1}{}};
( 2,10)*{\cdots};
( 4,10)*{\afvju4{}{\epsilon_h}};
( 5, 6)*{\aflabelleft A};
( 2, 5)*{\affr66};
( 0, 0)*{\afvjd4{\epsilon'_1}{}};
( 2, 0)*{\cdots};
( 4, 0)*{\afvjd4{}{\epsilon'_k}};
(-2, 0)*{\invisiblemark};
( 6, 0)*{\invisiblemark}}
\quad.
\]
In general, we represent atomic flows as directed-graph diagrams, except that the special vertices $\top$ and $\bot$ are not shown, and the labels of the vertices are explicitly shown as graphical elements. When we refer to the vertices of an atomic flow, we do not include $\top$ and $\bot$. Sometimes we identify vertices with their labels. 
\end{defi}

An atomic flow is a directed graph, whose edges are associated to atom occurrences in derivations, and the direction of the edges corresponds to the up-down direction in a derivation. Vertices are associated to points in the derivation where atom occurrences are created or destroyed, and the nature of each vertex is described by its label. Naturally, these graphs are acyclic (condition~\ref{ItAcycl}). The two special vertices $\top$ and $\bot$ represent the top and bottom of a derivation: we can consider $\top$ the vertex that creates all the atom occurrences in the premiss and $\bot$ the vertex that destroys all atom occurrences in the conclusion.

The polarity assignment condition (\ref{ItPol}) ensures that atoms in(co)contractions have the same polarity, and those in (co)interactions have dual polarities (as happens in derivations). Every atomic flow has $2^n$ polarity assignments, where $n$ is the number of connected components in the graph. We should not be worried about the apparent complexity of the polarity assignment condition: in fact, we could equivalently consider two sorts of (co)contraction and (co)weakening labels, the negative and the positive ones, and ask for vertices to be joined by respecting their polarities. This is clearly a locally checkable property, much simpler than, for example, some global correctness criterion for proof nets.

%---------------------------------------
\begin{exa}
Consider the atomic flow
\begin{align*}
A=(&\{\;\nu_1\;,\;\nu_2\;,\;\nu_3\;\},\\
   &\{\;\one\;,\;\two\;,\;\three\;,\;\four\;,\;\five\;\},\\
   &\{\;\nu_1\mapsto\aiu\;,\;\nu_2\mapsto\acu\;,\;\nu_3\mapsto\aiu\;\},\\
   &\{\;\one\mapsto\top\;,\;\two\mapsto\top\;,\;\three\mapsto\nu_2\;,\;
        \four\mapsto\nu_2\;,\;\five\mapsto\top\;\},\\
   &\{\;\one\mapsto\nu_1\;,\;\two\mapsto\nu_2\;,\;\three\mapsto\nu_1\;,\;
        \four\mapsto\nu_3\;,\;\five\mapsto\nu_3\;\})
\quad;
\end{align*}
the following are three of its possible representations:
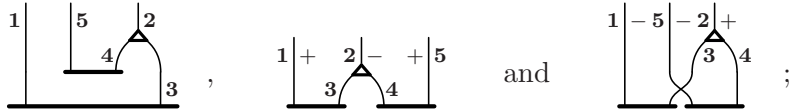
\[
\atomicflow{
(10,8)*{\afacu\four{}{}{}{}\two};
( 0,8)*{\afvjd8\one{}};
( 4,8)*{\afvjd8{}\five};
( 6,2)*{\afaiunw{}{}};
( 6,0)*{\afaiuex{}{}{}\three{}{}31}}
\quad,\qquad
\aflower{\atomicflow{
( 0  ,6)*{\afvjd{8}\one\ppl};
( 6  ,6)*{\afacu\three{}{}\four\two\pmi};
(12  ,6)*{\afvjd{8}\ppl\five};
(10  ,0)*{\afaiunw{}{}};
( 2  ,0)*{\afaiunw{}{}};
(-1.5,0)*{\invisiblemark};
(13.5,0)*{\invisiblemark}}}
\qquad\hbox{and}\qquad
\atomicflow{
( 8  ,10)*{\afacu{}\three{}\four\two\ppl};
( 0  , 8)*{\afvjd{12}\one\pmi};
( 4  ,10)*{\afvjd{8}\five\pmi};
( 5  , 4)*{\afex24};
(10  , 4)*{\afvj4};
( 2  , 0)*{\afaiunw{}{}};
( 8  , 0)*{\afaiunw{}{}};
(-1.5, 0)*{\invisiblemark};
(11.5, 0)*{\invisiblemark}}
\quad;
\]
\afnegspace
in the last two diagrams, we also indicated each of the two possible polarity assignments. This flow has one cocontraction and two cointeraction vertices; it has three upper edges, $\one$, $\two$ and $\five$, and no lower edges.
\end{exa}

%---------------------------------------
\afnegspace
\begin{exa}
The graph
$\atomicflow{
(0,4)*{\afaidnw{}{}};
(0,0)*{\afacd{}{}{}{}{}{}}}$
is not an atomic flow, for lack of a polarity assignment.
\end{exa}

We now define the mapping from derivations to atomic flows. As we said, the idea is that structural rules map to the respective atomic-flow vertices, and the edges trace the atoms between inference steps. We first state a fact, whose proof is immediate.

%---------------------------------------
\begin{prop}\label{PropUnFl}
Given an\/ $\SKS$ derivation\/ $\Phi$, there is a unique atomic flow $A$ (modulo isomorphisms) such that:
\begin{enumerate}
%-------------------
\item there is a surjective map between the set of atom occurrences of\/ $\Phi$ and the set of edges of $A$;
%-------------------
\item for each inference step $\vlsmash{\vlinf{\rho}{}{\xi\{\beta\}}{\xi\{\alpha\}}}$ of\/ $\Phi$, where $\rho\in\{\aid,\aiu,\awd,\awu,\acd,\acu\}$ and $\vlinf{\rho}{}{\beta}{\alpha}$ is a rule instance, all atom occurrences in $\xi\vlhole$ in the premiss are respectively mapped to the same edges of $A$ as the atom occurrences in $\xi\vlhole$ in the conclusion; the atom occurrences in $\vlinf{\rho}{}{\beta}{\alpha}$ are mapped to edges of $A$ such that the edges are related with vertices as indicated below, for each possible case of the inference step:
\[
\begin{array}{@{}ccc@{}ccc@{}}
\vlinf{\aid}{}{\vls[a^\one.{\bar a^\two}]}{\ttt}&\mbox{to\/}&
\vcenter{\afaid\one{}{}\two{}{}}
\quad,&\qquad
\vlinf{\aiu}{}{\fff}{\vls(a^\one.{\bar a^\two})}&\mbox{to\/}&
\vcenter{\afaiu\one{}{}\two{}{}}
\quad,\\
\noalign{\medskip}
\vlinf{\awd}{}{a^\one}{\fff}                    &\mbox{to\/}&
\vcenter{\afawd{}{}{}\one{}} 
\quad,&\qquad
\vlinf{\awu}{}{\ttt}{a^\one}                    &\mbox{to\/}&
\vcenter{\afawu{}{}{}\one{}}
\quad,\\
\noalign{\medskip}
\vlinf{\acd}{}{a^\three}{\vls[a^\one.a^\two]}   &\mbox{to\/}&
\vcenter{\afacd\one{}{}\two{}\three}
\quad,&\qquad
\vlinf{\acu}{}{\vls(a^\two.a^\three)}{a^\one}   &\mbox{to\/}&
\vcenter{\afacu\two{}{}\three{}\one}
\quad,\\
\end{array}
\]
where the mapping is indicated by small numerals.
%-------------------
\item for each inference step of\/ $\Phi$ of kind
\[\hss
\begin{array}{@{}r@{}l@{}}
\vlinf{\swi}{}{\xi\vlscn[(\alpha.\beta).\gamma]}
              {\xi\vlscn(\alpha.[\beta.\gamma])}           \quad,&\qquad
\vlinf{\med}{}{\xi\vlscn([\alpha.\gamma].[\beta.\delta])}
              {\xi\vlscn[(\alpha.\beta).(\gamma.\delta)]}  \quad,      \\
\noalign{\smallskip}
\vlinf={}{\xi\vlscn[\beta.\alpha]}{\xi\vlscn[\alpha.\beta]}\quad,&\qquad
\vlinf={}{\xi\vlscn(\beta.\alpha)}{\xi\vlscn(\alpha.\beta)}\quad,      \\
\noalign{\smallskip}
\vlinf={}{\xi\vlscn[\alpha.[\beta.\gamma]]}
         {\xi\vlscn[[\alpha.\beta].\gamma]}                \quad,&\qquad
\vlinf={}{\xi\vlscn[[\alpha.\beta].\gamma]}
         {\xi\vlscn[\alpha.[\beta.\gamma]]}                \quad,      \\
\noalign{\smallskip}
\vlinf={}{\xi\vlscn(\alpha.(\beta.\gamma))}
         {\xi\vlscn((\alpha.\beta).\gamma)}                \quad,&\qquad
\vlinf={}{\xi\vlscn((\alpha.\beta).\gamma)}
         {\xi\vlscn(\alpha.(\beta.\gamma))}                \quad,      \\
\noalign{\smallskip}
\vlinf={}{\xi\{\alpha\}}{\xi\vlscn[\alpha.\fff]}           \quad,\qquad
\vlinf={}{\xi\vlscn[\alpha.\fff]}{\xi\{\alpha\}}           \quad,&\qquad
\vlinf={}{\xi\{\alpha\}}{\xi\vlscn(\alpha.\ttt)}        \qquad\hbox{and\/}\qquad
\vlinf={}{\xi\vlscn(\alpha.\ttt)}{\xi\{\alpha\}}
\end{array}
\]
all the atom occurrences in $\xi\vlhole$, $\alpha$, $\beta$, $\gamma$ and $\delta$ in the premiss are respectively mapped to the same edges of $A$ as the atom occurrences in $\xi\vlhole$, $\alpha$, $\beta$, $\gamma$ and $\delta$ in the conclusion.
\end{enumerate}
\end{prop}

%---------------------------------------
\begin{defi}
Given a derivation $\Phi$, we say that the unique atomic flow $A$ defined in Proposition~\ref{PropUnFl} is the atomic flow \emph{associated with} the derivation $\Phi$. Sometimes, when an atom occurrence $a$ in $\Phi$ maps to an edge $\epsilon$ in $A$, we decorate $\epsilon$ with the label $a$.
\end{defi}

%---------------------------------------
\begin{exa}
Figure~\ref{FigExAF} has some examples of atomic flows associated with derivations.
\end{exa}

\newcommand{\RD}[1]{{\color{Red}#1}}
\newcommand{\GR}[1]{{\color{Green}#1}}
\newcommand{\DO}[1]{{\color{DarkOrchid}#1}}
\newcommand{\PB}[1]{{\color{ProcessBlue}#1}}
\newcommand{\MG}[1]{{\color{Magenta}#1}}
\newcommand{\SG}[1]{{\color{SpringGreen}#1}}
\newcommand{\RS}[1]{{\color{RawSienna}#1}}
\newcommand{\YO}[1]{{\color{YellowOrange}#1}}
\newcommand{\PW}[1]{{\color{Periwinkle}#1}}
%---------------------------------------
\begin{figure}[tbp]
\[
\begin{array}{@{}c@{}c@{}c@{}}
\vlderivation                                                  {
\vlin{=   }{}{\ttt                                  }         {
\vlin{\aiu}{}{\vls[\fff.\ttt]                       }        {
\vlin{=   }{}{\vls[(\GR{a}.\RD{\bar a}).\ttt]       }       {
\vlin{\swi}{}{\vls[[(\RD{\bar a}.\GR{a}).\ttt].\ttt]}      {
\vlin{=   }{}{\vls[(\RD{\bar a}.[\GR{a}.\ttt]).\ttt]}     {
\vlin{\swi}{}{\vls[([\GR{a}.\ttt].\RD{\bar a}).\ttt]}    {
\vlin{=   }{}{\vls([\GR{a}.\ttt].[\RD{\bar a}.\ttt])}   {
\vlin{\med}{}{\vls([\GR{a}.\ttt].[\ttt.\RD{\bar a}])}  {
\vlin{=   }{}{\vls[(\GR{a}.\ttt).(\ttt.\RD{\bar a})]} {
\vlin{\aid}{}{\vls[\GR{a}.\RD{\bar a}]              }{
\vlhy        {\ttt                                  }}}}}}}}}}}}
\qquad&
\vlderivation                                                              {
\vlin{\aiu}{}
   {\vls(\DO{a}.\fff)                                            }        {
\vlin{=   }{}
   {\vls(\DO{a}.(\PB{a}.\MG{\bar a}))                            }       {
\vlin{\acu}{}
   {\vls((\DO{a}.\PB{a}).\MG{\bar a})                            }      {
\vlin{=   }{}
   {\vls(\SG{a}.\MG{\bar a})                                     }     {
\vlin{\aiu}{}
   {\vls([\fff.\SG{a}].\MG{\bar a})                              }    {
\vlin{\acd}{}
   {\vls([(\RD{a}.\RS{\bar a}).\SG{a}].\MG{\bar a})              }   {
\vlin{\swi}{}
   {\vls([(\RD{a}.[\GR{\bar a}.\YO{\bar a}]).\SG{a}].\MG{\bar a})}  {
\vlin{=   }{}
   {\vls((\RD{a}.[[\GR{\bar a}.\YO{\bar a}].\SG{a}]).\MG{\bar a})} {
\vlin{\aid}{}
   {\vls(\RD{a}.[\GR{\bar a}.[\YO{\bar a}.\SG{a}]].\MG{\bar a})  }{
\vlhy        
   {\vls(\RD{a}.[\GR{\bar a}.\ttt].\MG{\bar a})                  }}}}}}}}}}}
\qquad&
\vlderivation                                                              {
\vlin{=   }{}{\vls(([\RS{a}.\YO{b}].\PW{c}).([\GR{a}.\DO{b}].\SG{c}))}    {
\vlin{\med}{}{\vls(([\RS{a}.\YO{b}].[\GR{a}.\DO{b}]).(\PW{c}.\SG{c}))}   {
\vlin{\acu}{}{\vls([(\RS{a}.\GR{a}).(\YO{b}.\DO{b})].(\PW{c}.\SG{c}))}  {
\vlin{\acu}{}{\vls([(\RS{a}.\GR{a}).(\YO{b}.\DO{b})].\MG{c})         } {
\vlin{\acu}{}{\vls([(\RS{a}.\GR{a}).\PB{b}].\MG{c})                  }{
\vlhy        {\vls([\RD{a}.\PB{b}].\MG{c})                           }}}}}}}\\
\atomicflow{
(0,0)*{\afaiucol{}{}{}{}{}{}{Green}{Red}{}};
(0,4)*{\afaidnw{}{}}}
\qquad&
\atomicflow{
( 2,14)*{\afvjcol4{Green}};
( 0,10)*{\afvjcol{12}{Red}};
(16,10)*{\afvjcol{12}{Magenta}};
( 4, 8)*{\afacdcol{}{}{}{}{}{}{Green}{YellowOrange}{RawSienna}};
(10, 8)*{\afacucol{}{}{}{}{}{}{DarkOrchid}{ProcessBlue}{SpringGreen}};
( 2, 2)*{\afaiunw{}{}};
( 8, 2)*{\afvjcol4{DarkOrchid}};
(14, 2)*{\afaiunw{}{}};
( 8,12)*{\afaidnw{}{}}}
\qquad&
\atomicflow{
( 0,0)*{\afacucol{}{}{}{}{}{}{RawSienna}{Green}{Red}};
(10,0)*{\afacucol{}{}{}{}{}{}{YellowOrange}{DarkOrchid}{ProcessBlue}};
(20,0)*{\afacucol{}{}{}{}{}{}{Periwinkle}{SpringGreen}{Magenta}}}
\end{array}
\]
\caption{Examples of atomic flows associated with derivations.}
\label{FigExAF}
\end{figure}

Inference rules are usually called linear when they do not `create' nor `destroy' atoms. Linear rules of $\SKS$ are switch, medial and (every equation defining) rule $=$. Note that linear inference rules do not introduce any vertices in atomic flows.

We now show that there is no such thing as an `invalid atomic flow', or, in other words, the mapping from derivations to atomic flows is surjective.

%---------------------------------------
\begin{thm}\label{TheoFlowDer}
Every atomic flow is associated with some derivation.
\end{thm}

%---------------------------------------
\proof
First, we construct a derivation scheme $\Phi$ that `glues together' any two atomic flows without introducing any vertices. For every $\alpha$, $\beta$, the atomic flow of the following derivation consists only of edges (there are no vertices):
\[\Psi=
\vlderivation                                                {
\vlin{=   }{}{\vls[(\alpha.\beta).\ttt]              }      {
\vlin{\swi}{}{\vls[[(\beta.\alpha).\ttt].\ttt]       }     {
\vlin{=   }{}{\vls[(\beta.[\alpha.\ttt]).\ttt]       }    {
\vlin{\swi}{}{\vls[[([\alpha.\ttt].\beta).\ttt].\ttt]}   {
\vlin{=   }{}{\vls[([\alpha.\ttt].[\beta.\ttt]).\ttt]}  {
\vlin{\med}{}{\vls[([\alpha.\ttt].[\ttt.\beta]).\ttt]} {
\vlin{=   }{}{\vls[[(\alpha.\ttt).(\ttt.\beta)].\ttt]}{
\vlhy        {\vls[[\alpha.\beta].\ttt]              }}}}}}}}}
\quad.
\]
We use $\Psi$ and $\sus$ (see Remark~\ref{RemSupSwitch}) to `move' an atom $a$ from one context $\zeta\vlhole$ to another context $\xi\vlhole$, again with an associated atomic flow that is free of vertices:
\[
\vlderivation                                                    {
\vlde{\Psi}{\{\swi,\med\}}{\vls[(\xi\{a\}.\zeta\{\fff\}).\ttt]} {
\vlin{\sus}{             }{\vls[[\xi\{a\}.\zeta\{\fff\}].\ttt]}{
\vlhy                     {\vls[(\xi\{\ttt\}.\zeta\{a\}).\ttt]}}}}
\quad.
\]
This construction can be used zero or more times to get the desired $\Phi$, for $h\ge0$:
\[
\vlderivation                                                             {
\vlde{\Phi}{\{\swi,\med\}}
     {\vls[(\xi\{a_1 \}\cdots\{a_h \}.\zeta\{\fff\}\cdots\{\fff\}).\ttt]}{
\vlhy{\vls[(\xi\{\ttt\}\cdots\{\ttt\}.\zeta\{a_1 \}\cdots\{a_h \}).\ttt]}}}
\quad.
\]
We can now prove the theorem by induction on the number of vertices of a given atomic flow $A$. The cases where $A$ only has zero or one vertex are trivial. Let us then suppose that $A$ has more than one vertex; then $A$ can be considered as composed of two flows $B$ and $C$, each with less vertices than $A$, as follows:
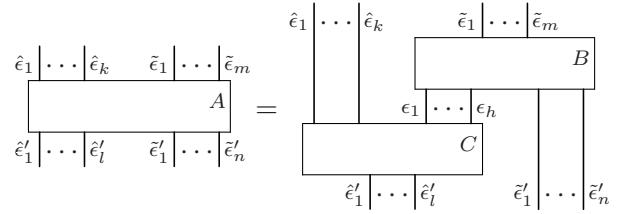
\[
\atomicflow{
( 0  ,10)*{\afvjd4{\hat\epsilon_1}{}};
( 2  ,10)*{\cdots};
( 4  ,10)*{\afvjd4{}{\hat\epsilon_k}};
(12  ,10)*{\afvjd4{\tilde\epsilon_1}{}};
(14  ,10)*{\cdots};
(16  ,10)*{\afvjd4{}{\tilde\epsilon_m}};
(17  , 6)*{\aflabelleft A};
( 8  , 5)*{\affr{18}6};
( 0  , 0)*{\afvju4{\hat\epsilon'_1}{}};
( 2  , 0)*{\cdots};
( 4  , 0)*{\afvju4{}{\hat\epsilon'_l}};
(12  , 0)*{\afvju4{\tilde\epsilon'_1}{}};
(14  , 0)*{\cdots};
(16  , 0)*{\afvju4{}{\tilde\epsilon'_n}};
(-2  , 0)*{\invisiblemark};
(18.5, 0)*{\invisiblemark}}
=
\atomicflow{
(  6  ,20)*{\afvjd4{}{\tilde\epsilon_m}};
(  4  ,20)*{\cdots};
(  2  ,20)*{\afvjd4{\tilde\epsilon_1}{}};
(-11  ,20)*{\cdots};
( 12  ,16)*{\aflabelleft B};
(  4  ,15)*{\affr{16}6};
( -9  ,15)*{\afvjd{14}{}{\hat\epsilon_k}};
(-13  ,15)*{\afvjd{14}{\hat\epsilon_1}{}};
(  1  ,10)*{\afvju4{}{\epsilon_h}};
( -1  ,10)*{\cdots};
( -3  ,10)*{\afvju4{\epsilon_1}{}};
(  2  , 6)*{\aflabelleft C};
( 11  , 5)*{\afvju{14}{}{\tilde\epsilon'_n}};
(  7  , 5)*{\afvju{14}{\tilde\epsilon'_1}{}};
( -6  , 5)*{\affr{16}6};
(  9  , 0)*{\cdots};
( -4  , 0)*{\afvju4{}{\hat\epsilon'_l}};
( -6  , 0)*{\cdots};
( -8  , 0)*{\afvju4{\hat\epsilon'_1}{}};
( 13  , 0)*{\invisiblemark};
(-15.5, 0)*{\invisiblemark}}
\quad,
\]
where $h,k,l,m,n\ge0$ (this can possibly be done in many different ways). By the inductive hypothesis, there exist derivations $\vlder{\Phi_B}{}{\zeta\{a_1^{\epsilon_1}\}\cdots\{a_h^{\epsilon_h}\}}{\gamma}$ and $\vlder{\Phi_C}{}{\delta}{\xi\{a_1^{\epsilon_1}\}\cdots\{a_h^{\epsilon_h}\}}$ whose flows are, respectively, $B$ and $C$. Using these, we can build
\[
\vlderivation                                                             {
\vlde{\vls[(\Phi_C.\zeta\{\fff\}\cdots\{\fff\}).\ttt]}{}
     {\vls[(\delta.\zeta\{\fff\}\cdots\{\fff\}).\ttt]                 }  {
\vlde{\Phi}{}
     {\vls[(\xi\{a_1^{\epsilon_1}\}\cdots\{a_h^{\epsilon_h}\}.
            \zeta\{\fff\}\cdots\{\fff\}).\ttt]                        } {
\vlde{\vls[(\xi\{\ttt\}\cdots\{\ttt\}.\Phi_B).\ttt]}{}
     {\vls[(\xi\{\ttt\}\cdots\{\ttt\}.
            \zeta\{a_1^{\epsilon_1}\}\cdots\{a_h^{\epsilon_h}\}).\ttt]}{
\vlhy{\vls[(\xi\{\ttt\}\cdots\{\ttt\}.\gamma).\ttt]                   }}}}}
\quad,
\]
whose flow is $A$.
\qed

The derivations in the constructions for the theorem above involve tautologies of the kind $\vls[\alpha.\ttt]$. There, the `logical content' of $\alpha$ does not matter because it is trivialized by $\ttt$, and we are free to build derivations with any premiss and conclusion. A challenging question is whether a derivation exists with a given associated flow and with given premiss and conclusion; in this case, we cannot resort to logical units to trivialize derivations. In the following, in particular in Sections~\ref{SectRed} and \ref{SectStreamAlg}, we reason about derivations whose premiss and conclusion are arbitrary and have to be preserved through transformations of derivations and their associated flows.

%-------------------------------------------------------------------------------
\subsection{Paths and Cycles}

We now define the notions of `$\ai$-path' and `$\ai$-cycle' in atomic flows. Paths are sequences of adjacent edges that only `go down' or only `go up'; $\ai$-paths are formed by joining paths at interaction or cointeraction vertices; $\ai$-cycles are circular $\ai$-paths. We also define the notion of `simple edge', \emph{i.e.}, an edge connecting an interaction and a cointeraction, as in the first diagram in this paper.

%---------------------------------------
\begin{defi}
Given an atomic flow $(V,E,\eta,\up,\lo)$ and $\epsilon_1,\dots,\epsilon_h\in E$ such that, for $1\le i<h$, we have $\lo(\epsilon_i)=\up(\epsilon_{i+1})$, $\up(\epsilon_1)=\nu$ and $\lo(\epsilon_h)=\nu'$, we say that $\epsilon_1,\dots,\epsilon_h$ is a \emph{path from $\nu$ to $\nu'$} and that $\epsilon_h,\dots,\epsilon_1$ is a \emph{path from $\nu'$ to $\nu$}; both paths have \emph{length} $h$. An \emph{$\ai$-path from $\nu$ to $\nu'$} of \emph{length} $h$ is either a path from $\nu$ to $\nu'$ of length $h$ or a sequence of edges $\epsilon_1,\dots,\epsilon_k,\epsilon_{k+1},\dots,\epsilon_h$ such that $\epsilon_k \ne \epsilon_{k+1}$ and, for some $\nu''\in V$ with $\eta(\nu'')\in\{\aid,\aiu\}$, we have that $\epsilon_1,\dots,\epsilon_k$ is an $\ai$-path from $\nu$ to $\nu''$ and $\epsilon_{k+1},\dots,\epsilon_h$ is an $\ai$-path from $\nu''$ to $\nu'$. An $\ai$-path of length $h$ is \emph{maximal} if no $\ai$-path containing its edges has length greater than $h$. An $\ai$-path from (resp., to) $\nu$ of length $h$ is a \emph{maximal\/ $\ai$-path from} (resp., \emph{to}) $\nu$ if no $\ai$-path from (resp., to) $\nu$ containing its edges has length greater than $h$. A path from an interaction to a cointeraction vertex or vice versa is called an \emph{$\ai$-connection}.
\end{defi}

%---------------------------------------
\begin{exa}
The atomic flow on the left has the $\ai$-paths on the right, and the paths are marked with an asterisk:
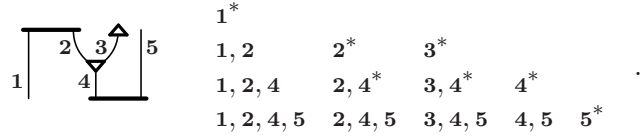
\[
\begin{array}{@{}c@{}c@{}}
\atomicflow{
( 2  ,10)*{\afaidnw{}{}};
( 0  , 6)*{\afvju8\one{}};
( 6  , 6)*{\afacd\two{}\three{}\four{}};
( 8  , 8)*{\afawunw{}{}};
(10  , 6)*{\afvjd8{}\five};
( 8  , 0)*{\afaiunw{}{}};
(-1.5, 0)*{\invisiblemark};
(11.5, 0)*{\invisiblemark}}
\qquad&
\begin{array}{@{}lllll@{}}
\one^*               &                &                  &           &       \\
\one,\two            &\two^*          &\three^*          &           &       \\
\one,\two,\four      &\two,\four^*    &\three,\four^*    &\four^*    &       \\
\one,\two,\four,\five&\two,\four,\five&\three,\four,\five&\four,\five&\five^*\\
\end{array}\quad.
\\
\end{array}
\]
In addition, the flow has the paths and $\ai$-paths obtained from the shown ones by inverting the order of edges, for example $\five,\four,\two,\one$ is an $\ai$-path. The $\ai$-paths from the interaction vertex are $\one$ and $\two$ and $\two,\four$ and $\two,\four,\five$; the $\ai$-paths to the contraction vertex are $\one,\two$ and $\two$ and $\three$ and $\four$ and $\five,\four$; of all the $\ai$-paths to the cointeraction vertex, $\two,\four$ is the only $\ai$-connection; the only other $\ai$-connection in the flow is $\four,\two$. The maximal $\ai$-paths are $\one,\two,\four,\five$ and $\three,\four,\five$ and their inverses. The maximal $\ai$-paths from the cointeraction vertex are $\four,\two,\one$ and $\four,\three$ and $\five$; the maximal $\ai$-paths to the contraction vertex are $\one,\two$ and $\three$ and $\five,\four$.
\end{exa}

Simple edges represent immediate causality relations between axioms and cuts. They play a crucial role in the following, in particular when they belong to $\ai$-cycles.

%---------------------------------------
\begin{defi}\label{DefSimple}
An $\ai$-connection consisting of a single edge is called a \emph{simple edge}. A \emph{clean path} is an $\ai$-path where every $\ai$-connection is a simple edge. An \emph{$\ai$-cycle} is an $\ai$-path from a vertex to itself, where no edge appears twice; we do not distinguish $\ai$-cycles that only differ for cyclic permutations of their edges or for inversion, so, if $\epsilon_1,\dots,\epsilon_h$ is an $\ai$-cycle, then $\epsilon_2,\dots,\epsilon_h,\epsilon_1$ and $\epsilon_h,\dots,\epsilon_1$ are the same $\ai$-cycle. A \emph{fragile cycle} is an $\ai$-cycle containing a simple edge. Atomic flows and $\ai$-paths are both called \emph{cycle-free} if they do not contain $\ai$-cycles.
\end{defi}

%---------------------------------------
\begin{exa}\label{ExExtrem}
In the following cycle-free flow all the $\ai$-paths are clean paths; $\one$, $\two$ and $\three$ are $\ai$-connections and simple edges:
\[
\atomicflow{
( 2,0)*{\afaiu{}{}\one{}{}{}};
( 0,4)*{\afacdnw{}{}{}{}};
( 6,4)*{\afaidnw{}{}};
(10,0)*{\afaiu\two{}\three{}{}{}};
(14,4)*{\afaidnw{}{}{}{}};
(18,0)*{\afacd{}{}{}{}{}{}}}
\quad.
\]
\end{exa}

%---------------------------------------
\begin{exa}
Consider the following atomic flow:
\nopagebreak[4]\medskip\afnegspace
\[
\atomicflow{
( 8,12)*{\afaidex\one{}{}\two{}{}41};
( 0, 6)*{\afacunw\three{}{}\four};
(16, 6)*{\afacunw\seven{}{}\eight};
( 8, 8)*{\afaid{}\five\six{}{}{}};
( 4, 2)*{\afaiunw{}{}};
(12, 2)*{\afaiunw{}{}};
( 8, 0)*{\afaiuex{}{}{}{}{}{}51}}
\quad.
\]
\afnegspace
The flow contains two $\ai$-cycles: $\one,\four,\five,\six,\seven,\two$ and $\one,\three,\eight,\two$. The first $\ai$-cycle contains the two simple edges $\five$ and $\six$, so it is a fragile cycle; the second $\ai$-cycle does not contain any simple edge. Note that the two $\ai$-cycles are `overlapping', in the sense that edges $\one$ and $\two$ belong to both.
\end{exa}

%---------------------------------------
\begin{rem}
If an $\ai$-path is maximal, then it is cycle-free.
\end{rem}

%-------------------------------------------------------------------------------
\subsection{Streamlined Derivations}

Intuitively, we can consider interaction and weakening as creators of paths, and we can consider cointeraction and coweakening as their destroyers. We call a derivation `streamlined' if no path is both created and destroyed.

%---------------------------------------
\begin{defi}
An $\SKS$ derivation is \emph{streamlined} if, in its associated atomic flow, there are no paths from interaction or weakening vertices to cointeraction or coweakening vertices.
\end{defi}

%---------------------------------------
\begin{rem}\label{RemStr}
It immediately follows from the definition that the diagram below describes the shape of a streamlined derivation; the boxes stand for flows obtained by freely composing edges and vertices whose labels are only those indicated on the boxes:
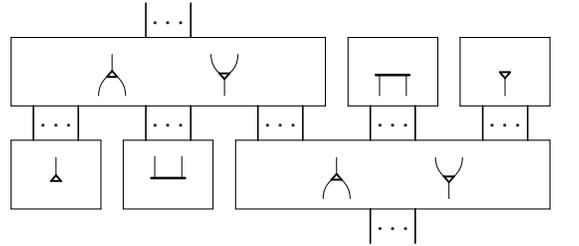
\[
\newbox\boxone\setbox\boxone=\hbox{$
   \divide\atflowunit by5\multiply\atflowunit by3\afsetunits
   \atomicflow{(0,0)*{\afacd{}{}{}{}{}{}}}$}
\newbox\boxtwo\setbox\boxtwo=\hbox{$
   \divide\atflowunit by5\multiply\atflowunit by3\afsetunits
   \atomicflow{(0,0)*{\afacu{}{}{}{}{}{}}}$}
\newbox\boxthree\setbox\boxthree=\hbox{$
   \divide\atflowunit by5\multiply\atflowunit by3\afsetunits
   \atomicflow{(0,0)*{\afaid{}{}{}{}{}{}}}$}
\newbox\boxfour\setbox\boxfour=\hbox{$
   \divide\atflowunit by5\multiply\atflowunit by3\afsetunits
   \atomicflow{(0,0)*{\afaiu{}{}{}{}{}{}}}$}
\newbox\boxfive\setbox\boxfive=\hbox{$
   \divide\atflowunit by5\multiply\atflowunit by3\afsetunits
   \atomicflow{(0,0)*{\afawd{}{}{}{}{}{}}}$}
\newbox\boxsix\setbox\boxsix=\hbox{$
   \divide\atflowunit by5\multiply\atflowunit by3\afsetunits
   \atomicflow{(0,0)*{\afawu{}{}{}{}{}{}}}$}
\atomicflow{
(-12,12)*{\afvj4};
(-10,12)*{\cdots};
( -8,12)*{\afvj4};
%---
(-15, 6)*{\copy\boxtwo};
(-10, 6)*{\affr{28}8};
( -5, 6)*{\copy\boxone};
( 10, 6)*{\copy\boxthree};
( 10, 6)*{\affr88};
( 20, 6)*{\copy\boxfive};
( 20, 6)*{\affr88};
%---
(-22, 0)*{\afvj4};
(-20, 0)*{\cdots};
(-18, 0)*{\afvj4};
(-12, 0)*{\afvj4};
(-10, 0)*{\cdots};
( -8, 0)*{\afvj4};
( -2, 0)*{\afvj4};
(  0, 0)*{\cdots};
(  2, 0)*{\afvj4};
(  8, 0)*{\afvj4};
( 10, 0)*{\cdots};
( 12, 0)*{\afvj4};
( 18, 0)*{\afvj4};
( 20, 0)*{\cdots};
( 22, 0)*{\afvj4};
%---
(-20,-6)*{\copy\boxsix};
(-20,-6)*{\affr88};
(-10,-6)*{\copy\boxfour};
(-10,-6)*{\affr88};
(  5,-6)*{\copy\boxtwo};
( 10,-6)*{\affr{28}8};
( 15,-6)*{\copy\boxone};
%---
(  8,-12)*{\afvj4};
( 10,-12)*{\cdots};
( 12,-12)*{\afvj4}}\quad.
\]
\end{rem}

%---------------------------------------
\begin{exa}\label{ExStream}
The first flow is not streamlined, the other two are streamlined:
\[
\aflower{\atomicflow{
(0,6)*{\afacd{}{}{}{}{}{}};
(4,6)*{\afawd{}{}{}{}};
(2,0)*{\afaiunw{}{}}}}
\quad,\qquad
\atomicflow{
( 2,0)*{\afacd{}{}{}{}{}{}};
( 0,4)*{\afawdnw{}{}};
( 6,8)*{\afacu{}{}{}{}{}{}};
(10,2)*{\afaiunw{}{}{}{}};
(12,8)*{\afvj8}}
\qquad\hbox{and}\qquad
\atomicflow{
(10,8)*{\afacu{}{}{}{}{}{}};
(16,8)*{\afvj8};
( 2,4)*{\afaidnw{}{}};
( 0,0)*{\afvj8};
(14,2)*{\afaiunw{}{}};
( 6,0)*{\afacd{}{}{}{}{}{}}}
\quad.
\]
\end{exa}

%---------------------------------------
\begin{rem}\label{RemStrSKS}
A streamlined $\SKS$ proof is cut-free. In fact, consider the diagram in Remark~\ref{RemStr}: there can be no $\ai$-cycles, and any maximal $\ai$-path from a cut goes all the way to the top. So, if there are cuts, there must be atoms in the premiss of the derivation, and it cannot be a proof. Note also that a cut-free proof is not necessarily streamlined, as it might have paths from interactions to coweakenings.
\end{rem}

%---------------------------------------
\begin{defi}
We say that an algorithm $P$ \emph{streamlines} a derivation $\Phi$ if the output of $P$ on $\Phi$ is a streamlined derivation that has the same premiss and conclusion as $\Phi$, in the same deductive system.
\end{defi}

We know that a Craig interpolant of formulae $\alpha$ and $\beta$ is a formula $\gamma$ such that $\alpha\vlim\gamma$ and $\gamma\vlim\beta$, and all the atoms in $\gamma$ appear in $\alpha$ and $\beta$. Consider a derivation whose premiss is $\alpha$ and conclusion is $\beta$ and consider all the $\ai$-paths from atoms of $\alpha$ that are not in $\beta$ and from atoms of $\beta$ that are not in $\alpha$. While, in general, each of these $\ai$-paths can be composed of an arbitrary number of paths, in a streamlined derivation they consist of at most two paths. We would then be tempted to exploit this simplification to try and read interpolants as some intermediate formulae in streamlined derivations. Unfortunately, this does not work so simply, as the following example shows.

%---------------------------------------
\begin{exa}
The two following streamlined derivations have the same premiss $\alpha=\vls(a.\bar a)$ and  conclusion $\beta=\vls[b.\bar b]$ and the same atomic flow (at the centre):
\[
\vlderivation                             {
\vlin{\aid}{}{\vls        [b.\bar b]}    {
\vlin{=   }{}{                 \ttt }   {
\vlin{\swi}{}{\vls[(\fff.\fff).\ttt]}  {
\vlin{=   }{}{\vls(\fff.[\fff.\ttt])} {
\vlin{\aiu}{}{     \fff             }{
\vlhy        {\vls(a.\bar a)        }}}}}}}
\qquad\qquad
\atomicflow{
(0,4)*{\afaiu{}{}{}{}{}{}};
(0,0)*{\afaid{}{}{}{}{}{}}}
\qquad\qquad
\vlderivation                                         {
\vlin{=   }{}{\vls                   [b.\bar b] }    {
\vlin{\aiu}{}{\vls[(\fff      .\fff).[b.\bar b]]}   {
\vlin{\swi}{}{\vls[((a.\bar a).\fff).[b.\bar b]]}  {
\vlin{\aid}{}{\vls((a.\bar a).[\fff.[b.\bar b]])} {
\vlin{=   }{}{\vls((a.\bar a).[\fff.\ttt      ])}{
\vlhy        {\vls (a.\bar a)                   }}}}}}}
\quad.
\]
The formulae that appear in the derivation on the left, apart from the premiss and conclusion, are interpolants of $\alpha$ and $\beta$. In the derivation on the right, the same is not true.
\end{exa}

The results we get in this paper for derivations share some characteristics with Craig's interpolation for formulae: 1) we focus on a normal form that essentially depends on the atoms in common between premiss and conclusion, 2) the cost of getting it is exponential (to the best of our knowledge), and 3) it is intimately related to cut elimination.

We are able to use the methods presented here in order to read interpolants from derivations, so overcoming the problem shown in the previous example. We start from streamlined derivations and then we perform some further constructions; this way, we obtain a normal form such that each inference step of the kind $\aiu$, $\awu$ or $\acu$ is above all inference steps of the kind $\aid$, $\awd$ or $\acd$. This normal form is obtained by Kai Br\"unnler in \cite{Brun:06:Deep-Inf:qy} by resorting to the sequent calculus, while we can obtain it directly in the calculus of structures. These results are presented in \cite{GuglGund:08:Normalis:yu}.

%===============================================================================
\section{Reductions of Atomic Flows}\label{SectRed}

We control normalisation of derivations by manipulating atomic flows, in the sense of graph rewriting. There are two kinds of flow reductions: local and global ones. In local reductions, a bounded subflow in a flow is substituted by another subflow that fits in the context. In global reductions, the entire flow is rewritten: normally, two slightly altered copies of a flow are connected together. In this section, we see local transformations, which are based on reduction rules; in Section~\ref{SectStreamAlg}, we deal with global reductions. It is convenient to classify reduction rules into those for weakening and those for contraction. After seeing flow reductions and tying them with derivations, in Subsection~\ref{SubsRed}, we explore some of their basic properties, in the two short Subsections~\ref{SubsWeak} and \ref{SubsContr}.

%-------------------------------------------------------------------------------
\subsection{Reductions}\label{SubsRed}

We introduce reductions for atomic flows, for which we define a concept of soundness. We will soon see that, corresponding to every sound reduction, there is a transformation on derivations that preserves premiss and conclusion. We will then be able to control complex proof transformations by simple atomic flow transformations.

We start by defining some flow reductions: they are relations $A\to B$ between flows, which we interpret as the possibility of replacing flow $A$ with flow $B$.

%---------------------------------------
\begin{defi}
In Figure~\ref{FigRed}, we define graphical expressions of the kind $r\colon A\to B$, where $r$ is a name and $A$ and $B$ are flows.
\end{defi}

\newcommand{\rwdcd}{{{\mathsf w}{\downarrow}{\hbox{-}}{\mathsf c}{\downarrow}}}
\newcommand{\rwdiu}{{{\mathsf w}{\downarrow}{\hbox{-}}{\mathsf i}{\uparrow  }}}
\newcommand{\rwdwu}{{{\mathsf w}{\downarrow}{\hbox{-}}{\mathsf w}{\uparrow  }}}
\newcommand{\rwdcu}{{{\mathsf w}{\downarrow}{\hbox{-}}{\mathsf c}{\uparrow  }}}
\newcommand{\rcuwu}{{{\mathsf c}{\uparrow  }{\hbox{-}}{\mathsf w}{\uparrow  }}}
\newcommand{\rcdwu}{{{\mathsf c}{\downarrow}{\hbox{-}}{\mathsf w}{\uparrow  }}}
\newcommand{\rcdiu}{{{\mathsf c}{\downarrow}{\hbox{-}}{\mathsf i}{\uparrow  }}}
\newcommand{\rcdcu}{{{\mathsf c}{\downarrow}{\hbox{-}}{\mathsf c}{\uparrow  }}}
\newcommand{\ridwu}{{{\mathsf i}{\downarrow}{\hbox{-}}{\mathsf w}{\uparrow  }}}
\newcommand{\ridcu}{{{\mathsf i}{\downarrow}{\hbox{-}}{\mathsf c}{\uparrow  }}}
%---------------------------------------
\begin{figure}[tbp]
\[
\begin{array}{@{}c@{}c@{}}
%-------------------
\rwdcd\colon\quad\afraise{\atomicflow{
( 0  ,0)*{\afacd{}{}{}\one{}\two};
(-2  ,4)*{\afawdnw{}{}};
(-3.5,0)*{\invisiblemark};
( 3.5,0)*{\invisiblemark}}}
\quad\to\quad
\atomicflow{
( 0  ,3.3)*{\aflabelright{\one,\two}};
( 0  ,3  )*{\afvj6};
(-1.5,0  )*{\invisiblemark};
( 3  ,0  )*{\invisiblemark}}
&\qquad
%-------------------
\rcuwu\colon\quad\aflower{\atomicflowinv{
( 0  ,0)*{\afacu{}{}{}\two{}\one};
(-2  ,6)*{\afawunw{}{}};
( 3.5,0)*{\invisiblemark}}}
\quad\to\quad
\atomicflow{
(0,3  )*{\afvj6};
(0,3.3)*{\aflabelright{\one,\two}};
(3,0  )*{\invisiblemark}}
\\
%-------------------
\rwdiu\colon\quad\atomicflow{
( 0  ,0)*{\afaiu{}{}{}\one{}{}};
(-2  ,4)*{\afawdnw{}{}};
(-3.5,0)*{\invisiblemark};
( 3.5,0)*{\invisiblemark}}
\quad\to\quad
\aflower{\atomicflow{
( 0  ,0)*{\afawu{}{}{}\one};
(-1.5,0)*{\invisiblemark};
( 1.5,0)*{\invisiblemark}}}
&\qquad
%-------------------
\ridwu\colon\quad\atomicflowinv{
( 0  ,0)*{\afaid{}{}{}\one{}{}};
(-2  ,6)*{\afawunw{}{}};
( 3.5,0)*{\invisiblemark}}
\quad\to\quad
\afraise{\atomicflowinv{
(0  ,0)*{\afawd{}{}{}\one};
(1.5,0)*{\invisiblemark}}}
\\
%-------------------
\multispan2{\hfil$
\rwdwu\colon\quad\atomicflow{
( 0  ,6)*{\afawd{}{}{}{}};
( 0  ,0)*{\afawunw{}{}{}{}};
(-1.5,0)*{\invisiblemark};
( 1.5,0)*{\invisiblemark}}
\quad\to\quad
\atomicflow{}
$\hfil}\\
%-------------------
\rwdcu\colon\quad\afraise{\atomicflow{
( 0  ,-4)*{\afacu\one{}{}\two{}{}};
( 0  , 0)*{\afawdnw{}{}};
(-3.5, 0)*{\invisiblemark};
( 3.5, 0)*{\invisiblemark}}}
\quad\to\quad
\atomicflow{
(-2  ,-4)*{\afawd{}{}\one{}};
( 2  ,-4)*{\afawd{}{}{}\two};
(-3.5, 0)*{\invisiblemark};
( 3.5, 0)*{\invisiblemark}}
&\qquad
%-------------------
\rcdwu\colon\quad\aflower{\atomicflowinv{
( 0  ,-4)*{\afacd\one{}{}\two{}{}};
( 0  , 2)*{\afawunw{}{}};
(-3.5, 0)*{\invisiblemark};
( 3.5, 0)*{\invisiblemark}}}
\quad\to\quad
\atomicflowinv{
(-2  ,-4)*{\afawu{}{}\one{}};
( 2  ,-4)*{\afawu{}{}{}\two};
(-3.5, 0)*{\invisiblemark};
( 3.5, 0)*{\invisiblemark}}
\\
%-------------------
\rcdiu\colon\quad\aflower{\atomicflow{
(   6, 6)*{\afvjd4{}\three{}{}};
(   3, 0)*{\afaiuex{}{}{}{}{}{}32};
(   0, 4)*{\afacdnw\one{}{}\two};
(-3.5, 0)*{\invisiblemark};
( 7.5, 0)*{\invisiblemark}}}
\quad\to\quad
\aflower{\atomicflow{
(  10,8)*{\afacu{}{}{}{}{}\three};
(   0,8)*{\afvjd8\one{}};
(   4,8)*{\afvjd8{}\two};
(   6,2)*{\afaiunw{}{}};
(   6,0)*{\afaiuex{}{}{}{}{}{}31}}}
&\qquad
%-------------------
\ridcu\colon\quad\afraise{\atomicflowinv{
(   6,6)*{\afvju4{}\three{}{}};
(   3,0)*{\afaidex {}{}{}{}{}{}32};
(   0,6)*{\afacunw\one{}{}\two};
(-3.5,0)*{\invisiblemark};
( 7.5,0)*{\invisiblemark}}}
\quad\to\quad
\afraise{\atomicflowinv{
(  10,8)*{\afacd{}{}{}{}{}\three};
(   0,8)*{\afvju8\one{}};
(   4,8)*{\afvju8{}\two};
(   6,4)*{\afaidnw{}{}};
(   6,0)*{\afaidex{}{}{}{}{}{}31}}}
\\
%-------------------
\multispan2{\hfil$
\rcdcu\colon\quad\atomicflow{
( 0,6)*{\afacd\one{}{}\two{}{}};
( 0,0)*{\afacunw\three{}{}\four};
(-4,0)*{\invisiblemark};
( 4,0)*{\invisiblemark}}
\quad\to\quad
\atomicflow{
(0,12)*{\afacu{}{}{}{}\one{}};
(6,12)*{\afacu{}{}{}{}{}\two};
( 0,0)*{\afacd{}{}{}{}\three{}};
( 6,0)*{\afacd{}{}{}{}{}\four};
(-2,6)*{\afvj4};
( 8,6)*{\afvj4};
( 3,6)*{\afex24}}
$\hfil}\\
\end{array}
\]
\caption{Atomic-flow reduction rules.}
\label{FigRed}
\end{figure}%

We would like to use the reductions in Figure~\ref{FigRed} as rules for rewriting inside generic atomic flows. To do so, in general, we should have matching upper and lower edges in the flows that participate in the reduction, and the reductions in the figure clearly do so. However, we also have to pay attention to polarities, not to disrupt atomic flows. In fact, consider the following example.

%---------------------------------------
\begin{exa}
The `reduction' on the left, when used inside a larger atomic flow, might create a situation as on the right:
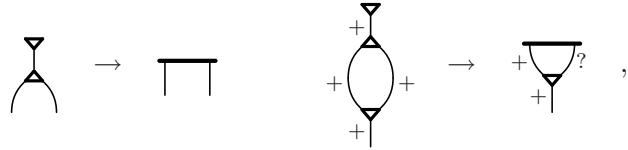
\nopagebreak[4]\medskip\afnegspace
\[
\atomicflow{
( 0  ,0)*{\afacu{}{}{}{}{}{}};
( 0  ,4)*{\afawdnw{}{}}}
\quad\to\quad
\atomicflow{
( 0  ,0)*{\afaid{}{}{}{}{}{}}}
\qquad\qquad
\atomicflow{
( 0  , 8)*{\afacu{}{}{}{}\ppl{}};
( 0  ,12)*{\afawdnw{}{}};
(-2  , 4)*{\aflabelleft\ppl};
( 2  , 4)*{\aflabelright\ppl};
( 0  , 0)*{\afacd{}{}{}{}\ppl{}};
(-3.5, 0)*{\invisiblemark};
( 3.5, 0)*{\invisiblemark}}
\quad\to\quad
\atomicflow{
( 0  ,4)*{\afaidnw{}{}};
( 0  ,0)*{\afacd\ppl{}{}{\scriptstyle?}\ppl{}};
(-3.5,0)*{\invisiblemark};
( 3.5,0)*{\invisiblemark}}
\quad,
\] 
where the graph at the right is not an atomic flow, for lack of a polarity assignment.
\end{exa}

This prompts us to define reduction rules and reductions for atomic flows as follows.

%---------------------------------------
\begin{defi}
An (\emph{atomic-flow}) \emph{reduction rule $r$ from flow $A$ to flow $B$} is a quadruple $(A,B,f,g)$ such that:
\begin{enumerate}
\item $f$ is a one-to-one map from the upper edges of $A$ to the upper edges of $B$,
\item $g$ is a one-to-one map from the lower edges of $A$ to the lower edges of $B$,
\item for every polarity assignment $\pi$ for $A$, there is a polarity assignment $\pi'$ for $B$ such that $\pi'(f(\epsilon))=\pi(\epsilon)$ and $\pi'(g(\epsilon'))=\pi(\epsilon')$, for any upper edge $\epsilon$ and any lower edge $\epsilon'$ of $A$;
\end{enumerate}
we define reduction rules with graphical expressions $r\colon A\to B$, where $f$ and $g$ are indicated by labelling edges. A binary relation $R$ on the set of atomic flows is called an (\emph{atomic-flow}) \emph{reduction} if, whenever $C\mathrel{R}D$, there is a one-to-one map from the upper edges of $C$ to the upper edges of $D$ and a one-to-one map from the lower edges of $C$ to the lower edges of $D$. For every reduction rule $r\colon A\to B$, the reduction ${\to_r}$ is defined, such that $C\to_r D$ if and only if $A$ appears as a subgraph in $C$ and we obtain $D$ by replacing $A$ with $B$ in $C$, while respecting the correspondence of edges; we call this operation a \emph{reduction by $r$}.
\end{defi}

\begin{rem}
The condition on polarity assignments for a reduction rule $r$ guarantees that the $D$ in $C\to_r D$ is a proper atomic flow, if $C$ is one.
\end{rem}

\begin{rem}
Because of the condition on polarity assignments for reduction rules, two distinct connected components in a flow cannot be connected by a reduction. To see that this is impossible, consider the following `reduction rule', which violates the condition on polarity assignments:
\[
\aflower{\atomicflow{
(-2,0)*{\afawu{}{}{}{}};
( 2,0)*{\afawu{}{}{}{}}}}
\quad\to\quad
\aflower{\atomicflow{
( 0  , 2)*{\afaiu{}{}{}{}{}{}}}}
\quad.
\]
\afnegspace
For this `reduction rule' there exist both valid (left) and invalid (right) polarity assignments:
\[
\aflower{\atomicflow{
(-2  ,0)*{\afawu{}{}\ppl{}};
( 2  ,0)*{\afawu{}{}{}\pmi};
(-3.5,0)*{\invisiblemark};
( 3.5,0)*{\invisiblemark}}}
\quad\to\quad
\aflower{\atomicflow{
( 0  ,2)*{\afaiu\ppl{}{}\pmi{}{}};
(-3.5,0)*{\invisiblemark};
( 3.5,0)*{\invisiblemark}}}
\qquad\qquad
\aflower{\atomicflow{
(-2,0)*{\afawu{}{}\ppl{}};
( 2,0)*{\afawu{}{}{}\ppl};
(-3.5,0)*{\invisiblemark};
( 3.5,0)*{\invisiblemark}}}
\quad\to\quad
\aflower{\atomicflow{
( 0  , 2)*{\afaiu\ppl{}{}{\scriptstyle?}{}{}};
(-3.5, 0)*{\invisiblemark};
( 3.5, 0)*{\invisiblemark}}}
\quad.
\]
\afnegspace
\end{rem}

It is immediate to check:

\begin{prop}
The graphical expressions in Figure~\ref{FigRed} are atomic-flow reduction rules.
\end{prop}

Our reduction rules bear a striking resemblance to many rewriting systems on graphs, and in particular with interaction nets \cite{Lafo:97:Interact:jb}. It is certainly possible that some interesting connections with other formalisms can be drawn at some point, but, at this time, we are not aware of any. We note that the resemblance might simply be due to there being very few things that we can do with atoms: we can carry them through, delete them or duplicate them, and it is difficult to think of anything else. We might think of making several copies of an atom at once, instead of just two, and this indeed has some uses in the fight against the bureaucracy related to associativity of (co)contraction.

What is peculiar to our work is the fact that reducing flows corresponds to transforming derivations in a very direct way. The correspondence is captured by the notion of soundness, which we now define.

%---------------------------------------
\begin{defi}
A reduction $R$ is \emph{sound} if, for every $A$ and $B$ such that $A\mathrel{R}B$ and for every derivation $\Phi$ with atomic flow $A$, there is a derivation $\Psi$ with atomic flow $B$ such that $\Phi$ and $\Psi$ have the same premiss and conclusion; in this case we write $\Phi\mathrel{R}\Psi$. A reduction rule $r$ is \emph{sound} if $\to_r$ is sound.
\end{defi}

The proof of the following theorem is essentially contained in Figures~\ref{FigRedW} and \ref{FigRedC}.

%---------------------------------------
\begin{figure}[tbp]
\[
\begin{array}{@{}l@{}c@{}}
%---------------------------------------
\rwdcd\colon\hfil\afraise{\atomicflow{
( 0  ,0)*{\afacd\three{}{}\one{}\two};
(-2  ,4)*{\afawdnw{}{}};
(-3.5,0)*{\invisiblemark};
( 3.5,0)*{\invisiblemark}}}
\quad\to\quad
\atomicflow{
( 0  ,3.3)*{\aflabelright{\one,\two}};
( 0  ,3  )*{\afvj6};
(-1.5,0  )*{\invisiblemark};
( 3  ,0  )*{\invisiblemark}}
&\qquad
\vlderivation                                 {
\vlin{\acd}{}{\zeta\{a^\two\}             }  {
\vlde{\Phi}{}{\zeta\vlscn[a^\three.a^\one]} {
\vlin{\awd}{}{\xi\{a^\three\}             }{
\vlhy        {\xi\{\fff\}                 }}}}}
\quad\to_\rwdcd\quad
\vlderivation                                                      {
\vlin{=                      }{}{\zeta\{a^{\one,\two}\}         } {
\vlde{\Phi\{a^\three\ot\fff\}}{}{\zeta\vlscn[\fff.a^{\one,\two}]}{
\vlhy                           {\xi\{\fff\}                    }}}}
\\
\noalign{\bigskip}
%---------------------------------------
\rwdiu\colon\hfil\atomicflow{
( 0  ,0)*{\afaiu\two {}{}\one{}{}};
(-2  ,4)*{\afawdnw{}{}};
(-3.5,0)*{\invisiblemark};
( 3.5,0)*{\invisiblemark}}
\quad\to\quad
\aflower{\atomicflow{
( 0  ,0)*{\afawu{}{}{}\one};
(-1.5,0)*{\invisiblemark};
( 1.5,0)*{\invisiblemark}}}
&\qquad
\vlderivation                                    {
\vlin{\aiu}{}{\zeta\{\fff\}                  }  {
\vlde{\Phi}{}{\zeta\vlscn(a^\two.\bar a^\one)} {
\vlin{\awd}{}{\xi\{a^\two\}                  }{
\vlhy        {\xi\{\fff\}                    }}}}}
\quad\to_\rwdiu\quad
\vlderivation                                                   {
\vlin{=                    }{}{\zeta\{\fff\}                }  {
\vlin{\awu                 }{}{\zeta\vlscn(\fff.\ttt)       } {
\vlde{\Phi\{a^\two\ot\fff\}}{}{\zeta\vlscn(\fff.\bar a^\one)}{
\vlhy                         {\xi\{\fff\}                  }}}}}
\\
\noalign{\bigskip}
%---------------------------------------
\rwdwu\colon\hfil\atomicflow{
( 0  ,6)*{\afawd{}{}{}\one};
( 0  ,0)*{\afawunw{}{}{}{}};
(-1.5,0)*{\invisiblemark};
( 1.5,0)*{\invisiblemark}}
\quad\to\quad
\atomicflow{}
&\qquad
\vlderivation                    {
\vlin{\awu}{}{\zeta\{\ttt  \}}  {
\vlde{\Phi}{}{\zeta\{a^\one\}} {
\vlin{\awd}{}{\xi\{a^\one\}  }{
\vlhy        {\xi\{\fff  \}  }}}}}
\quad\to_\rwdwu\quad
\vlderivation                                                    {
\vlin{=                    }{}{\zeta\{\ttt\}                }   {
\vlin{\swi                 }{}{\zeta\vlscn[(\fff.\fff).\ttt]}  {
\vlin{=                    }{}{\zeta\vlscn(\fff.[\fff.\ttt])} {
\vlde{\Phi\{a^\one\ot\fff\}}{}{\zeta\{\fff\}                }{
\vlhy                         {\xi\{\fff\}                  }}}}}}
\\
\noalign{\bigskip}
%---------------------------------------
\rwdcu\colon\hfil\afraise{\atomicflow{
( 0  ,-4)*{\afacu\one{}{}\two{}\three};
( 0  , 0)*{\afawdnw{}{}};
(-3.5, 0)*{\invisiblemark};
( 3.5, 0)*{\invisiblemark}}}
\quad\to\quad
\atomicflow{
(-2  ,-4)*{\afawd{}{}\one{}};
( 2  ,-4)*{\afawd{}{}{}\two};
(-3.5, 0)*{\invisiblemark};
( 3.5, 0)*{\invisiblemark}}
&\qquad
\vlderivation                               {
\vlin{\acu}{}{\zeta\vlscn(a^\one.a^\two)}  {
\vlde{\Phi}{}{\zeta\{a^\three\}         } {
\vlin{\awd}{}{\xi\{a^\three\}           }{
\vlhy        {\xi\{\fff\}               }}}}}
\quad\to_\rwdcu\quad
\vlderivation                                                   {
\vlin{\awd                   }{}{\zeta\vlscn(a^\one.a^\two)}   {
\vlin{\awd                   }{}{\zeta\vlscn(a^\one.\fff)  }  {
\vlin{=                      }{}{\zeta\vlscn(\fff.\fff)    } {
\vlde{\Phi\{a^\three\ot\fff\}}{}{\zeta\{\fff\}             }{
\vlhy                           {\xi\{\fff\}               }}}}}}
\\
%---------------------------------------
\end{array}
\]
\caption{`Downwards' reduction rules for weakening and their soundness.}
\label{FigRedW}
\end{figure}%

%---------------------------------------
\begin{figure}[tbp]
\[
\begin{array}{@{}l@{}c@{}}
%---------------------------------------
\rcdiu\colon\aflower{\atomicflow{
(   6, 6)*{\afvjd4{}\three{}{}};
(   3, 0)*{\afaiuex{}\four{}{}{}{}32};
(   0, 4)*{\afacdnw\one{}{}\two};
(-3.5, 0)*{\invisiblemark};
( 7.5, 0)*{\invisiblemark}}}
\quad\to\quad
\aflower{\atomicflow{
(10,8)*{\afacu{}{}{}{}{}\three};
( 0,8)*{\afvjd8\one{}};
( 4,8)*{\afvjd8{}\two};
( 6,2)*{\afaiunw{}{}};
( 6,0)*{\afaiuex{}{}{}{}{}{}31}}}
&\qquad
\vlderivation                                       {
\vlin{\aiu}{}{\zeta\{\fff\}                     }  {
\vlde{\Phi}{}{\zeta\vlscn(a^\four.\bar a^\three)} {
\vlin{\acd}{}{\xi\{a^\four\}                    }{
\vlhy        {\xi\vlscn[a^\one.a^\two]          }}}}}
\quad\to_\rcdiu\quad
\vlderivation                                                     {
\vlin{\aiu}{}{\zeta\{\fff\}                               }      {
\vlin{=   }{}{\zeta\vlscn(a^\one.\bar a)                  }     {
\vlin{\aiu}{}{\zeta\vlscn([\fff.a^\one].\bar a)           }    {
\vlin{\swi}{}{\zeta\vlscn([(\bar a.a^\two).a^\one].\bar a)}   {
\vlin{=   }{}{\zeta\vlscn((\bar a.[a^\two.a^\one]).\bar a)}  {
\vlin{\acu}{}{\zeta\vlscn([a^\one.a^\two].(\bar a.\bar a))} {
\vlde{\Phi\{a^\four\ot\vls[a.a]\}}
           {}{\zeta\vlscn([a^\one.a^\two].\bar a^\three)  }{
\vlhy        {\xi\vlscn[a^\one.a^\two]                    }}}}}}}}}
\\
\noalign{\bigskip}
%---------------------------------------
\rcdcu\colon\hfil\atomicflow{
( 0,6)*{\afacd\one{}{}\two{}\five};
( 0,0)*{\afacunw\three{}{}\four};
(-4,0)*{\invisiblemark};
( 4,0)*{\invisiblemark}}
\quad\to\quad
\atomicflow{
( 0,12)*{\afacu{}{}{}{}\one{}};
( 6,12)*{\afacu{}{}{}{}{}\two};
( 0, 0)*{\afacd{}{}{}{}\three{}};
( 6, 0)*{\afacd{}{}{}{}{}\four};
(-2, 6)*{\afvj4};
( 8, 6)*{\afvj4};
( 3, 6)*{\afex24}}
&\qquad
\vlderivation                                  {
\vlin{\acu}{}{\zeta\vlscn(a^\three.a^\four)}  {
\vlde{\Phi}{}{\zeta\{a^\five\}             } {
\vlin{\acd}{}{\xi\{a^\five\}               }{
\vlhy        {\xi\vlscn[a^\one.a^\two]     }}}}}
\quad\to_\rcdcu\quad
\vlderivation                                                            {
\vlin{\acd                       }{}{\zeta\vlscn(a^\three.a^\four)}     {
\vlin{\acd                       }{}{\zeta\vlscn(a^\three.[a.a])  }    {
\vlin{\med                       }{}{\zeta\vlscn([a.a].[a.a])     }   {
\vlin{\acu                       }{}{\zeta\vlscn[(a.a).(a.a)]     }  {
\vlin{\acu                       }{}{\zeta\vlscn[a^\one.(a.a)]    } {
\vlde{\Phi\{a^\five\ot\vls[a.a]\}}{}{\zeta\vlscn[a^\one.a^\two]   }{
\vlhy                               {\xi\vlscn[a^\one.a^\two]     }}}}}}}}
\\
%---------------------------------------
\end{array}
\]
\caption{`Downwards' reduction rules for contraction and their soundness.}
\label{FigRedC}
\end{figure}%

%---------------------------------------
\begin{thm}\label{TheoSound}
The reduction rules\/ $\rwdcd$, $\rwdiu$, $\rwdwu $, $\rwdcu$, $\rcdiu$, $\rcdcu$, $\rcuwu$, $\ridwu$, $\rcdwu$ and\/ $\ridcu$ are sound.
\end{thm}

\proof
For $r\in\{\rwdcd,\rwdiu,\rwdwu,\rwdcu,\rcdiu,\rcdcu\}$ and $r\colon A\to B$ as in the left columns of Figures~\ref{FigRedW} and \ref{FigRedC}, for every $C$ and $D$ such that $C\to_r D$ and for every $\Psi$ with flow $C$, the right columns of the tables provide reductions $\Psi\to_r\Psi'$, where $\Psi'$ has flow $D$, as follows. If $\Phi'\to_r\Phi''$ is the reduction provided by the table, then
\[
\Psi=
\vlderivation              {
\vlde{\Psi_2}{}{\beta  }  {
\vlde{\Phi' }{}{\beta' } {
\vlde{\Psi_1}{}{\alpha'}{
\vlhy          {\alpha }}}}}
\qquad\hbox{and}\qquad
\Psi'=
\vlderivation              {
\vlde{\Psi_2}{}{\beta  }  {
\vlde{\Phi''}{}{\beta' } {
\vlde{\Psi_1}{}{\alpha'}{
\vlhy          {\alpha }}}}}
\quad.
\]
We can deal with the remaining rules by employing dual derivations to the ones shown.
\qed

%---------------------------------------
\begin{rem}\label{RemIndep}
The previous soundness theorem only depends on the switch and medial rules for the reductions in Figure~\ref{FigRedC}. Any system obtained from $\SKS$ by replacing $\swi$ and $\med$ with linear rules that can derive them would support a soundness theorem like the one above, for the same reduction rules. For example, we could think of replacing $\swi$ with the rule $\vlinf{\swi'}{}{\vls[(a.c).[b.d]]}{\vls([a.b].[c.d])}$, from which $\swi$ is derivable.
\end{rem}

%---------------------------------------
\begin{defi}
A finite set of reduction rules is a \emph{flow rewriting system}. For every flow rewriting system $F=\{r_1,\dots,r_h\}$ we define ${\to_F}={\to_{r_1}\cup\cdots\cup{\to_{r_h}}}$. The reflexive transitive closure of $\to_F$ is denoted by $\to_F^\star$. Given a set of atomic flows $S$, we say that a flow rewriting system $F$ is \emph{terminating on $S$} if there is no infinite chain $A_1\to_F A_2\to_F\cdots$, for every $A_1\in S$; if $F$ is terminating on the set of atomic flows, we say that it is \emph{terminating}. We say that atomic flow $A$ is \emph{normal} for flow rewriting system $F$ if there is no atomic flow $B$ such that $A\to_F B$.
\end{defi}

%-------------------------------------------------------------------------------
\subsection{Weakening and Coweakening}\label{SubsWeak}

The reduction rules for weakening and coweakening make for a very simple flow rewriting system. They are very `friendly' rules, because they greatly simplify atomic flows and associated derivations.

\newcommand{\frw}{{\mathsf w}}
%---------------------------------------
\begin{defi}
The following flow rewriting system is called $\frw$:
\[
\{\;\rwdcd\;,\;\rcuwu\;,\;\rwdiu\;,\;\ridwu\;,\;\rwdwu\;,\;\rwdcu\;,\;\rcdwu\;\}
\quad.
\]
\end{defi}

%---------------------------------------
\begin{thm}\label{TheoWTerm}
Flow rewriting system\/ $\frw$ is terminating.
\end{thm}

%---------------------------------------
\proof
At every reduction, either the number of vertices decreases, or it stays the same but the number of contraction and cocontraction vertices decreases.
\qed

%---------------------------------------
\begin{rem}\label{RemNDCW}
If flow $A$ is normal for $\frw$, then there is no $\ai$-path from a weakening or coweakening vertex to another vertex in $A$.
\end{rem}

Since reducing by $\frw$ does not introduce new edges, we have:

%---------------------------------------
\begin{prop}\label{PropWLF}
If $A$ is cycle-free and $A\to_\frw^\star B$ then $B$ is cycle-free.
\end{prop}

%-------------------------------------------------------------------------------
\subsection{Contraction and Cocontraction}\label{SubsContr}

The reduction rules for contraction and cocontraction are much less `friendly' than weakening/coweakening ones, mainly because they create infinite reduction chains. A judicious use of these rules is the key to success for our normalisation methods.

\newcommand{\frc}{{\mathsf c}}
%---------------------------------------
\begin{defi}
The following flow rewriting system is called $\frc$:
\[
\{\;\rcdiu\;,\;\ridcu\;,\;\rcdcu\;\}\quad.
\]
\end{defi}

%---------------------------------------
\begin{rem}\label{RemCycle}
Flow rewriting system $\frc$ is not terminating:
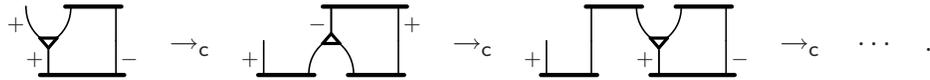
\nopagebreak[4]\medskip\afnegspace
\[
\atomicflow{
( 4  ,8)*{\afaidnw{}{}};
( 6  ,6)*{\afvj4};
( 3  ,0)*{\afaiuex{\ppl}{}{}{\pmi}{}{}32};
( 0  ,4)*{\afacd{\ppl}{}{}{}{}{}};
(-3.5,0)*{\invisiblemark};
( 7.5,0)*{\invisiblemark}}
\quad\to_\frc\quad
\atomicflow{
( 9  ,12)*{\afaidex{\pmi}{}{}{\ppl}{}{}32};
( 6  , 8)*{\afacu{}{}{}{}{}{}};
(12  , 6)*{\afvj4};
( 0  , 6)*{\afvjd4\ppl{}};
( 2  , 2)*{\afaiunw{}{}};
(10  , 2)*{\afaiunw{}{}{}{}{}{}};
(-1.5, 0)*{\invisiblemark};
(13.5, 0)*{\invisiblemark}}
\quad\to_\frc\quad
\atomicflow{
( 6  ,8)*{\afaidnw{}{}};
(14  ,8)*{\afaidnw{}{}};
( 4  ,6)*{\afvj4};
(16  ,6)*{\afvj4};
( 2  ,0)*{\afaiu{\ppl}{}{}{}{}{}};
(13  ,0)*{\afaiuex{\ppl}{}{}{\pmi}{}{}32};
(10  ,4)*{\afacd{}{}{}{}{}{}};
(-1.5,0)*{\invisiblemark};
(17.5,0)*{\invisiblemark}}
\quad\to_\frc\quad\cdots\quad.
\]
\afnegspace
We see that if a contraction vertex belongs to an $\ai$-cycle, reductions by $\frc$ make it `bounce' in the $\ai$-cycle and create a trail; while bouncing, the vertex alternates between contraction and cocontraction; if we assign a polarity to the flow, the vertex alternates between being positive and negative.
\end{rem}

Through a simple argument by contradiction, we have:

%---------------------------------------
\begin{prop}\label{PropCLF}
If $A$ is cycle-free and $A\to_\frc^\star B$ then $B$ is cycle-free.
\end{prop}

Again, reasoning by contradiction, we have:

%---------------------------------------
\begin{prop}\label{PropCOW}
If an atomic flow is normal for\/ $\frc$ then all its\/ $\ai$-paths are clean paths.
\end{prop}

The previous proposition could be rephrased by saying that if an atomic flow is normal for $\frc$ then all its $\ai$-connections are simple edges.

Since reducing by $\frw$ does not introduce new vertices, we have:

%---------------------------------------
\begin{prop}\label{PropWC}
If $A$ is normal for\/ $\frc$ and $A\to_\frw^\star B$ then $B$ is normal for\/ $\frc$.
\end{prop}

By contradiction and a simple case analysis, we have:

%---------------------------------------
\begin{prop}\label{PropCW}
If $A$ is normal for\/ $\frw$ and $A\to_\frc^\star B$ then $B$ is normal for\/ $\frw$.
\end{prop}

Maximal $\ai$-paths provide for a measure when dealing with the termination of $\frc$.

%---------------------------------------
\begin{rem}\label{RemRank}
A simple inspection to the reduction rules of $\frc$ convinces us that reducing by $\frc$ does not change the number and length of the maximal $\ai$-paths of a flow. The same holds for the maximal $\ai$-paths to or from vertices that are not involved in a given reduction.
\end{rem}

%---------------------------------------
\begin{thm}\label{TheoCTerm}
Flow rewriting system\/ $\frc$ is terminating on the set of cycle-free atomic flows.
\end{thm}

\newcommand{\rk}{{\mathsf r}}
%-------------------
\proof
Let $A$ be a cycle-free flow. We associate to each contraction (resp., cocontraction) vertex $\nu$ its \emph{rank} $\rk_\nu=\sum_{p_i\in I_\nu}h_i$, where $I_\nu$ is the set of all maximal $\ai$-paths $p_i=\epsilon^i_1,\dots,\epsilon^i_{h_i}$ from $\nu$, such that $\epsilon^i_1$ is the lower (resp., upper) edge of $\nu$ (so, the rank of a vertex is the sum of the lengths of certain maximal $\ai$-paths from it). Note that every (co)contraction vertex has non-zero rank. We prove that a reduction of $A$ by $\frc$ decreases the sum of the ranks of the (co)contraction vertices of $A$. First note that the rank of the vertices not involved in the reduction step stays the same (see Remark~\ref{RemRank}). We then need to show that the sum of the ranks decreases for the vertices involved. There are three cases, depending on the reduction rule:
\begin{itemize}
%-------------------
\item[$\rcdiu$:]a contraction vertex $\nu$ is replaced by a cocontraction vertex $\nu'$, and $\rk_{\nu'}=\rk_\nu-n$, where $n>0$ is the number of maximal $\ai$-paths from $\nu$ whose first edge is the lower edge of $\nu$;
%-------------------
\item[$\ridcu$:]this is dual to the previous case;
%-------------------
\item[$\rcdcu$:]a contraction vertex $\nu$ and a cocontraction vertex $\nu'$ are replaced by two contraction vertices $\nu_1$ and $\nu_2$ and two cocontraction vertices $\nu'_1$ and $\nu'_2$; we have $\rk_{\nu_1}+\rk_{\nu_2}=\rk_\nu-n$, where $n>0$ is the number of maximal $\ai$-paths from $\nu$ whose first edge is the lower edge of $\nu$; analogously, we have $\rk_{\nu'_1}+\rk_{\nu'_2}=\rk_{\nu'}-n'$, where $n'>0$ is the number of maximal $\ai$-paths from $\nu'$ whose first edge is the upper edge of $\nu'$.
\qed
\end{itemize}

%---------------------------------------
\begin{rem}\label{RemExpC}
Normalising by $\frc$ can blow the size of atomic flows exponentially, in particular in a situation like the following (noted by Lutz Stra{\ss}burger):
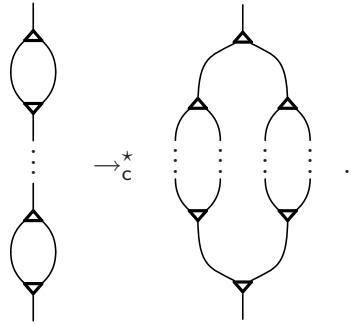
\[
\atomicflow{
(0,29  )*{\afacu{}{}{}{}{}{}};
(0,21  )*{\afacd{}{}{}{}{}{}};
(0,15.8)*{\vdots};
(0, 8  )*{\afacu{}{}{}{}{}{}};
(0, 0  )*{\afacd{}{}{}{}{}{}}}
\quad\to_\frc^\star\quad
\atomicflow{
( 0,29  )*{\afacuexsq{}{}{}{}{}{}21};
(-4,21  )*{\afacunw{}{}{}{}{}{}};
( 4,21  )*{\afacunw{}{}{}{}{}{}};
(-6,17.8)*{\vdots};
(-2,17.8)*{\vdots};
( 2,17.8)*{\vdots};
( 6,17.8)*{\vdots};
( 0, 4  )*{\afacdexsq{}{}{}{}{}{}21};
(-4,10  )*{\afacdnw{}{}{}{}{}{}};
( 4,10  )*{\afacdnw{}{}{}{}{}{}}}
\quad.
\]
In fact, if there are $n$ couples cocontraction/contraction like the two shown above on the left, then there are $2^n$ maximal $\ai$-paths, and their number (and length) is conserved by $\to_\frc^\star$ (see Remark~\ref{RemRank}). Exactly one $\ai$-path passes through each edge in the middle portion of the flow on the right.
\end{rem}

Normalising flows via the rewriting system $\frc$ is important, and we know from Remark~\ref{RemCycle} and Theorem~\ref{TheoCTerm} that we can only have termination for cycle-free atomic flows. It turns out that we can `break' $\ai$-cycles if we manage to move contractions and cocontractions away from at least one $\ai$-connection per $\ai$-cycle, so to create a simple edge. As we saw in Definition~\ref{DefSimple}, we call the $\ai$-cycles exhibiting this property `fragile'. We now see how to transform every $\ai$-cycle into a fragile cycle, and Subsection~\ref{SubsBFCECP} shows how to break fragile cycles.

%---------------------------------------
\begin{thm}\label{TheoCycleFragile}
For every atomic flow $A$ there is an atomic flow $B$ such that $A\to_\frc^\star B$ and all the\/ $\ai$-cycles in $B$ are fragile cycles.
\end{thm}

%---------------------------------------
\proof
Given any polarity assignment for $A$, consider all contractions and cocontractions with some edge belonging to an $\ai$-cycle and mapping to $\pmi$. Apply the rules of $\frc$ to these vertices until they change polarity. Now, every $\ai$-cycle contains at least one $\ai$-connection mapping to $\pmi$, and this cannot be anything else than a simple edge.
\qed

%---------------------------------------
\begin{exa}
The negative $\ai$-connection $\one,\two,\three$ can be made into a negative simple edge:
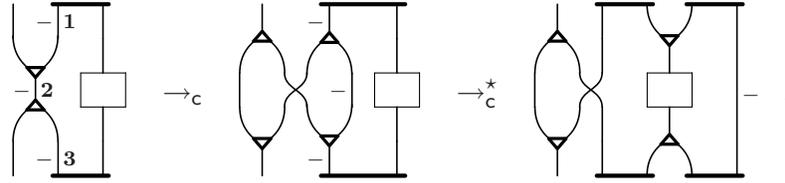
\nopagebreak[4]\vskip-12pt
\nopagebreak[4]\medskip\afnegspace
\[
\atomicflow{
(-4,  8)*{\afvj4};
( 2, 10)*{\afaid{\pmi}{\one}{}{}{}{}};
(-2, -2)*{\afacu{}{}{}{}{\pmi}{\two}};
(-2,  2)*{\afacdnw{}{}{}{}};
(-4, -8)*{\afvj4};
( 2,-10)*{\afaiu{\pmi}{\three}{}{}{}{}};
( 4,  4)*{\afvj4};
( 4,  0)*{\affr44};
( 4, -4)*{\afvj4}}
\quad\to_\frc\quad
\atomicflow{
( 3, 10)*{\afaidex{\pmi}{}{}{}{}{}32};
( 0, 4)*{\afacunw{}{}{}{}};
(-6, 6)*{\afacu{}{}{}{}{}{}};
(-6,-6)*{\afacd{}{}{}{}{}{}};
(-8, 0)*{\afvj4};
( 2, 0)*{\afvju4{\pmi}{}};
(-3, 0)*{\afex24};
(3,-10)*{\afaiuex{\pmi}{}{}{}{}{}32};
(0, -6)*{\afacdnw{}{}{}{}};
(6,  4)*{\afvj4};
(6,  0)*{\affr44};
(6, -4)*{\afvj4}}
\quad\to_\frc^\star\quad
\atomicflow{
(-4,  6)*{\afacu{}{}{}{}{}{}};
(-4, -6)*{\afacd{}{}{}{}{}{}};
(-6,  0)*{\afvj4};
(-1,  0)*{\afex24};
( 0,  4)*{\afvj4};
( 0, -4)*{\afvj4};
( 2, 10)*{\afaidnw{}{}{}{}};
(10, 10)*{\afaidnw{}{}{}{}};
( 0,  8)*{\afvj4};
( 6,  6)*{\afacd{}{}{}{}{}{}};
( 6,  0)*{\affr44};
(12,  0)*{\afvj{20}};
(12,  0)*{\aflabelright{\pmi}};
( 6, -6)*{\afacu{}{}{}{}{}{}};
( 0, -8)*{\afvj4};
( 2,-12)*{\afaiunw{}{}{}{}};
(10,-12)*{\afaiunw{}{}{}{}};
(13.5,0)*{\invisiblemark}}
\quad.
\]
\afnegspace
\end{exa}

%---------------------------------------
\begin{exa}
The following atomic flow reduction (also present in Figure~\ref{FigExStrRed} on page \pageref{FigExStrRed}) shows another application of the previous theorem:
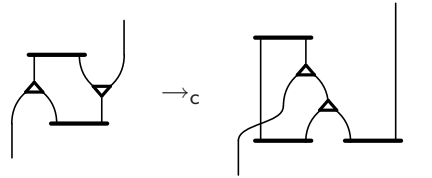
\[
\atomicflow{
( 5,12)*{\afvj4};
( 3, 6)*{\afacd{}{}{}{}{}{}};
(-3, 6)*{\afacu{}{}{}{}{}{}};
( 1, 0)*{\afaiunw{}{}};
(-5, 0)*{\afvj4};
(-1,10)*{\afaidnw{}{}}}
\quad\to_\frc\quad
\atomicflow{
( -2,10)*{\afacu{}{}{}{}{}{}};
( -6, 8)*{\afvj{12}};
(  6,10)*{\afvj{16}};
( -6, 4)*{\afcjr44};
(  0, 4)*{\afacunw{}{}{}{}};
(  4, 0)*{\afaiunw{}{}};
( -4, 0)*{\afaiunw{}{}};
( -8, 0)*{\afvj4};
( -4,14)*{\afaidnw{}{}}}
\quad.
\]
\end{exa}

%===============================================================================
\section{Streamlining Algorithms}\label{SectStreamAlg}

In this section we prove our main results. The basis of our normalisation technique is the elimination of a simple edge, which means eliminating also the axiom and cut it connects. This entails a duplication of the entire flow/derivation, and the repeated process generates an exponential growth in size. Of course, this is what we expect from an algorithm that generalises cut elimination for propositional logic.

We start by studying a single elimination of a simple edge. We will then study two algorithms, based on simple-edge elimination, that, in combination, achieve the desired normalisation result. The final part of this section presents two normalisation algorithms, of varying strength, and discusses possible variations.

%-------------------------------------------------------------------------------
\subsection{Elimination of a Simple Edge}

Definition~\ref{DefRedS} contains a graphical representation of the main idea of this section. It is a reduction of atomic flows and associated derivations, whose purpose is eliminating one simple edge, and at the same time removing one interaction and one cointeraction vertex. Contrary to previous reduction rules, we are dealing here with a global transformation of atomic flows and associated derivations: a reduction can involve an entire flow/derivation and not just a local subflow/subderivation.

The reduction of simple edges that we are about to define is not strictly necessary for getting our results, because we could jump directly to the more general Definition~\ref{DefCBFlow}. However, it is important to understand this reduction in isolation, both because it is simpler, and also because it provides the basis for different reduction strategies from those that we discuss in this paper.

\newcommand{\frse}{{\mathsf{se}}}
%---------------------------------------
\begin{defi}\label{DefRedS}
We define the reduction $\to_\frse$ (where $\frse$ stands for \emph{simple edge}) as follows, for every atomic flow $A$:
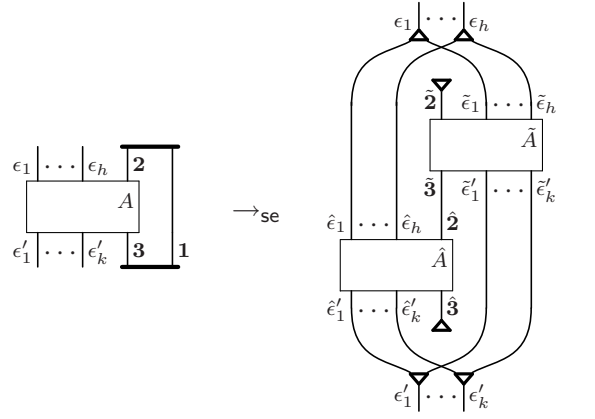
\[
\atomicflow{
(10,16)*{\afaidnw{}{}};
( 0,14)*{\afvju4{\epsilon_1}{}};
( 2,14)*{\cdots};
( 4,14)*{\afvju4{}{\epsilon_h}};
( 8,14)*{\afvju4{}\two};
( 9,10)*{\aflabelleft A};
( 4, 9)*{\affr{10}6};
(12,11)*{\afvj{10}};
( 0, 4)*{\afvjd4{\epsilon'_1}{}};
( 2, 4)*{\cdots};
( 4, 4)*{\afvjd4{}{\epsilon'_k}};
(10, 2)*{\afaiu{}\three{}\one{}{}};
(-3, 0)*{\invisiblemark};
(14, 0)*{\invisiblemark}}
\quad\to_\frse\quad
\atomicflow{
( 0,24   )*{\afvju{16}{\hat\epsilon_1}{}};
( 2,18   )*{\cdots};
( 4,24   )*{\afvju{16}{}{\hat\epsilon_h}};
( 8,18.25)*{\aflabelright{\hat\two}};
( 9,14   )*{\aflabelleft {\hat A}};
( 4,13   )*{\affr{10}6};
( 0, 8.4 )*{\afvjd2{\hat\epsilon'_1}{}};
( 2, 8   )*{\cdots};
( 4, 8.4 )*{\afvjd2{}{\hat\epsilon'_k}};
(12,32   )*{\afvju2{\tilde\epsilon_1}{}};
(14,32   )*{\cdots};
(16,32   )*{\afvju2{}{\tilde\epsilon_h}};
(17,28   )*{\aflabelleft {\tilde A}};
(12,27   )*{\affr{10}6};
(12,16   )*{\afvjd{16}{\tilde\epsilon'_1}{}};
(14,22   )*{\cdots}; 
(16,16   )*{\afvjd{16}{}{\tilde\epsilon'_k}};
( 8,22.25)*{\aflabelleft{\tilde\three}};
( 8,34   )*{\afawd{}{}{\tilde\two}{}};
( 8,42   )*{\cdots};
( 6,38   )*{\afacuexsq{}{}{}{}{\epsilon_1}{}31};
(10,38   )*{\afacuexsq{}{}{}{}{}{\epsilon_h}31};
( 8,20   )*{\afvj8};
( 8, 6   )*{\afawu{}{}{}{\hat\three}};
( 6, 2   )*{\afacdexsq{}{}{}{}{\epsilon'_1}{}31};
(10, 2   )*{\afacdexsq{}{}{}{}{}{\epsilon'_k}31};
( 8,-2   )*{\cdots};
(-3, 0   )*{\invisiblemark};
(19, 0   )*{\invisiblemark}}
\quad,
\]
where $h,k\ge0$, edges have been renamed with $\hat{\enspace}$ and $\tilde{\enspace}$ accents, flows $\tilde A$ and $\hat A$ are both isomorphic to $A$, and edges $\hat\two$ and $\tilde\three$ are identified.
\end{defi}

A simple inspection of the definition of $\to_\frse$ suffices to prove the following statement, about $\to_\frse$ not introducing any $\ai$-cycles.

%---------------------------------------
\begin{prop}
If atomic flow $B$ is cycle-free and $B\to_\frse C$, then $C$ is cycle-free.
\end{prop}

%---------------------------------------
\begin{thm}\label{ThSESound}
Reduction\/ $\to_\frse$ is sound.
\end{thm}

%---------------------------------------
\proof
Let $\Phi$ be a derivation with flow $B$, such that $B\to_\frse C$. We show that there exists a derivation $\Psi$ with flow $C$ and with the same premiss and conclusion as $\Phi$. In the following, we refer to the figure in Definition~\ref{DefRedS}. We assume that $\Phi$ has premiss $\xi\{\mk\ttt\}$ and conclusion $\zeta\{\mk\fff\}$, where the evidenced and labelled $\mk\ttt$ and $\mk\fff$ can be traced to the interaction and cointeraction vertices eliminated by $\to_\frse$, respectively (this can always be done by using switches and unit equations). Intuitively, we can think of $\mk\ttt$ and $\mk\fff$ as mapping to special `unit edges', which can be substituted just like normal edges. So, we assume that $\Phi$ is
\[
\vlderivation                                           {
\vlde{\Phi_3}{}{\zeta\{\mk \fff\}                 }    {
\vlin{\aiu  }{}{\zeta'\{\mk \fff\}                }   {
\vlde{\Phi_2}{}{\zeta'\vlscn(\bar a^\three.a^\one)}  {
\vlin{\aid  }{}{\xi'\vlscn[\bar a^\two.a^\one]    } {
\vlde{\Phi_1}{}{\xi'\{\mk \ttt\}                  }{
\vlhy          {\xi\{\mk \ttt\}                   }}}}}}}
\quad.
\]
We obtain the two derivations $\Psi'$ and $\Psi''$ from $\Phi$ as follows:
\[
\Psi'=\;\;
\vlderivation                                                  {
\vlde{\Phi_3\{\mk\fff\ot\bar a\}}
             {}{\zeta\{\bar a^{\tilde\three}\}          }     {
\vlin{=     }{}{\zeta'\{\bar a^{\tilde\three}\}         }    {
\vlde{\Phi_2\{a^\one\ot\ttt\}}
             {}{\zeta'\vlscn(\bar a^{\tilde\three}.\ttt)}   {
\vlin{\awd  }{}{\xi'\vlscn[\bar a^{\tilde\two}.\ttt]    }  {
\vlin{=     }{}{\xi'\vlscn[\fff.\ttt]                   } {
\vlde{\Phi_1}{}{\xi'\{\ttt\}                            }{
\vlhy          {\xi\{\ttt\}                             }}}}}}}}
\qquad\hbox{and}\qquad
\Psi''=\;\;
\vlderivation                                                                  {
\vlde{\Phi_3                  }{}{\zeta\{\fff\}                         }     {
\vlin{=                       }{}{\zeta'\{\fff\}                        }    {
\vlin{\awu                    }{}{\zeta'\vlscn(\ttt.\fff)               }   {
\vlde{\Phi_2\{a^\one\ot \fff\}}{}{\zeta'\vlscn(\bar a^{\hat\three}.\fff)}  {
\vlin{=                       }{}{\xi'\vlscn[\bar a^{\hat\two}.\fff]    } {
\vlde{\Phi_1\{\mk\ttt\ot\bar a\}}
                               {}{\xi'\{\bar a^{\hat\two}\}             }{
\vlhy                            {\xi\{\bar a^{\hat\two}\}              }}}}}}}}
\quad.
\]
Derivation $\Psi'$ has flow $B'$ and $\Psi''$ has $B''$:
\nopagebreak[4]\bigskip\afnegspace
\[
B'=\raise2\atflowelheight\hbox{$
\atomicflow{
( 4,14)*{\afvju4{\tilde\epsilon_1}{}};
( 6,14)*{\cdots};
( 8,14)*{\afvju4{}{\tilde\epsilon_h}};
( 0, 4)*{\afvjd4{\tilde\three}{}};
( 9,10)*{\aflabelleft {\tilde A}};
( 4, 9)*{\affr{10}6};
( 4, 4)*{\afvjd4{\tilde\epsilon'_1}{}};
( 6, 4)*{\cdots};
( 8, 4)*{\afvjd4{}{\tilde\epsilon'_k}};
( 0,16)*{\afawd{}{}{\tilde\two}{}};
(-2, 4)*{\invisiblemark};
(11, 4)*{\invisiblemark}}$}
\qquad\hbox{and}\qquad
B''=\lower2\atflowelheight\hbox{$
\atomicflow{
( 0,14)*{\afvju4{\hat\epsilon_1}{}};
( 2,14)*{\cdots};
( 4,14)*{\afvju4{}{\hat\epsilon_h}};
( 8,14)*{\afvju4{}{\hat\two}};
( 9,10)*{\aflabelleft {\hat A}};
( 4, 9)*{\affr{10}6}; 
( 0, 4)*{\afvjd4{\hat\epsilon'_1}{}};
( 2, 4)*{\cdots};
( 4, 4)*{\afvjd4{}{\hat\epsilon'_k}};
( 8, 2)*{\afawu{}{}{}{\hat\three}};
(-3, 4)*{\invisiblemark}}$}
\quad.
\]
\afnegspace
We combine $\Psi'$ and $\Psi''$ to get the desired derivation $\Psi$ with flow $C$ and the same premiss and conclusion as $\Phi$:
\[
\Psi=\;\;
\vlderivation                                                              {
\vlin{\cod                      }{}{\zeta\{\fff\}                    }    {
\vlde{\vls[\Psi''.\zeta\{\fff\}]}{}{\vls[\zeta\{\fff\}.\zeta\{\fff\}]}   {
\vlin{\sus                      }{}{\vls[\xi\{\bar a\}.\zeta\{\fff\}]}  {
\vlde{\vls(\xi\{\ttt\}.\Psi')   }{}{\vls(\xi\{\ttt\}.\zeta\{\bar a\})} {
\vlin{\cou                      }{}{\vls(\xi\{\ttt\}.\xi\{\ttt\})    }{
\vlhy                              {\xi\{\ttt\}                      }}}}}}}
\quad,
\]
where $\sus$, $\cod$ and $\cou$ are `macro' rules introduced in Remarks~\ref{RemSupSwitch} and \ref{RemGenContr}.
\qed

The reduction $\to_\frse$ on atomic flows is symmetric and deterministic, and this might be a surprise, because we tend to expect a non-deterministic choice in the elimination of a cut in classical propositional logic. The corresponding construction on derivations might appear symmetric and deterministic, too, but, subtly, it is not. In fact, we have two ways of realising the macro rule $\vldownsmash{\vlinf{\sus}{}{\vls[\xi\{\alpha\}.\zeta\{\fff\}]}{\vls(\xi\{\ttt\}.\zeta\{\alpha\})}}$ by using $\swi$; informally, either we `put $\xi\{\ttt\}$ inside $\zeta\vlhole$' or we `put $\zeta\{\alpha\}$ inside $\xi\vlhole$' (see Remark~\ref{RemSupSwitch}). Note, however, that $\vlinf{\swi}{}{\vls[(\beta.\alpha).\gamma]}{\vls(\beta.[\alpha.\gamma])}$ is a special case of $\sus$ (modulo $=$); so, we might replace $\swi$ with $\sus$ inside $\SKS$, and we would eliminate this non-determinism. Contrary to $\swi$, the rule $\sus$ is not local (in the sense that it cannot be checked in time bounded by a constant); we do not know whether a local rule exists that could make for a symmetric and deterministic derivation corresponding to flows obtained by $\to_\frse$.

Repeated $\to_\frse$ reductions can lead to infinite chains, because at each step a flow is duplicated. In fact, if the flow contains more than one simple edge the number of simple edges after each $\to_\frse$ reduction increases.

%-------------------------------------------------------------------------------
\subsection{Breaking Fragile Cycles and Eliminating Clean Paths}\label{SubsBFCECP}

Since indiscriminately composing $\to_\frse$ reductions might lead to infinite chains, in order to normalise proofs, we have to impose some discipline on $\to_\frse$ reductions. Moreover, it is not guaranteed that $\to_\frse$ can achieve normalisation without preliminarily `preparing' the atomic flow. All that said, we are interested in as free as possible reduction mechanisms.

It turns out that we can achieve our objectives by defining two similar, recursive reductions based on $\to_\frse$. The main idea is to constrain $\to_\frse$ into a very simple binary recursion scheme, controlled by two different conditions on simple edges. These conditions provide induction measures that basically amount to counting the number of simple edges in the flow to reduce. As we are about to see, we are able to keep the symmetry already exhibited by $\to_\frse$.

\newcommand{\frbc}{{\mathsf{bc}}}
\newcommand{\frex}{{\mathsf{ex}}}
The following two reductions, $\to_\frbc$ and $\to_\frex$, should be understood as `one-shot' reductions, meaning that they will only be needed once each in the normalisation process. They consist, basically, of chains of $\to_\frse$ reductions. In order to define $\to_\frex$, we need the notion of `extremal simple edge'.

%---------------------------------------
\begin{defi}
If $\epsilon_1,\dots,\epsilon_k,\dots,\epsilon_h$ is a maximal clean path (\emph{i.e.}, a clean path that is also a maximal $\ai$-path), $\epsilon_k$ is a simple edge and the edges $\epsilon_{k+1}$, \dots, $\epsilon_h$ are not simple edges, then $\epsilon_k$ is an \emph{extremal simple edge}.
\end{defi}

%---------------------------------------
\begin{rem}
In every maximal clean path there are at most two extremal simple edges, one for each `direction of the path'. In Example~\ref{ExExtrem} (page~\pageref{ExExtrem}) there are two maximal clean paths, and edges $\one$ and $\three$ are the extremal simple edges of both.
\end{rem}

%---------------------------------------
\begin{defi}\label{DefCBFlow}
We inductively define the reductions $\to_\frbc$ and $\to_\frex$ (where $\frbc$ stands for \emph{break\/ $\ai$-cycles} and $\frex$ for \emph{eliminate extremal simple edges}) as follows. Given a flow $B$, the base cases are:
\begin{itemize}
\item if there are no fragile cycles in $B$ then $B\to_\frbc B$;
\item if there are no extremal simple edges in $B$ then $B\to_\frex B$.
\end{itemize}
For the inductive cases, let us suppose that 
\[
B=\atomicflow{
(10,16)*{\afaidnw{}{}};
( 0,14)*{\afvju4{\epsilon_1}{}};
( 2,14)*{\cdots};
( 4,14)*{\afvju4{}{\epsilon_h}};
( 8,14)*{\afvju4{}\two};
( 9,10)*{\aflabelleft A};
( 4, 9)*{\affr{10}6};
(12,11)*{\afvj{10}};
( 0, 4)*{\afvjd4{\epsilon'_1}{}};
( 2, 4)*{\cdots};
( 4, 4)*{\afvjd4{}{\epsilon'_k}};
(10, 2)*{\afaiu{}\three{}\one{}{}};
(-3, 0)*{\invisiblemark};
(14, 0)*{\invisiblemark}}
\qquad\hbox{and}\qquad
C=
\atomicflow{
( 0,24   )*{\afvju{16}{\hat\epsilon_1}{}};
( 2,18   )*{\cdots};
( 4,24   )*{\afvju{16}{}{\hat\epsilon_h}};
( 8,18.25)*{\aflabelright{\hat\two}};
( 9,14   )*{\aflabelleft{D''}};
( 4,13   )*{\affr{10}6};
( 0, 8.4 )*{\afvjd2{\hat\epsilon'_1}{}};
( 2, 8   )*{\cdots};
( 4, 8.4 )*{\afvjd2{}{\hat\epsilon'_k}};
(12,32   )*{\afvju2{\tilde\epsilon_1}{}};
(14,32   )*{\cdots};
(16,32   )*{\afvju2{}{\tilde\epsilon_h}};
(17,28   )*{\aflabelleft{D'}};
(12,27   )*{\affr{10}6};
(12,16   )*{\afvjd{16}{\tilde\epsilon'_1}{}};
(14,22   )*{\cdots}; 
(16,16   )*{\afvjd{16}{}{\tilde\epsilon'_k}};
( 8,22.25)*{\aflabelleft{\tilde\three}};
( 8,42   )*{\cdots};
( 6,38   )*{\afacuexsq{}{}{}{}{\epsilon_1}{}31};
(10,38   )*{\afacuexsq{}{}{}{}{}{\epsilon_h}31};
( 8,20   )*{\afvj8};
( 6, 2   )*{\afacdexsq{}{}{}{}{\epsilon'_1}{}31};
(10, 2   )*{\afacdexsq{}{}{}{}{}{\epsilon'_k}31};
( 8,-2   )*{\cdots};
(-3, 0   )*{\invisiblemark};
(19, 0   )*{\invisiblemark}}
\quad,
\]
where $h,k\ge0$, and let
\nopagebreak[4]\medskip\afnegspace
\[
B'=\raise2\atflowelheight\hbox{$
\atomicflow{
( 4,14)*{\afvju4{\tilde\epsilon_1}{}};
( 6,14)*{\cdots};
( 8,14)*{\afvju4{}{\tilde\epsilon_h}};
( 0, 4)*{\afvjd4{\tilde\three}{}};
( 9,10)*{\aflabelleft{\tilde A}};
( 4, 9)*{\affr{10}6};
( 4, 4)*{\afvjd4{\tilde\epsilon'_1}{}};
( 6, 4)*{\cdots};
( 8, 4)*{\afvjd4{}{\tilde\epsilon'_k}};
( 0,16)*{\afawd{}{}{\tilde\two}{}};
(-2, 4)*{\invisiblemark};
(11, 4)*{\invisiblemark}}$}
\qquad\hbox{and}\qquad
B''=\lower2\atflowelheight\hbox{$
\atomicflow{
( 0,14)*{\afvju4{\hat\epsilon_1}{}};
( 2,14)*{\cdots};
( 4,14)*{\afvju4{}{\hat\epsilon_h}};
( 8,14)*{\afvju4{}{\hat\two}};
( 9,10)*{\aflabelleft{\hat A}};
( 4, 9)*{\affr{10}6};
( 0, 4)*{\afvjd4{\hat\epsilon'_1}{}};
( 2, 4)*{\cdots};
( 4, 4)*{\afvjd4{}{\hat\epsilon'_k}};
( 8, 2)*{\afawu{}{}{}{\hat\three}};
(-3, 4)*{\invisiblemark};
(10, 4)*{\invisiblemark}}$}
\quad,
\]
\afnegspace
where the correspondence of edges has been indicated by adding accents to their labels. We have that:
\begin{itemize}
\item if $\one$ is an edge belonging to an $\ai$-cycle, $B'\to_\frbc D'$ and $B''\to_\frbc D''$ then $B\to_\frbc C$;
\item if $\one$ is an extremal simple edge, $B'\to_\frex D'$ and $B''\to_\frex D''$ then $B\to_\frex C$.
\end{itemize}
\end{defi}

%---------------------------------------
\begin{exa}
Consider the atomic flow to the left in Figure~\ref{FigExRed}. Assuming that the two evidenced simple edges both belong to $\ai$-cycles and that the box $A$ stands for a cycle-free flow, then the atomic flow on the right is the result of a $\to_\frbc$ reduction. Similarly, if the two evidenced simple edges are extremal simple edges, and the box stands for a flow that contains no simple edges, then the atomic flow on the right is the result of a $\to_\frex$ reduction.
\end{exa}

%---------------------------------------
\begin{figure}[tbp]
\[
%-------------------
\atomicflow{
(14,24)*{\afaidex{}{}{}{}{}{}31};
( 0,20)*{\afvju8{}{}};
( 4,20)*{\afvju8{}{}};
( 8,18)*{\afvju4{\one}{}};
(14,20)*{\afaidnw{}{}};
(12,18)*{\afvju4{\three}{}};
( 2,18)*{\cdots};
(16,16)*{\afvj8};
( 6,14)*{\affr{14}4};
(13,15)*{\aflabelleft A};
(20,14)*{\afvj{12}};
( 2,10)*{\cdots};
( 8,10)*{\afvjd4{\two}{}};
(14, 8)*{\afaiu{\four}{}{}{}{}{}};
( 0, 8)*{\afvjd8{}{}};
( 4, 8)*{\afvjd8{}{}};
(14, 4)*{\afaiuex{}{}{}{}{}{}31}}
\quad\to_\frbc\quad
%-------------------
\atomicflow{
( 8,30  )*{\afacunwexsq{}{}{}{}41};
(12,30  )*{\afacunwexsq{}{}{}{}41};
(16,30  )*{\afacunwexsq{}{}{}{}41};
(12,28  )*{\afawd{}{}{\three}{}};
(16,26  )*{\afvju2{}{}};
(20,26  )*{\afvju2{}{}};
(24,26  )*{\afvju2{\one}{}};
(18,26  )*{\cdots};
( 0,21  )*{\afvju{10}{}{}};
( 4,21  )*{\afvju{10}{}{}};
( 8,21  )*{\afvju{10}{\one}{}};
(18,22  )*{\affr{14}4};
(25,23  )*{\aflabelleft A};
(28,21  )*{\afvj{38}};
(32,21  )*{\afvj{38}};
( 2,18  )*{\cdots};
(12,18  )*{\afvju4{\three}{\four}};
(18,18  )*{\cdots};
(24,16  )*{\afawu{}{}{\two}{}};
( 6,14  )*{\affr{14}4};
(13,15  )*{\aflabelleft A};
(16,15  )*{\afvjd{10}{}{}};
(20,15  )*{\afvjd{10}{}{}};
( 2,10  )*{\cdots};
( 0,10.4)*{\afvjd2{}{}};
( 4,10.4)*{\afvjd2{}{}};
( 8, 8  )*{\afawu{}{}{\two}{}};
(12, 8  )*{\afawu{}{}{\four}{}};
%---------
(20,82)*{\cdots};
(18,78)*{\afacuexsq{}{}{}{}{}{}51};
(22,78)*{\afacuexsq{}{}{}{}{}{}51};
(20,37)*{\afcjr85};
(18,-4)*{\afacdexsq{}{}{}{}{}{}51};
(22,-4)*{\afacdexsq{}{}{}{}{}{}51};
(20,-8)*{\cdots};
%---------
( 8,4)*{\afacdnwexsq{}{}{}{}41};
(12,4)*{\afacdnwexsq{}{}{}{}41};
%---------
(16,38)="A";
"A"+(12,30  )*{\afacunwexsq{}{}{}{}41};
"A"+(16,30  )*{\afacunwexsq{}{}{}{}41};
"A"+(12,28  )*{\afawd{}{}{\three}{}};
"A"+(16,28  )*{\afawd{}{}{}{\one}};
"A"+(20,26  )*{\afvju2{}{}};
"A"+(24,26  )*{\afvju2{}{}};
"A"+(22,26  )*{\cdots};
"A"+( 4,21  )*{\afvju{10}{}{}};
"A"+( 8,21  )*{\afvju{10}{}{}};
"A"+(18,22  )*{\affr{14}4};
"A"+(25,23  )*{\aflabelleft A};
"A"+( 0,20  )*{\afawd{}{}{\one}{}};
"A"+( 6,18  )*{\cdots};
"A"+(12,18  )*{\afvju4{\three}{\four}};
"A"+(22,18  )*{\cdots};
"A"+( 6,14  )*{\affr{14}4};
"A"+(13,15  )*{\aflabelleft A};
"A"+(16,15  )*{\afvjd{10}{}{\two}};
"A"+(20,15  )*{\afvjd{10}{}{}};
"A"+(24,15  )*{\afvjd{10}{}{}};
"A"+( 6,10  )*{\cdots};
"A"+( 0,10.4)*{\afvjd2{\two}{}};
"A"+( 4,10.4)*{\afvjd2{}{}};
"A"+( 8,10.4)*{\afvjd2{}{}}; 
"A"+(12, 8  )*{\afawu{}{}{\four}{}};
"A"+( 8, 4  )*{\afacdnwexsq{}{}{}{}41};
"A"+(12, 4  )*{\afacdnwexsq{}{}{}{}41};
"A"+(16, 4  )*{\afacdnwexsq{}{}{}{}41};
%---------
( 8,53)*{\afvj{38}};
(12,53)*{\afvj{38}}}
\]
\caption{Example of a two-step $\to_\frbc$ (or $\to_\frex$) reduction.}
\label{FigExRed}
\end{figure}

Notice that the flow in Figure~\ref{FigExRed} represents the `external' shape of any flow after eliminating any two simple edges. Eliminating more simple edges would follow the same pattern.

%---------------------------------------
\begin{rem}\label{RemGenRed}
It is possible to generalise the construction in Figure~\ref{FigExRed} to any number $n$ of simple edges: for any $n$, there is an atomic flow of the same nature as the one at the right of the figure, with $2^n$ boxes. So, a $\to_\frbc$ or $\to_\frex$ reduction can be `executed' in one step if the simple edges involved, and their order, is known in advance.
\end{rem}

The following two theorems guarantee properties of flows after reducing by $\to_\frbc$ and $\to_\frex$, provided that the flow to which we apply them meets some conditions. These conditions can be achieved by a careful use of the flow rewriting system $\frc$ for contractions, as we see later on.

%---------------------------------------
\begin{thm}\label{ThALSE}
If all the\/ $\ai$-cycles in atomic flow $B$ are fragile cycles then there exists a cycle-free flow $C$ such that $B\to_\frbc C$.
\end{thm}

%---------------------------------------
\proof
By induction on the number of $\ai$-cycles: follow Definition~\ref{DefCBFlow} and note that when composing $D'$ and $D''$ in $C$, no $\ai$-cycles are created.
\qed

%---------------------------------------
\begin{exa}
Notice that removing one simple edge might break more than one $\ai$-cycle. The following flows have, respectively, two, three and two fragile cycles. Removing edges $\one$, $\two$, $\three$ or $\four$ breaks two $\ai$-cycles each:
\nopagebreak[4]\medskip\afnegspace
\[
\atomicflow{
(2,  8)*{\afaidnw{}{}{}{}};
(0, -4)*{\afacd{}{}{}{}{}{}};
(0,  4)*{\afacu{}{}{}{}{}{}};
(2,-10)*{\afaiunw{}{}{}{}};
(4,  0)*{\afvju{16}\one{}}}
\quad,\qquad
\atomicflow{
(-6, 0)*{\afvju{12}\two{}};
(-4, 6)*{\afaidnw{}{}{}{}};
(-4,-8)*{\afaiunw{}{}{}{}};
( 0, 2)*{\afacd{}{}{}{}{}{}};
( 0,-4)*{\afacunw{}{}{}{}};
( 4, 6)*{\afaidnw{}{}{}{}};
( 4,-8)*{\afaiunw{}{}{}{}};
( 6, 0)*{\afvju{12}{}\three}}
\qquad\hbox{and}\qquad
\atomicflow{
(-10, 0)*{\afvj8{}{}{}{}};
( -8, 4)*{\afaidnw{}{}{}{}};
( -4, 0)*{\afacd{}{}{}{}{}{}};
( -2,-6)*{\afaiunw{}{}{}{}};
( -4,-8)*{\afaiuex{}{}{}{}{}{}31};
(  0, 0)*{\afvju8\four{}};
( 10, 0)*{\afvj8{}{}{}{}};
(  8,-6)*{\afaiunw{}{}{}{}};
(  4,-0)*{\afacu{}{}{}{}{}{}};
(  2, 4)*{\afaidnw{}{}{}{}};
(  4, 8)*{\afaidex{}{}{}{}{}{}31}}
\quad.
\]
\afnegspace
\end{exa}

%---------------------------------------
\begin{rem}
Reductions $\to_\frse$ and $\to_\frbc$ can introduce new simple edges along the edge labelled $\hat\two$ and $\tilde\three$ in Definitions~\ref{DefRedS} and \ref{DefCBFlow}. This phenomenon is the reason for having introduced the notion of extremal simple edge. As an example, consider the following application of $\to_\frse$ to a non-extremal simple edge. The problem occurs when the two copies of the given atomic flow (neither of which contains a simple edge) are combined, and the edge
that connects them becomes a simple edge:
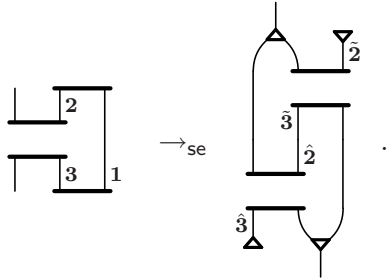
\[
\atomicflow{
(0  , 2)*{\afaiu{}{}{}{\two}{}{}};
(4  , 6)*{\afaidnw{}{}{}{}};
(6  , 0)*{\afvju{12}{}{\one}{}{}};
(4  ,-8)*{\afaiunw{}{}{}{}};
(0  ,-2)*{\afaid{}{}{}{\three}{}{}};
(7.5, 0)*{\invisiblemark}}
\quad\to_\frse\quad
\atomicflow{
(-2, 12)*{\afacu{}{}{}{}{}{}};
( 4, 12)*{\afawd{}{}{}{\tilde\two}};
( 2,  6)*{\afaiunw{}{}{}{}};
( 2,  4)*{\afaid{\tilde\three}{}{}{}{}{}};
( 4, -4)*{\afvju{8}{}{}{}{}};
(-2, -4)*{\afaiu{}{}{}{\hat\two}{}{}};
(-2, -8)*{\afaidnw{}{}{}{}};
(-4,  4)*{\afvju{8}{}{}{}{}};
(-4,-12)*{\afawu{}{}{\hat\three}{}};
( 2,-12)*{\afacd{}{}{}{}{}{}}}
\quad.
\]
Reducing over extremal simple edges avoids the problem.
\end{rem}

The algorithms we show in this paper would terminate even if we did not insist on reducing over extremal simple edges, but doing so simplifies the induction measure.

%---------------------------------------
\begin{thm}\label{ThALSP}
If all the\/ $\ai$-paths in cycle-free atomic flow $B$ are clean paths then there exists a flow $C$ such that no\/ $\ai$-connections appear in it and $B\to_\frex C$.
\end{thm}

%---------------------------------------
\proof
By induction on the number of simple edges. We follow Definition~\ref{DefCBFlow} and its notation. We have to verify that, when composing $D'$ and $D''$ in $C$, no $\ai$-connections are created. This could only happen if, in $A$, edge $\two$ were upper edge of a cointeraction and edge $\three$ were lower edge of an interaction; this is impossible, because $\one$ is an extremal simple edge.
\qed

%---------------------------------------
\begin{thm}\label{ThLBSound}
Reductions\/ $\to_\frbc$ and\/ $\to_\frex$ are sound.
\end{thm}

%---------------------------------------
\proof
The proof is almost identical for $\to_\frbc$ and $\to_\frex$, so, to fix ideas, we prove the theorem for $\to_\frbc$, by induction on its definition. The base case is trivial. For the inductive case, the construction is similar to the one in the proof of Theorem~\ref{ThSESound}, and we see it briefly. Consider
\[
\Phi=\!
\vlderivation                                           {
\vlde{\Phi_3}{}{\zeta\{\mk \fff\}                 }    {
\vlin{\aiu  }{}{\zeta'\{\mk \fff\}                }   {
\vlde{\Phi_2}{}{\zeta'\vlscn(\bar a^\three.a^\one)}  {
\vlin{\aid  }{}{\xi'\vlscn[\bar a^\two.a^\one]    } {
\vlde{\Phi_1}{}{\xi'\{\mk \ttt\}                  }{
\vlhy          {\xi\{\mk \ttt\}                   }}}}}}}
\quad,\qquad
\Psi'=\;\;
\vlderivation                                                  {
\vlde{\Phi_3\{\mk\fff\ot\bar a\}}
             {}{\zeta\{\bar a^{\tilde\three}\}          }     {
\vlin{=     }{}{\zeta'\{\bar a^{\tilde\three}\}         }    {
\vlde{\Phi_2\{a^\one\ot\ttt\}}
             {}{\zeta'\vlscn(\bar a^{\tilde\three}.\ttt)}   {
\vlin{\awd  }{}{\xi'\vlscn[\bar a^{\tilde\two}.\ttt]    }  {
\vlin{=     }{}{\xi'\vlscn[\fff.\ttt]                   } {
\vlde{\Phi_1}{}{\xi'\{\ttt\}                            }{
\vlhy          {\xi\{\ttt\}                             }}}}}}}}
\qquad\hbox{and}\qquad
\Psi''=\;\;
\vlderivation                                              {
\vlde{\Phi_3}
           {}{\zeta\{\fff\}                         }     {
\vlin{=   }{}{\zeta'\{\fff\}                        }    {
\vlin{\awu}{}{\zeta'\vlscn(\ttt.\fff)               }   {
\vlde{\Phi_2\{a^\one\ot \fff\}}
           {}{\zeta'\vlscn(\bar a^{\hat\three}.\fff)}  {
\vlin{=   }{}{\xi'\vlscn[\bar a^{\hat\two}.\fff]    } {
\vlde{\Phi_1\{\mk\ttt\ot\bar a\}}
           {}{\xi'\{\bar a^{\hat\two}\}             }{
\vlhy        {\xi\{\bar a^{\hat\two}\}              }}}}}}}}
\quad,
\]
where $\Phi$ has flow $B$, such that $B\to_\frbc C$, and $\Psi'$ and $\Psi''$ are obtained from $\Phi$ and have flows $B'$ and $B''$, respectively, as per Definition~\ref{DefCBFlow}. By induction hypothesis, there are derivations $\hat\Psi'$ and $\hat\Psi''$ with flows $D'$ and $D''$ such that $B'\to_\frbc D'$ and $B''\to_\frbc D''$ and with the same premisses and conclusions as $\Psi$ and $\Psi'$. We can compose these derivations into the following, whose flow is $C$:
\[
\vlderivation                                        {
\vlin{\cod}{}{\zeta\{\fff\}                    }    {
\vlde{\vls[\hat\Psi''.\zeta\{\fff\}]}
           {}{\vls[\zeta\{\fff\}.\zeta\{\fff\}]}   {
\vlin{\sus}{}{\vls[\xi\{\bar a\}.\zeta\{\fff\}]}  {
\vlde{\vls(\xi\{\ttt\}.\hat\Psi')}
           {}{\vls(\xi\{\ttt\}.\zeta\{\bar a\})} {
\vlin{\cou}{}{\vls(\xi\{\ttt\}.\xi\{\ttt\})    }{
\vlhy        {\xi\{\ttt\}                      }}}}}}}
\quad.
\]
\qed

%---------------------------------------
\begin{rem}
Similarly to what we do in Remark~\ref{RemIndep}, we observe here that the previous soundness theorem holds for any proof system containing the same structural rules as $\SKS$ and whose logical rules are such that the $\cou$, $\sus$ and $\cod$ rules are derivable. So, the soundness theorem depends only very loosely on the choice of logical rules.
\end{rem}

%---------------------------------------
\begin{rem}\label{RemExpGr}
After a $\to_\frbc$ or $\to_\frex$ reduction, the size of an atomic flow grows by a factor of $O(2^n)$, where $n$ is the number of simple edges involved in the reduction (see Remark~\ref{RemGenRed}). This accounts for an equally exponential growth in the corresponding derivation.
\end{rem}

\newcommand{\BC}{{\mathsf{BC}}}
\newcommand{\EX}{{\mathsf{EX}}}
Out of the two reductions $\to_\frbc$ and $\to_\frex$, we define two algorithms, $\BC$ and $\EX$, which use flow rewriting system $\frc$ to achieve the necessary preconditions for the two reductions.

%---------------------------------------
\begin{defi}\label{DefCB}
For every $\SKS$ derivation $\Phi_0$ with atomic flow $A_0$, \emph{algorithm\/ $\BC$} (for \emph{break\/ $\ai$-cycles}) performs the following steps: 
\begin{enumerate}
\item\emph{Make\/ $\ai$-cycles fragile.} Transform $\Phi_0$ into $\Phi_1$, whose flow $A_1$ is obtained by assigning a polarity to $A_0$ and applying $\frc$ over negative contraction and cocontraction vertices belonging to $\ai$-cycles, until they are all positive.
\item\emph{Break\/ $\ai$-cycles.} Transform $\Phi_1$ into the algorithm's output $\Phi_2$, whose flow $A_2$ is such that $A_1\to_\frbc A_2$.
\end{enumerate}
\end{defi}

%---------------------------------------
\begin{thm}\label{TheoElCycles}
Given any\/ $\SKS$ derivation, algorithm\/ $\BC$ produces a derivation with the same premiss and conclusion and such that its atomic flow is cycle-free.
\end{thm}

%---------------------------------------
\proof
We refer to Definition~\ref{DefCB}. Step 1 transforms the given derivation $\Phi_0$ into $\Phi_1$, which, by Theorems~\ref{TheoCycleFragile} and \ref{TheoSound}, exists and all its $\ai$-cycles are fragile cycles and its premiss and conclusion are the same as those of $\Phi_0$. Then, by Theorems~\ref{ThALSE} and \ref{ThLBSound} we conclude that $\Phi_2$, obtained from $\Phi_1$ in Step~2, exists, is cycle-free and its premiss and conclusion are the same as those of $\Phi_1$.
\qed

The first and second reduction steps in Figure~\ref{FigExStrRed} provide an example of application of algorithm $\BC$.

%---------------------------------------
\begin{defi}\label{DefCP}
For every $\SKS$ derivation $\Phi_2$, with cycle-free atomic flow $A_2$, \emph{algorithm\/ $\EX$} (for \emph {eliminate extremal simple edges}) performs the following steps: 
\begin{enumerate}
\item\emph{Make $\ai$-paths clean.} Transform $\Phi_2$ into $\Phi_3$, whose flow $A_3$ is normal for $\frc$.
\item\emph{Remove simple edges.} Transform $\Phi_3$ into the algorithm's output $\Phi_4$, whose flow $A_4$ is such that $A_3\to_\frex A_4$.
\end{enumerate}
\end{defi}

%---------------------------------------
\begin{thm}\label{TheoElPaths}
Given any\/ $\SKS$ derivation whose atomic flow is cycle-free, algorithm\/ $\EX$ produces a derivation with the same premiss and conclusion and such that, in its atomic flow, there are no\/ $\ai$-connections.
\end{thm}

%---------------------------------------
\proof
We refer to Definition~\ref{DefCP}. Step 1 transforms the given derivation $\Phi_2$ into $\Phi_3$, which, by Theorems~\ref{TheoCTerm} and \ref{TheoSound} and by Propositions~\ref{PropCLF} and \ref{PropCOW}, exists, is cycle-free and all its $\ai$-paths are clean paths and its premiss and conclusion are the same as those of $\Phi_2$. Then, by Theorems~\ref{ThALSP} and \ref{ThLBSound} we conclude that $\Phi_4$, obtained from $\Phi_3$ in Step~2, exists, there are no $\ai$-connections in its atomic flow, and its premiss and conclusion are the same as for $\Phi_3$.
\qed

The output of $\BC$ and $\EX$ is not normal for $\frc$, in general, as a quick look at Figure~\ref{FigExRed} shows. 

%-------------------------------------------------------------------------------
\subsection{Streamlining Algorithms}\label{SubsStrPr}

We obtain here our main result: normalisation algorithms for propositional logic that normalise generic derivations, and that entail cut elimination on proofs. We show two (similar) such algorithms, but it is clear that there are many possible variations. These algorithms are obtained as combinations of the reductions $\to_\frw$, $\to_\frc$, $\to_\frbc$ and $\to_\frex$. We strengthen the notion of streamlined derivation introduced in the beginning, in order to appreciate the variations in the design of the normalising algorithms.

\newcommand{\Str}{\mathsf{Str}}
%---------------------------------------
\begin{defi}\label{DefStr}
For every $\SKS$ derivation $\Phi_0$, \emph{algorithm\/ $\Str$} (for \emph{streamlining}) performs the following steps:
\begin{enumerate}
\item\emph{Break\/ $\ai$-cycles.} Apply $\BC$ to $\Phi_0$ and obtain $\Phi_2$.
\item\emph{Eliminate\/ $\ai$-connections.} Apply $\EX$ to $\Phi_2$ and obtain $\Phi_4$.
\item\emph{Move away weakenings and coweakenings.} Transform $\Phi_4$ into algorithm's output $\Phi_5$, whose flow $A_5$ is normal for $\frw$.
\end{enumerate}
\end{defi}

%---------------------------------------
\begin{exa}
In Figure~\ref{FigExStrRed}, algorithm $\Str$ is applied. For each flow the corresponding derivation is shown. Inference rules are used modulo the relation $=$ for formulae. Trivial rule instances are omitted. The atomic flow we start out with has an $\ai$-cycle, but no simple edge. The first step is to move a contraction to create a simple edge. The second step removes the simple edge. The remaining steps are weakening reductions. The reader can apply the construction in the proof of Theorem~\ref{ThSESound} (or \ref{ThLBSound}) in order to generate the third derivation. We recall that adding some switch instances can be necessary in order to `pull' the units involved in interactions and cointeractions to the premiss and conclusion of the derivation to be reduced.
\end{exa}

\newcommand{\LD}[1]{{\color{Lavender}#1}}
\newcommand{\AM}[1]{{\color{Aquamarine}#1}}
\newcommand{\RR}[1]{{\color{RubineRed}#1}}
\newcommand{\LG}[1]{{\color{LimeGreen}#1}}
\newcommand{\OG}[1]{{\color{OliveGreen}#1}}
%---------------------------------------
\begin{figure}[tbp]
\[
\begin{array}{@{}c@{}c@{}}
\atomicflow{
( 5,12)*{\afvjcol4{Red}};
( 3, 6)*{\afacdcol{}{}{}{}{}{}{Green}{Red}{Lavender}};
(-3, 6)*{\afacucol{}{}{}{}{}{}{DarkOrchid}{YellowOrange}{OliveGreen}};
( 1, 0)*{\afaiunw{}{}};
(-5, 0)*{\afvjcol4{DarkOrchid}};
(-1,10)*{\afaidnw{}{}}}\\
{}_\frc\raise3pt\hbox{$\downarrow$}\\
\atomicflow{
( -2,10)*{\afacucol{}{}{}{}{}{}{DarkOrchid}{YellowOrange}{OliveGreen}};
( -6, 8)*{\afvjcol{12}{Green}};
(  6,10)*{\afvjcol{16}{Red}};
( -6, 4)*{\afcjrcol44{DarkOrchid}};
(  0, 4)*{\afacunwcol{}{}{}{}{LimeGreen}{ProcessBlue}};
(  4, 0)*{\afaiunw{}{}};
( -4, 0)*{\afaiunw{}{}};
( -8, 0)*{\afvjcol4{DarkOrchid}};
( -4,14)*{\afaidnw{}{}}}
&\quad\to_\frbc\quad
\end{array}
\atomicflow{
( 7,20)*{\afvjcol{28}{RawSienna}};
(-3, 8)*{\afvjcol{4}{RubineRed}};
(-2,20)="A";
"A"+( 5, 8)*{\afvjcol{12}{Green}};
"A"+(-3,10)*{\afacucol{}{}{}{}{}{}{Periwinkle}{YellowOrange}{OliveGreen}};
"A"+(-1, 4)*{\afacunwcol{}{}{}{}{LimeGreen}{ProcessBlue}};
"A"+( 3, 0)*{\afaiunw{}{}};
"A"+(-5,-4)*{\afvjcol{20}{Periwinkle}}; 
"A"+(-3,14)*{\afawdnw{}{}};
(-3,18)*{\afcjlcol48{LimeGreen}};
(-1,12)*{\afacunwcol{}{}{}{}{RubineRed}{Aquamarine}};
( 1, 8)*{\afacunwcol{}{}{}{}{Lavender}{SpringGreen}};
( 5, 4)*{\afaiunw{}{}};
(-1, 4)*{\afawunw{}{}};
( 5,38)*{\afacucol{}{}{}{}{}{}{Green}{RawSienna}{Red}};
(-5, 2)*{\afacdcol{}{}{}{}{}{}{Periwinkle}{RubineRed}{DarkOrchid}}}
\quad\to_\frw^\star
\atomicflow{
( 7,20)*{\afvjcol{28}{RawSienna}};
( 3,10)*{\afawdcol{}{}{}{}{}{SpringGreen}};
(-1,10)*{\afawdcol{}{}{}{}{}{Lavender}};
( 5, 4)*{\afaiunw{}{}};
(-1, 4)*{\afawunw{}{}};
( 1,20)="A";
"A"+( 2,8)*{\afvjcol{12}{Green}};
"A"+(-2,6)*{\afawdcol{}{}{}{}{}{ProcessBlue}};
"A"+( 0,0)*{\afaiunw{}{}};
( 5,38)*{\afacucol{}{}{}{}{}{}{Green}{RawSienna}{Red}};
(-5, 2)*{\afacdcol{}{}{}{}{}{}{Periwinkle}{RubineRed}{DarkOrchid}};
(-7, 6)*{\afawdnw{}{}};
(-3, 6)*{\afawdnw{}{}}}
\quad\to_\frw^\star\quad
\atomicflow{
(2,16)*{\afacucol{}{}{}{}{}{}{Green}{RawSienna}{Red}};
(0,10)*{\afawunw{}{}};
(4,10)*{\afawunw{}{}};
(2, 6)*{\afawdcol{}{}{}{}{}{DarkOrchid}}}
\quad\to_\frw\quad
\atomicflow{
(0,6)*{\afawucol{}{}{}{}{}{Red}};
(0,0)*{\afawdcol{}{}{}{}{}{DarkOrchid}}}
\]
\bigskip
%---------------------------------------
\[
%---------------------------------------
\newbox\derone\setbox\derone=\hbox{$
\vlderivation                                                      {
\vlin{\aiu}{}{\vls [(\DO{a}.  \fff              ).\ttt   ] }      {
\vlin{\swi}{}{\vls [(\DO{a}.  \YO{a}.\LD{\bar a}).\ttt   ] }     {
\vlin{\acu}{}{\vls([(\DO{a}.\YO{a}).\ttt    ].\LD{\bar a}) }    {
\vlin{\acd}{}{\vls([\OG{a}             .\ttt].\LD{\bar a}) }   {
\vlin{\med}{}{\vls([\OG{a}.\ttt].[\GR{\bar a}.\RD{\bar a}])}  {
\vlin{\swi}{}{\vls [(\OG{a}.\RD{\bar a}     ).\GR{\bar a}] } {
\vlin{\aid}{}{\vls([ \OG{a}.\GR{\bar a}].\RD{\bar a}     ) }{
\vlhy        {                           \RD{\bar a}       }}}}}}}}}
$}
%---------------------------------------
\newbox\dertwo\setbox\dertwo=\hbox{$
\vlderivation                                                                 {
\vlin{\aiu}{}{\vls[(\DO{a}.\fff              ).\ttt]                }        {
\vlin{\aiu}{}{\vls[(\DO{a}.\PB{a}.\RD{\bar a}).\ttt]                }       {
\vlin{\swi}{}{\vls[(\DO{a}.\PB{a}.[(\GR{\bar a}.\LG{a}).\RD{\bar a}]).\ttt]
                                                                    }      {
\vlin{\acu}{}{\vls[(\DO{a}.\LG{a}.\PB{a}.[\GR{\bar a}.\RD{\bar a}]).\ttt]
                                                                    }     {
\vlin{\swi}{}{\vls [(\DO{a}.\YO{a}.[\GR{\bar a}.\RD{\bar a}]).\ttt] }    {
\vlin{\acu}{}{\vls([(\DO{a}.\YO{a}).\ttt].[\GR{\bar a}.\RD{\bar a}])}   {
\vlin{\med}{}{\vls([ \OG{a}        .\ttt].[\GR{\bar a}.\RD{\bar a}])}  {
\vlin{\swi}{}{\vls  [\GR{\bar a}.(\OG{a} .\RD{\bar a})]             } {
\vlin{\aid}{}{\vls ([\GR{\bar a}. \OG{a}].\RD{\bar a})              }{
\vlhy        {                            \RD{\bar a}               }}}}}}}}}}}
$}
%---------------------------------------
\newbox\derthree\setbox\derthree=\hbox{$
\vlderivation                                                          {
\vlin{\cod}{}{\vls[(\DO{a}.\fff).\ttt]             }                  {
\vlin{\aiu}{}{\vls[(\PW{a}.\fff).\ttt.(\RR{a}.\fff)]
                                                   }                 {
\vlin{\awu}{}{\vls[(\PW{a}.\fff).\ttt.(\RR{a}.\SG{a}.\RS{\bar a})]
                                                   }                {
\vlin{\swi}{}{\vls[(\PW{a}.\fff).\ttt.
       (\RR{a}.\SG{a}.[(\LD{a}.\fff).\RS{\bar a}])]}               {
\vlin{\acu}{}{\vls[(\PW{a}.\fff).\ttt.(\RR{a}.\LD{a}.\SG{a}.\RS{\bar a})]
                                                   }              {
\vlin{\swi}{}{\vls[(\PW{a}.\fff).\ttt.(\RR{a}.\AM{a}.\RS{\bar a})]
                                                   }             {
\vlin{\acu}{}{\vls[(\PW{a}.\fff).\ttt.([(\RR{a}.\AM{a}).\ttt].\RS{\bar a})]
                                                   }            {
\vlin{\med}{}{\vls[(\PW{a}.\fff).\ttt.([\LG{a}.\ttt].\RS{\bar a})]
                                                   }           {
\vlin{\sus}{}{\vls[(\PW{a}.\fff).\ttt.(\LG{a}.\RS{\bar a})]
                                                   }          {
\vlin{\aiu}{}{\vls([(\PW{a}.\LG{a}).\ttt].\RS{\bar a})
                                                   }         {
\vlin{\swi}{}{\vls([(\PW{a}.[\LG{a}.(\PB{a}.\GR{\bar a})]).\ttt].\RS{\bar a})
                                                   }        {
\vlin{\swi}{}{\vls([(\PW{a}.\PB{a}.[\LG{a}.\GR{\bar a}]).\ttt].\RS{\bar a})
                                                   }       {
\vlin{\acu}{}{\vls([(\PW{a}.\LG{a}.\PB{a}.[\ttt.\GR{\bar a}]).\ttt].\RS{\bar a})
                                                   }      {
\vlin{\swi}{}{\vls([(\PW{a}.\YO{a}.[\ttt.\GR{\bar a}]).\ttt].\RS{\bar a})
                                                   }     {
\vlin{\acu}{}{\vls([(\PW{a}.\YO{a}).\ttt].[\ttt.\GR{\bar a}].\RS{\bar a})
                                                   }    {
\vlin{\med}{}{\vls([\OG{a}.\ttt].[\ttt.\GR{\bar a}].\RS{\bar a})
                                                   }   {
\vlin{\swi}{}{\vls([\ttt.(\OG{a}.\GR{\bar a})].\RS{\bar a})
                                                   }  {
\vlin{\awd}{}{\vls([\ttt.\OG{a}].\GR{\bar a}.\RS{\bar a})
                                                   } {
\vlin{\cou}{}{\vls(\GR{\bar a}.\RS{\bar a})        }{
\vlhy        {                 \RD{\bar a}         }}}}}}}}}}}}}}}}}}}}}
$}
%---------------------------------------
\newbox\derfour\setbox\derfour=\hbox{$
\vlderivation                                                          {
\vlin{\cod}{}{\vls[(\DO{a}.\fff).\ttt]              }                 {
\vlin{\awd}{}{\vls[(\PW{a}.\fff).\ttt.(\RR{a}.\fff)]}                {
\vlin{\awd}{}{\vls[(\PW{a}.\fff).\ttt]              }               {
\vlin{\aiu}{}{\vls{\ttt}                            }              {
\vlin{\awd}{}{\vls[\ttt.(\fff.\SG{a}.\RS{\bar a})]  }             {
\vlin{\awu}{}{\vls[\ttt.(\fff.\RS{\bar a})]         }            {
\vlin{\awd}{}{\vls[\ttt.(\fff.[(\LD{a}.\fff).\RS{\bar a}])]
                                                    }           {
\vlin{\swi}{}{\vls[\ttt.(\fff.\RS{\bar a})]         }          {
\vlin{\med}{}{\vls[\ttt.\RS{\bar a}]                }         {
\vlin{\sus}{}{\vls[\ttt.(\fff.\RS{\bar a})]         }        {
\vlin{\aiu}{}{\vls{\RS{\bar a}}                     }       {
\vlin{\awd}{}{\vls([(\fff.\PB{a}.\GR{\bar a}).\ttt].\RS{\bar a})
                                                    }      {
\vlin{\swi}{}{\vls([(\fff.\GR{\bar a}).\ttt].\RS{\bar a})
                                                    }     {
\vlin{\swi}{}{\vls([(\fff.[\ttt.\GR{\bar a}]).\ttt].\RS{\bar a})
                                                    }   {
\vlin{\med}{}{\vls([\ttt.\GR{\bar a}].\RS{\bar a})  }  {
\vlin{\swi}{}{\vls([\ttt.(\fff.\GR{\bar a})].\RS{\bar a})
                                                    } {
\vlin{\cou}{}{\vls(\GR{\bar a}.\RS{\bar a})         }{
\vlhy        {                 \RD{\bar a}          }}}}}}}}}}}}}}}}}}}
$}
%---------------------------------------
\newbox\derfive\setbox\derfive=\hbox{$
\vlderivation                                   {
\vlin{\awd}{}{\vls[(\DO{a}.\fff).\ttt]     }   {
\vlin{\awu}{}{\vls{\ttt}                   }  {
\vlin{\awu}{}{\vls{\RS{\bar a}}            } {
\vlin{\cou}{}{\vls(\GR{\bar a}.\RS{\bar a})}{
\vlhy        {                 \RD{\bar a} }}}}}}
$}
%---------------------------------------
\newbox\dersix\setbox\dersix=\hbox{$
\vlderivation                            {
\vlin{\awd}{}{\vls[(\DO{a}.\fff).\ttt]} {
\vlin{\awu}{}{\ttt                    }{
\vlhy        {\RD{\bar a}             }}}}
$}
%---------------------------------------
\xy
(  0  ,110  )*{\box\derone};
(  1.5, 82  )*{{}_\frc\raise3pt\hbox{$\downarrow$}};
(  0  , 50  )*{\box\dertwo};
( 20  , 70  )*{\to_\frbc};
( 50  , 43.5)*{\box\derthree};
( 83  , 70  )*{\to_\frw^\star};
(110  , 81  )*{\box\derfour};
( 97.5, 23  )*{{}_\frw\raise3pt\hbox{$\downarrow$}_\star};
( 91  ,  0  )*{\box\derfive};
(106  ,  5  )*{\to_\frw};
(120  ,  0  )*{\box\dersix}
\endxy
\]
\caption{Example of application of algorithm $\Str$.}
\label{FigExStrRed}
\end{figure}

We are about to prove that algorithm $\Str$ does more than just streamlining derivations: it turns out that Step~3 performs unnecessary reductions on weakening and coweakening vertices connected with contraction and cocontraction ones. However, the shape of the derivations produced by $\Str$ is notable, so we prefer to present $\Str$ this way, and it is understood that its definition can be weakened if strict streamlining is all that is required: we only need to apply the weakening reductions that are necessary to get streamlined derivations, and no more.

We now make precise the notion of streamlining obtained by $\Str$ (called `super-stream\-lin\-ing'), and a further strengthening that we explore in the rest of this section.

%---------------------------------------
\begin{defi}
An $\SKS$ derivation is \emph{super-streamlined} if it is streamlined and its associated atomic flow is normal for $\frw$. An $\SKS$ derivation is \emph{hyper-streamlined} if it is super-streamlined and its associated atomic flow is normal for $\frc$. We say that an algorithm $P$ \emph{super-streamlines} (resp., \emph{hyper-streamlines}) a derivation $\Phi$ if the output of $P$ on $\Phi$ is a super-streamlined (resp., hyper-streamlined) derivation that has the same premiss and conclusion as $\Phi$.
\end{defi}

%---------------------------------------
\begin{exa}
In Example~\ref{ExStream} (page~\pageref{ExStream}), the second flow is streamlined, but not super-stream\-lined, and the third is hyper-streamlined (so, super-streamlined as well). 
\end{exa}

%---------------------------------------
\begin{rem}
Consider the following figure
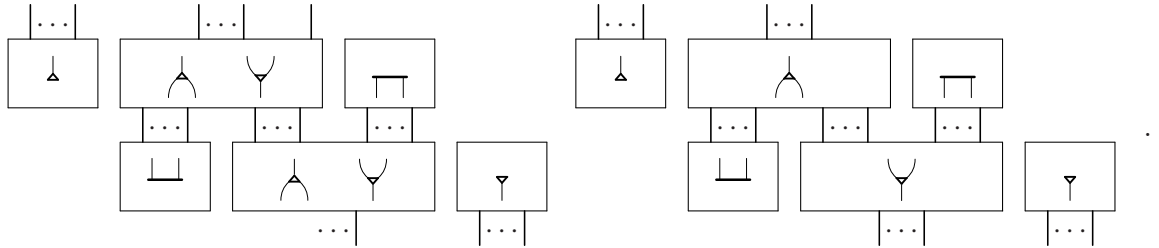
\newbox\boxone\setbox\boxone=\hbox{$
   \divide\atflowunit by5\multiply\atflowunit by3\afsetunits
   \atomicflow{(0,0)*{\afacd{}{}{}{}{}{}}}$}
\newbox\boxtwo\setbox\boxtwo=\hbox{$
   \divide\atflowunit by5\multiply\atflowunit by3\afsetunits
   \atomicflow{(0,0)*{\afacu{}{}{}{}{}{}}}$}
\newbox\boxthree\setbox\boxthree=\hbox{$
   \divide\atflowunit by5\multiply\atflowunit by3\afsetunits
   \atomicflow{(0,0)*{\afaid{}{}{}{}{}{}}}$}
\newbox\boxfour\setbox\boxfour=\hbox{$
   \divide\atflowunit by5\multiply\atflowunit by3\afsetunits
   \atomicflow{(0,0)*{\afaiu{}{}{}{}{}{}}}$}
\newbox\boxfive\setbox\boxfive=\hbox{$
   \divide\atflowunit by5\multiply\atflowunit by3\afsetunits
   \atomicflow{(0,0)*{\afawd{}{}{}{}{}{}}}$}
\newbox\boxsix\setbox\boxsix=\hbox{$
   \divide\atflowunit by5\multiply\atflowunit by3\afsetunits
   \atomicflow{(0,0)*{\afawu{}{}{}{}{}{}}}$}
\[
\atomicflow{
( -22, 12)*{\afvj4};
( -20, 12)*{\cdots};
( -18, 12)*{\afvj4};
(  -7, 12)*{\afvj4};
(  -5, 12)*{\cdots};
(  -3, 12)*{\afvj4};
( -20,  6)*{\copy\boxsix};
( -20,  6)*{\affr88};
(-8.5,  6)*{\copy\boxtwo};
(-1.5,  6)*{\copy\boxone};
(  -5,  6)*{\affr{18}8};
(  10,  6)*{\copy\boxthree};
(  10,  6)*{\affr88};
( -12,  0)*{\afvj4};
( -10,  0)*{\cdots};
(  -8,  0)*{\afvj4};
(  -2,  0)*{\afvj4};
(   0,  0)*{\cdots};
(   2,  0)*{\afvj4};
(   8,  0)*{\afvj4};
(  10,  0)*{\cdots};
(  12,  0)*{\afvj4};
(  20, -6)*{\copy\boxfive};
(  20, -6)*{\affr88};
( 8.5, -6)*{\copy\boxone};
( 1.5, -6)*{\copy\boxtwo};
(   5, -6)*{\affr{18}8};
( -10, -6)*{\copy\boxfour};
( -10, -6)*{\affr88};
(   3, 12)*{\afvj4};
(   5,-12)*{\cdots};
(   7,-12)*{\afvj4};
(  18,-12)*{\afvj4};
(  20,-12)*{\cdots};
(  22,-12)*{\afvj4}}
\quad
\atomicflow{
(-22, 12)*{\afvj4};
(-20, 12)*{\cdots};
(-18, 12)*{\afvj4};
( -7, 12)*{\afvj4};
( -5, 12)*{\cdots};
( -3, 12)*{\afvj4};
(-20,  6)*{\copy\boxsix};
(-20,  6)*{\affr88};
( -5,  6)*{\copy\boxtwo};
( -5,  6)*{\affr{18}8};
( 10,  6)*{\copy\boxthree};
( 10,  6)*{\affr88};
(-12,  0)*{\afvj4};
(-10,  0)*{\cdots};
( -8,  0)*{\afvj4};
( -2,  0)*{\afvj4};
(  0,  0)*{\cdots};
(  2,  0)*{\afvj4};
(  8,  0)*{\afvj4};
( 10,  0)*{\cdots};
( 12,  0)*{\afvj4};
( 20, -6)*{\copy\boxfive};
( 20, -6)*{\affr88};
(  5, -6)*{\copy\boxone};
(  5, -6)*{\affr{18}8};
(-10, -6)*{\copy\boxfour};
(-10, -6)*{\affr88};
(  3,-12)*{\afvj4};
(  5,-12)*{\cdots};
(  7,-12)*{\afvj4};
( 18,-12)*{\afvj4};
( 20,-12)*{\cdots};
( 22,-12)*{\afvj4}}
\quad.
\]
We see on the left the shape of flows of super-streamlined derivations and on the right that of flows of hyper-streamlined derivations, where the boxes represent flows obtained by freely composing edges and vertices whose labels are only those indicated on the boxes. Compare to Remark~\ref{RemStr}. This is the shape of the atomic flow of a hyper-streamlined proof:
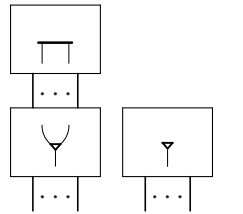
\[
\newbox\boxone\setbox\boxone=\hbox{$
   \divide\atflowunit by5\multiply\atflowunit by3\afsetunits
   \atomicflow{(0,0)*{\afacd{}{}{}{}{}{}}}$}
\newbox\boxthree\setbox\boxthree=\hbox{$
   \divide\atflowunit by5\multiply\atflowunit by3\afsetunits
   \atomicflow{(0,0)*{\afaid{}{}{}{}{}{}}}$}
\newbox\boxfive\setbox\boxfive=\hbox{$
   \divide\atflowunit by5\multiply\atflowunit by3\afsetunits
   \atomicflow{(0,0)*{\afawd{}{}{}{}{}{}}}$}
\atomicflow{ 
( 10,  6)*{\copy\boxthree};
( 10,  6)*{\affr88};
(  8,  0)*{\afvj4};
( 10,  0)*{\cdots};
( 12,  0)*{\afvj4};
( 20, -6)*{\copy\boxfive};
( 20, -6)*{\affr88};
( 10, -6)*{\copy\boxone};
( 10, -6)*{\affr88};
(  8,-12)*{\afvj4};
( 10,-12)*{\cdots};
( 12,-12)*{\afvj4};
( 18,-12)*{\afvj4};
( 20,-12)*{\cdots};
( 22,-12)*{\afvj4}}
\quad.
\]
Note that a hyper-streamlined $\SKS$ proof is a $\KS$ proof.
\end{rem}

%---------------------------------------
\begin{thm}\label{TheoStr}
Algorithm\/ $\Str$ super-streamlines (and so, streamlines) every\/ $\SKS$ derivation.
\end{thm}

%---------------------------------------
\proof
We refer to Definition~\ref{DefStr}. Theorems~\ref{TheoElCycles} and \ref{TheoElPaths} guarantee that there are no $\ai$-connections in $\Phi_4$'s flow $A_4$. By Theorem~\ref{TheoWTerm}, Step 3 terminates. Since no new $\ai$-paths are created by Step 3, and because of Remark~\ref{RemNDCW}, we can conclude that $\Phi_5$ is super-streamlined. Since the reductions of $\frw$ are sound (Theorem~\ref{TheoSound}), the premiss and conclusion of $\Phi_5$ are the same as those of $\Phi_0$.
\qed

Cut elimination immediately follows (see also Remark~\ref{RemStrSKS}).

%---------------------------------------
\begin{cor}[Cut elimination]
Algorithm\/ $\Str$, when applied to an\/ $\SKS$ proof, produces a cut-free proof.
\end{cor}

%---------------------------------------
\begin{rem}
Cut-free proofs obtained from $\Str$ do not exhibit any coweakening rule, but there might be cocontraction rules, since Step 2 of algorithm $\EX$ might introduce some.
\end{rem}

We can easily define a stronger algorithm than $\Str$.

\newcommand{\HStr}{\mathsf{HStr}}
%---------------------------------------
\begin{defi}
For every $\SKS$ derivation $\Phi_0$, \emph{algorithm\/ $\HStr$} (for \emph{hyper-streamlining}) performs the following steps:
\begin{enumerate}
\item\emph{Super-streamline.} Apply $\Str$ to $\Phi_0$ and obtain $\Phi_5$.
\item\emph{Move away contractions and cocontractions.} Transform $\Phi_5$ into algorithm's output $\Phi_6$, whose flow $A_6$ is normal for $\frc$.
\end{enumerate}
\end{defi}

%---------------------------------------
\begin{thm}
Algorithm\/ $\HStr$ hyper-streamlines every\/ $\SKS$ derivation; moreover, when applied to an\/ $\SKS$ proof, $\HStr$ produces a\/ $\KS$ proof.
\end{thm}

%---------------------------------------
\proof
We appeal to Theorem~\ref{TheoStr} and we note that, by Proposition~\ref{PropCW}, the flow of the output derivation is normal for $\frw$ and $\frc$. 
\qed

%---------------------------------------
\begin{rem}
Thanks to algorithm $\HStr$, we can easily obtain a `decomposition' result (as this kind of theorems are called in the deep-inference literature). Any $\SKS$ derivation can be reduced to a derivation of the form
\[
\vlderivation                         {
\vlde{}{\{\aiu,\awd,\acd\}}{\delta}  {
\vlde{}{\{\swi,\med     \}}{\gamma} {
\vlde{}{\{\aid,\awu,\acu\}}{\beta }{
\vlhy                      {\alpha}}}}}
\quad.
\]
In fact, a given derivation is transformed into a hyper-streamlined one, and then it is reduced to the form above by (obvious) permutations of inference steps. Actually, we can obtain more detailed normal forms than the above, for example, any of the following:
\[
\vlderivation                           {
\vlde{}{\{\aiu     \}}{\alpha'_1}      {
\vlde{}{\{\awd     \}}{\alpha'_2}     {
\vlde{}{\{\acd     \}}{\alpha'_3}    {
\vlde{}{\{\swi,\med\}}{\alpha'_4}   {
\vlde{}{\{\acu     \}}{\alpha_4 }  {
\vlde{}{\{\awu     \}}{\alpha_3 } {
\vlde{}{\{\aid     \}}{\alpha_2 }{
\vlhy                 {\alpha_1 }}}}}}}}}
\quad,\qquad
\vlderivation                           {
\vlde{}{\{\awd     \}}{\alpha'_1}      {
\vlde{}{\{\aiu     \}}{\alpha'_2}     {
\vlde{}{\{\acd     \}}{\alpha'_3}    {
\vlde{}{\{\swi,\med\}}{\alpha'_4}   {
\vlde{}{\{\acu     \}}{\alpha_4 }  {
\vlde{}{\{\aid     \}}{\alpha_3 } {
\vlde{}{\{\awu     \}}{\alpha_2 }{
\vlhy                 {\alpha_1 }}}}}}}}}
\quad,\qquad
\vlderivation                           {
\vlde{}{\{\awd     \}}{\alpha'_1}      {
\vlde{}{\{\acd     \}}{\alpha'_2}     {
\vlde{}{\{\aiu     \}}{\alpha'_3}    {
\vlde{}{\{\swi,\med\}}{\alpha'_4}   {
\vlde{}{\{\aid     \}}{\alpha_4 }  {
\vlde{}{\{\acu     \}}{\alpha_3 } {
\vlde{}{\{\awu     \}}{\alpha_2 }{
\vlhy                 {\alpha_1 }}}}}}}}}
\quad\hbox{, etc.}\qquad
\]
In \cite{BrunTiu:01:A-Local-:mz,Brun:04:Deep-Inf:rq}, this result is obtained as a consequence of cut elimination for $\SKS$. We note that direct attempts at proving decomposition can face hard termination problems when dealing with the permutation of $\acd$ over $\aiu$ inference steps (and the dual case). The problem boils down to taming infinite loops of the same nature as those mentioned in Remark~\ref{RemCycle}. The streamlined forms of derivations, of course, take care of this. For an instructive example of this phenomenon, in a more complicated logic to which we hope to extend atomic flows, we refer the reader to \cite{GuglStra:02:A-Non-co:dq}.
\end{rem}

The complexity of $\Str$ and $\HStr$ is exponential on the size of the input derivation, mainly due to the recursive duplication of derivations (see Remark~\ref{RemExpGr}). This is in line with the expectations for propositional logic normalisation. As we said in Remark~\ref{RemExpC}, normalising via $\to_\frc$ also bears a potentially exponential cost, due to the contraction/cocontraction reduction. Normalising via $\to_\frw$, instead, reduces the complexity of derivations, as Figure~\ref{FigExStrRed} shows. These algorithms can be optimised, in the sense of complexity of flows and proofs, in many ways. For example, a simple method for keeping the complexity low is to always reduce weakenings as soon as possible.

There is a source of exponential speed-ups that is worth exploring in the future. The definition of $\to_\frbc$ and $\to_\frex$ is based on recursively copying and stitching together copies of an entire flow. This is so because we need to make sure that, corresponding to the flows that we stitch together, there are proper corresponding derivations. However, it is not necessary to operate on the entire flow, it would be enough to copy and stitch only the minimal flow containing a simple edge and for which a subderivation can be found. In other words, we are duplicating entire flows for the only reason that we know that they correspond to derivations, but a more sophisticated approach could be able to do better.

We know that there is no strongly normalising algorithm based on $\to_\frse$, if by `strongly normalising' we mean a completely liberal use of $\to_\frse$. However, natural constraints on $\to_\frse$ might lead to strong normalisation. For example, it might be possible to succeed in this if we adopt the notion of `minimal derivation containing a simple edge' mentioned above.

%===============================================================================
\section{Conclusions}

We have shown a novel method of control for normalisation algorithms. It is based on a simple, graphical formalism, called `atomic flows', similar both to proof nets and flow graphs. Atomic flows appear to capture an important aspect of the normalisation process. With a novel technique supported by atomic flows, we proved a symmetric generalisation of cut elimination for derivations, which seems to be desperately complicated when subjected to some traditional syntactic analysis.

We and other researchers, in particular Fran\c{c}ois Lamarche and Lutz Stra{\ss}burger, found several times that dealing with loops similar in nature to those of Remark~\ref{RemCycle} is very challenging (see, for example, \cite{GuglStra:02:A-Non-co:dq,LamaStra:05:Construc:qq,LamaStra:05:Naming-P:ov,Stra:06:On-the-A:jy}). The similarities suggest that the techniques of the present paper are probably applicable to other normalisation algorithms and to categorical axiomatisations of classical logic.

We have exploited the possibility, peculiar to deep inference, of designing proof systems whose structural rules are all atomic and whose logical rules are all linear. This allows a great simplification in the geometric study of dependencies between inference rules. In comparison, the sequent calculus, for example, cannot exhibit an atomic contraction rule; this means that contractions on single atom occurrences are not independent, and so, substituting atom occurrences with more complex formulae becomes practically impossible.

Our success stems from the possibility of easily manipulating complex graphs by accessing single, independent occurrences of atoms. We expect to quickly broaden the range of applications of our methods, because all the major logics enjoy presentations in deep inference with atomic structural rules and linear logical ones, contrary to what is possible in any other known formalism. However, our methods so far rely on contraction (and cocontraction for symmetric normalisations), so, ironically, it is not obvious how to use them in the case of pure linear logics.

This work shows that cut elimination is far less a delicate property than it is usually assumed. As a matter of fact, any choice of logical inference rules, provided it makes for an implicationally complete system (so that we can recover the switch and medial rules), would leave the results presented here intact. We argued that the normalisation algorithm itself is not delicate, and enjoys a vast range of possibilities for optimisation. We find it interesting that, by adopting a more liberal syntactic discipline than that of non-deep inference, we correspondingly obtain more freedom for normalisation. Some would expect the opposite to happen. We interpret this as further evidence that traditional proof theory is too syntactic, to the point that syntactic artefacts obscure a deeper and simpler reality. 

Much of this paper has been about the flow rewriting systems $\frw$ and $\frc$, which allow reductions of flows based on local reduction rules. We then developed global reduction mechanisms, based on the global reduction $\to_\frse$, only enough to show our normalisation theorem, with the minimum conceptual effort. However, there is a multitude of interesting global reduction mechanisms, arising from the atomic flow perspective, that are worth exploring. The paper \cite{GuglGund:08:Normalis:yu} is devoted to these global mechanisms, especially in connection with possible computational interpretations.

We are currently investigating, together with Michel Parigot, the use of atomic flows in a computational interpretation of normalisation in $\SKS$. We are also trying to extend our technique to the case of intuitionistic \cite{Tiu:06:A-Local-:gf} and modal logics \cite{GoreTiu:06:Classica:uq,Stou:06:A-Deep-I:rt}. Preliminary investigations with Alwen Tiu lead us to think that the algorithm can be adapted to his presentation of intuitionistic logic in the calculus of structures; in particular, flow polarities match Tiu's polarities, which is an interesting phenomenon that we do not fully understand at this point. We anticipate that modalities will play a crucial role in understanding constructivity beyond intuitionist logic. 

With Michel Parigot, we are using atomic flows to design and test the properties of a new bureaucracy-free formalism that we are temporarily calling \emph{formalism\/ $\mathsf B$}. Just to give an example, in formalism $\mathsf B$, the derivation on the right in Figure~\ref{FigExAF} (page~\pageref{FigExAF}) becomes
\[
\vlgoodsyntax\vlnosmallbrackets
\vls(
\vlderivation                                     {
\vlin{\med}{}{\vls([a.b].[a.b])                 }{
\vlhy        {\vls[\vlinf{\acu}{}{\vls(a.a)}a.
                   \vlinf{\acu}{}{\vls(b.b)}b]}}}
.
\vlinf{\acu}{}{\vls(c.c)}c)
\quad,
\]
where disjunctions and conjunctions are considered modulo commutativity and associativity. Derivations in formalism $\mathsf B$ are, basically, atomic flows enriched with logical information. Consequently, in this formalism, the behaviour of normalisation algorithms controlled by atomic flows is much more natural than the same for the calculus of structures.

In the future, we want to explore the relations between atomic flows and proof nets, as in \cite{LamaStra:05:Construc:qq,LamaStra:05:Naming-P:ov,LamaStra:06:From-Pro:et,Stra:06:On-the-A:jy,StraLama:04:On-Proof:ec}. An interesting problem is to find simple combinatorial conditions that decide whether an atomic flow is associated to some derivation with a given premiss and conclusion. Our results might help in restricting this problem to cases where no structural rules are involved, \emph{i.e.}, to the purely linear fragment consisting of the rules switch and medial. The problem for switch in isolation is solved by the correctness criteria for multiplicative linear logic, see Christian Retor\'e's \cite{Reto:03:Handsome:fj}. Lutz Stra\ss burger, in \cite{Stra:07:A-Charac:fk}, found a criterion for the system containing only medial. For the combination of switch and medial, no criterion is known.

We are also interested in exploring the connections between our work and Craig's interpolation theorem, especially in relation with the size of interpolants (see Alessandra Carbone's work \cite{Carb:97:Interpol:fk} for a possibly related approach to ours). There appear to be connections between atomic flows and Dominic Hughes' combinatorial proofs \cite{Hugh:06:Proofs-W:fk}; we are especially interested in generating combinatorial proofs by manipulating atomic flows. Atomic flows could be given an algebraic characterisation with Albert Burroni's polygraphs \cite{Burr:93:Higher-D:fu}; Yves Guiraud investigated similar constructions in his paper \cite{Guir:06:The-Thre:qt}, where he analyses this way the structural bureaucracy of $\SKS$ derivations. Finally, there might be connections with the theory of matings, as in the works \cite{Andr:81:Theorem-:yf} by Peter Andrews and \cite{Bibe:81:On-Matri:wq} by Wolfgang Bibel; that theory also can be considered an abstract characterisation of proofs.

%===============================================================================
\section*{Acknowledgements}

We thank Michel Parigot for feedback during the early development of atomic flows, and for extensively testing the ideas presented here and suggesting improvements to this paper. We thank Kai Br\"unnler, Paola Bruscoli and Lutz Stra{\ss}burger for comments on a draft and Alwen Tiu for stimulating discussions and for finding a mistake in an earlier version of this work. We are very grateful to the anonymous referees for their splendid work and the many suggestions they had for improving this paper.

\iflmcs\else\let\oldurl\url\renewcommand{\url}[1]{\hfill\break\oldurl{#1}}\fi

\bibliographystyle{alpha}
\bibliography{biblio}

\begin{thebibliography}{Kah07b}

\bibitem[And81]{Andr:81:Theorem-:yf}
Peter~B. Andrews.
\newblock Theorem proving via general matings.
\newblock {\em Journal of the {ACM}}, 28(2):193--214, 1981.

\bibitem[BG08]{BrusGugl:07:On-the-P:fk}
Paola Bruscoli and Alessio Guglielmi.
\newblock On the proof complexity of deep inference.
\newblock {\em ACM Transactions on Computational Logic}, 2008.
\newblock In press. \url{http://cs.bath.ac.uk/ag/p/PrComplDI.pdf}.

\bibitem[Bib81]{Bibe:81:On-Matri:wq}
Wolfgang Bibel.
\newblock On matrices with connections.
\newblock {\em Journal of the {ACM}}, 28(4):633--645, 1981.

\bibitem[Bru02]{Brus:02:A-Purely:wd}
Paola Bruscoli.
\newblock A purely logical account of sequentiality in proof search.
\newblock In Peter~J. Stuckey, editor, {\em Logic Programming, 18th
  International Conference}, volume 2401 of {\em Lecture Notes in Computer
  Science}, pages 302--316. Springer-Verlag, 2002.
\newblock \url{http://cs.bath.ac.uk/pb/bvl/bvl.pdf}.

\bibitem[Br{\"u}03]{Brun:03:Atomic-C:oz}
Kai Br{\"u}nnler.
\newblock Atomic cut elimination for classical logic.
\newblock In M.~Baaz and J.~A. Makowsky, editors, {\em CSL 2003}, volume 2803
  of {\em Lecture Notes in Computer Science}, pages 86--97. Springer-Verlag,
  2003.
\newblock \url{http://www.iam.unibe.ch/~kai/Papers/ace.pdf}.

\bibitem[Br{\"u}04]{Brun:04:Deep-Inf:rq}
Kai Br{\"u}nnler.
\newblock {\em Deep Inference and Symmetry in Classical Proofs}.
\newblock Logos Verlag, Berlin, 2004.
\newblock \url{http://www.iam.unibe.ch/~kai/Papers/phd.pdf}.

\bibitem[Br{\"u}06a]{Brun:06:Cut-Elim:cq}
Kai Br{\"u}nnler.
\newblock Cut elimination inside a deep inference system for classical
  predicate logic.
\newblock {\em Studia Logica}, 82(1):51--71, 2006.
\newblock \url{http://www.iam.unibe.ch/~kai/Papers/q.pdf}.

\bibitem[Br{\"u}06b]{Brun:06:Deep-Inf:qy}
Kai Br{\"u}nnler.
\newblock Deep inference and its normal form of derivations.
\newblock In Arnold Beckmann, Ulrich Berger, Benedikt L{\"o}we, and John~V.
  Tucker, editors, {\em Computability in Europe 2006}, volume 3988 of {\em
  Lecture Notes in Computer Science}, pages 65--74. Springer-Verlag, July 2006.
\newblock \url{http://www.iam.unibe.ch/~kai/Papers/n.pdf}.

\bibitem[Br{\"u}06c]{Brun::Deep-Seq:ay}
Kai Br{\"u}nnler.
\newblock Deep sequent systems for modal logic.
\newblock In Guido Governatori, Ian Hodkinson, and Yde Venema, editors, {\em
  Advances in Modal Logic}, volume~6, pages 107--119. College Publications,
  2006.
\newblock \url{http://www.aiml.net/volumes/volume6/Bruennler.ps}.

\bibitem[Br{\"u}06d]{Brun:06:Locality:zh}
Kai Br{\"u}nnler.
\newblock Locality for classical logic.
\newblock {\em Notre Dame Journal of Formal Logic}, 47(4):557--580, 2006.
\newblock \url{http://www.iam.unibe.ch/~kai/Papers/LocalityClassical.pdf}.

\bibitem[BT01]{BrunTiu:01:A-Local-:mz}
Kai Br{\"u}nnler and Alwen~Fernanto Tiu.
\newblock A local system for classical logic.
\newblock In R.~Nieuwenhuis and A.~Voronkov, editors, {\em LPAR 2001}, volume
  2250 of {\em Lecture Notes in Artificial Intelligence}, pages 347--361.
  Springer-Verlag, 2001.
\newblock \url{http://www.iam.unibe.ch/~kai/Papers/lcl-lpar.pdf}.

\bibitem[Bur93]{Burr:93:Higher-D:fu}
Albert Burroni.
\newblock Higher dimensional word problems with applications to equational
  logic.
\newblock {\em Theoretical Computer Science}, 115(1):43--62, (1993).

\bibitem[Bus91]{Buss:91:The-Unde:uq}
Samuel~R. Buss.
\newblock The undecidability of k-provability.
\newblock {\em Annals of Pure and Applied Logic}, 53(1):75--102, 1991.

\bibitem[Car97]{Carb:97:Interpol:fk}
Alessandra Carbone.
\newblock Interpolants, cut elimination and flow graphs for the propositional
  calculus.
\newblock {\em Annals of Pure and Applied Logic}, 83:249--299, 1997.

\bibitem[CR79]{CookReck:79:The-Rela:mf}
Stephen~A. Cook and Robert~A. Reckhow.
\newblock The relative efficiency of propositional proof systems.
\newblock {\em Journal of Symbolic Logic}, 44(1):36--50, 1979.

\bibitem[DG04]{Di-G:04:Structur:wy}
Pietro Di~Gianantonio.
\newblock Structures for multiplicative cyclic linear logic: {D}eepness vs
  cyclicity.
\newblock In J.~Marcinkowski and A.~Tarlecki, editors, {\em CSL 2004}, volume
  3210 of {\em Lecture Notes in Computer Science}, pages 130--144.
  Springer-Verlag, 2004.
\newblock
  \url{http://www.dimi.uniud.it/~pietro/papers/Soft-copy-ps/scll.ps.gz}.

\bibitem[Gen69]{Gent:69:Investig:xi}
Gerhard Gentzen.
\newblock Investigations into logical deduction.
\newblock In M.~E. Szabo, editor, {\em The Collected Papers of Gerhard
  Gentzen}, pages 68--131. North-Holland, Amsterdam, 1969.

\bibitem[GG08]{GuglGund:08:Normalis:yu}
Alessio Guglielmi and Tom Gundersen.
\newblock Normalisation control in deep inference via atomic flows {II}.
\newblock \url{http://cs.bath.ac.uk/ag/p/NormContrDIAtFl2.pdf}, 2008.

\bibitem[Gir87]{Gira:87:Linear-L:wm}
Jean-Yves Girard.
\newblock Linear logic.
\newblock {\em Theoretical Computer Science}, 50:1--102, 1987.

\bibitem[GS01]{GuglStra:01:Non-comm:rp}
Alessio Guglielmi and Lutz Stra{\ss}burger.
\newblock Non-commutativity and {MELL} in the calculus of structures.
\newblock In L.~Fribourg, editor, {\em CSL 2001}, volume 2142 of {\em Lecture
  Notes in Computer Science}, pages 54--68. Springer-Verlag, September 2001.
\newblock \url{http://cs.bath.ac.uk/ag/p/NoncMELLCoS.pdf}.

\bibitem[GS02]{GuglStra:02:A-Non-co:lq}
Alessio Guglielmi and Lutz Stra{\ss}burger.
\newblock A non-commutative extension of {MELL}.
\newblock In M.~Baaz and A.~Voronkov, editors, {\em LPAR 2002}, volume 2514 of
  {\em Lecture Notes in Artificial Intelligence}, pages 231--246.
  Springer-Verlag, October 2002.
\newblock \url{http://www.lix.polytechnique.fr/~lutz/papers/NEL.pdf}.

\bibitem[GS07]{GuglStra:02:A-Non-co:dq}
Alessio Guglielmi and Lutz Stra{\ss}burger.
\newblock A system of interaction and structure {IV}: {T}he exponentials.
\newblock In the second round of revision for Mathematical Structures in
  Computer Science.
  \url{http://www.lix.polytechnique.fr/~lutz/papers/NELbig.pdf}, 2007.

\bibitem[GT07]{GoreTiu:06:Classica:uq}
Rajeev Gor{\'e} and Alwen Tiu.
\newblock Classical modal display logic in the calculus of structures and
  minimal cut-free deep inference calculi for {S5}.
\newblock {\em Journal of Logic and Computation}, 17(4):767--794, 2007.
\newblock \url{http://users.rsise.anu.edu.au/~tiu/papers/cmdl.pdf}.

\bibitem[Gug07]{Gugl:06:A-System:kl}
Alessio Guglielmi.
\newblock A system of interaction and structure.
\newblock {\em ACM Transactions on Computational Logic}, 8(1):1--64, 2007.
\newblock \url{http://cs.bath.ac.uk/ag/p/SystIntStr.pdf}.

\bibitem[Gui06]{Guir:06:The-Thre:qt}
Yves Guiraud.
\newblock The three dimensions of proofs.
\newblock {\em Annals of Pure and Applied Logic}, 141(1-2):266--295, 2006.
\newblock \url{http://www.loria.fr/~guiraudy/recherche/cos1.pdf}.

\bibitem[Hug06]{Hugh:06:Proofs-W:fk}
Dominic~J.D. Hughes.
\newblock Proofs without syntax.
\newblock {\em Annals of Mathematics}, 164(3):1065--1076, 2006.

\bibitem[Kah06]{Kahr:06:Reducing:hc}
Ozan Kahramano{\u g}ullar{\i}.
\newblock Reducing nondeterminism in the calculus of structures.
\newblock In M.~Hermann and A.~Voronkov, editors, {\em LPAR 2006}, volume 4246
  of {\em Lecture Notes in Artificial Intelligence}, pages 272--286.
  Springer-Verlag, 2006.
\newblock \url{http://www.doc.ic.ac.uk/~ozank/Papers/reducingNondet.pdf}.

\bibitem[Kah07a]{Kahr:07:Maude-as:lr}
Ozan Kahramano{\u g}ullar{\i}.
\newblock Maude as a platform for designing and implementing deep inference
  systems.
\newblock In {\em RULE 2007---The Eighth International Workshop on Rule-Based
  Programming}, Electronic Notes in Theoretical Computer Science. Elsevier,
  2007.
\newblock In press. \url{http://www.doc.ic.ac.uk/~ozank/Papers/rule07.pdf}.

\bibitem[Kah07b]{Kahr:07:System-B:fk}
Ozan Kahramano{\u g}ullar{\i}.
\newblock System {BV} is {NP}-complete.
\newblock {\em Annals of Pure and Applied Logic}, 2007.
\newblock In press.
  \url{http://www.doc.ic.ac.uk/~ozank/Papers/bv_npc_apal.pdf}.

\bibitem[Laf97]{Lafo:97:Interact:jb}
Yves Lafont.
\newblock Interaction combinators.
\newblock {\em Information and Computation}, 137:69--101, 1997.

\bibitem[LS05a]{LamaStra:05:Construc:qq}
Fran{\c c}ois Lamarche and Lutz Stra{\ss}burger.
\newblock Constructing free boolean categories.
\newblock In Prakash Panangaden, editor, {\em 20th Annual IEEE Symposium on
  Logic in Computer Science}, pages 209--218. IEEE, 2005.
\newblock \url{http://www.lix.polytechnique.fr/~lutz/papers/FreeBool-long.pdf}.

\bibitem[LS05b]{LamaStra:05:Naming-P:ov}
Fran{\c c}ois Lamarche and Lutz Stra{\ss}burger.
\newblock Naming proofs in classical propositional logic.
\newblock In Pawe{\l} Urzyczyn, editor, {\em Typed Lambda Calculi and
  Applications}, volume 3461 of {\em Lecture Notes in Computer Science}, pages
  246--261. Springer-Verlag, 2005.
\newblock
  \url{http://www.lix.polytechnique.fr/~lutz/papers/namingproofsCL.pdf}.

\bibitem[LS06]{LamaStra:06:From-Pro:et}
Fran{\c c}ois Lamarche and Lutz Stra{\ss}burger.
\newblock From proof nets to the free *-autonomous category.
\newblock {\em Logical Methods in Computer Science}, 2(4):3:1--44, 2006.
\newblock \url{http://arxiv.org/pdf/cs.LO/0605054}.

\bibitem[Ret03]{Reto:03:Handsome:fj}
Christian Retor{\'e}.
\newblock Handsome proof-nets: {P}erfect matchings and cographs.
\newblock {\em Theoretical Computer Science}, 294(3):473--488, 2003.

\bibitem[SL04]{StraLama:04:On-Proof:ec}
Lutz Stra{\ss}burger and Fran{\c c}ois Lamarche.
\newblock On proof nets for multiplicative linear logic with units.
\newblock In J.~Marcinkowski and A.~Tarlecki, editors, {\em CSL 2004}, volume
  3210 of {\em Lecture Notes in Computer Science}, pages 145--159.
  Springer-Verlag, 2004.
\newblock \url{http://www.lix.polytechnique.fr/~lutz/papers/multPN.pdf}.

\bibitem[Sto07]{Stou:06:A-Deep-I:rt}
Phiniki Stouppa.
\newblock A deep inference system for the modal logic {S5}.
\newblock {\em Studia Logica}, 85(2):199--214, 2007.
\newblock
  \url{http://www.iam.unibe.ch/til/publications/pubitems/pdfs/sto07.pdf}.

\bibitem[Str02]{Stra:02:A-Local-:ul}
Lutz Stra{\ss}burger.
\newblock A local system for linear logic.
\newblock In M.~Baaz and A.~Voronkov, editors, {\em LPAR 2002}, volume 2514 of
  {\em Lecture Notes in Artificial Intelligence}, pages 388--402.
  Springer-Verlag, 2002.
\newblock \url{http://www.lix.polytechnique.fr/~lutz/papers/lls-lpar.pdf}.

\bibitem[Str03a]{Stra:03:Linear-L:lp}
Lutz Stra{\ss}burger.
\newblock {\em Linear Logic and Noncommutativity in the Calculus of
  Structures}.
\newblock PhD thesis, Technische Universit{\"a}t Dresden, 2003.
\newblock \url{http://www.lix.polytechnique.fr/~lutz/papers/dissvonlutz.pdf}.

\bibitem[Str03b]{Stra:03:MELL-in-:oy}
Lutz Stra{\ss}burger.
\newblock {MELL} in the calculus of structures.
\newblock {\em Theoretical Computer Science}, 309:213--285, 2003.
\newblock \url{http://www.lix.polytechnique.fr/~lutz/papers/els.pdf}.

\bibitem[Str07a]{Stra:07:A-Charac:fk}
Lutz Stra{\ss}burger.
\newblock A characterisation of medial as rewriting rule.
\newblock In Franz Baader, editor, {\em RTA 2007}, volume 4533 of {\em Lecture
  Notes in Computer Science}, pages 344--358. Springer-Verlag, 2007.
\newblock \url{http://www.lix.polytechnique.fr/~lutz/papers/CharMedial.pdf}.

\bibitem[Str07b]{Stra:06:On-the-A:jy}
Lutz Stra{\ss}burger.
\newblock On the axiomatisation of boolean categories with and without medial.
\newblock {\em Theory and Applications of Categories}, 18(18):536--601, 2007.
\newblock \url{http://www.lix.polytechnique.fr/~lutz/papers/medial.pdf}.

\bibitem[Tiu06a]{Tiu:06:A-Local-:gf}
Alwen Tiu.
\newblock A local system for intuitionistic logic.
\newblock In M.~Hermann and A.~Voronkov, editors, {\em LPAR 2006}, volume 4246
  of {\em Lecture Notes in Artificial Intelligence}, pages 242--256.
  Springer-Verlag, 2006.
\newblock \url{http://users.rsise.anu.edu.au/~tiu/localint.pdf}.

\bibitem[Tiu06b]{Tiu:06:A-System:ai}
Alwen Tiu.
\newblock A system of interaction and structure {II}: {T}he need for deep
  inference.
\newblock {\em Logical Methods in Computer Science}, 2(2):4:1--24, 2006.
\newblock \url{http://arxiv.org/pdf/cs.LO/0512036}.

\end{thebibliography}

\vskip-50 pt

\end{document}